\documentclass{elsarticle}

\usepackage{lineno,hyperref}
\usepackage{graphicx}
\usepackage{{ntheorem}}
\usepackage[T1]{fontenc}
\usepackage{appendix}
\usepackage{amssymb}
\usepackage[dvipsnames]{xcolor}
\usepackage{amsfonts}
\usepackage{epstopdf}
\usepackage{multicol}
\usepackage{anyfontsize}
\usepackage[utf8]{inputenc}
\usepackage{calligra}
\usepackage{xcolor}
\usepackage{pdfpages}
\usepackage{nicefrac}
\usepackage{subfig}
\usepackage{amsmath}
\usepackage{natbib}
\usepackage{esint}
\usepackage{algorithm}
\usepackage{algorithmic}
\usepackage{comment}
\usepackage{mathtools}

\modulolinenumbers[5]

\journal{Journal of \LaTeX\ Templates}
\def\bcdot{{*}}

\newcommand{\ha}{\frac{1}{2}}

\newcommand{\bi}{\mathbf{i}}
\newcommand{\bj}{\mathbf{j}}

\newcommand{\bhalf}{\mathbf{1/2}}
\newcommand{\bef}{\mathbf{e}_1}
\newcommand{\bes}{\mathbf{e}_2}

\newcommand{\bone}{\mathbf{1}}
\newcommand{\R}{\mathbb{R}}
\newcommand{\A}{\mathcal{A}}
\newcommand{\F}{\mathcal{F}}
\newcommand{\I}{\mathcal{I}}

\newcounter{dctr}[section]
\numberwithin{equation}{section}
\numberwithin{figure}{section}

\newtheorem{theorem}{Theorem}[section]
\newtheorem{remark}{Remark}[section]

\begin{document}
	
	\begin{frontmatter}
		
		\title{\textbf{\large Well-balanced adaptive compact approximate Taylor methods\\ for systems of balance laws}}
		
		\author{
        H. Carrillo$^{\hspace{0.5mm}\textrm{a}}$,
        E. Macca$^{\hspace{0.5mm}\textrm{b}}$,
        Carlos Parés$^{\hspace{0.5mm}\textrm{c}}$,
        G. Russo$^{\hspace{0.5mm}\textrm{d}}$,
        }
        
         \address{$^\textrm{a}$  \emph{hugo.carrillo@rcarbonifera.tecm.mx}, TecNm region carbonifera,Coahuila México. \\
         $^\textrm{b}$ \emph{emanuele.macca@phd.unict.it}, Department of Mathematics and Computer Scince, University of Catania, Catania Italy.\\
         $^\textrm{c}$ \emph{pares@uma.es}, Department of Applied Mathematics, University of Malaga, Malaga, Spain. \\
         $^\textrm{d}$ \emph{russo@dmi.unict.it}, Department of Mathematics and Computer Scince, University of Catania, Catania Italy.
         \vspace{0pt}
         }
    	
    	\begin{abstract}
    	    Compact Approximate Taylor (CAT) methods for systems of conservation laws were introduced by Carrillo and Par\'es in 2019. These methods,  based on a strategy that allows one to extend high-order Lax-Wendroff methods to nonlinear systems without using the Cauchy-Kovalevskaya procedure, have arbitrary even order of accuracy $2p$ and use $(2p +1)$-point stencils, where $p$ is an arbitrary positive integer. More recently in 2021 Carrillo, Macca, Par\'es, Russo and Zor\'{\i}o introduced a strategy to get rid of the spurious oscillations close to discontinuities produced by CAT methods. This strategy led to the so-called Adaptive CAT (ACAT) methods, in which  the order of accuracy -- and thus the width of the stencils -- is adapted to the local smoothness of the solution. The goal of this paper is to extend CAT and ACAT methods to systems of balance laws. To do this, the source term is written as the derivative of its indefinite integral that is formally treated as a flux function. The well-balanced property of the methods is discussed and a variant that allows in principle to preserve any stationary solution is presented. The resulting methods are then applied to a number of systems going from a linear scalar conservation law to the 2D Euler equations with gravity, passing by the Burgers equations with source term and the 1D shallow water equations: the order and well-balanced properties are checked in several numerical tests. 
    	\end{abstract}
	\end{frontmatter}
	
	
	\section{Introduction.}
	This paper deals with the design of high-order well-balanced methods for hyperbolic quasi-linear systems of balance laws
     \begin{equation}
        \label{bal_sis}
        U_t + F(U)_x = S(U)H_x,
    \end{equation}   
 	with initial condition $U(x,0) = U_0(x)$, where $U:\R\times[0,+\infty)\rightarrow \R^d$ is the unknown vector field;  $F:\R^d \rightarrow \R^d$ is the flux function; $S:\mathbb{R}^d\rightarrow\mathbb{R}^d$ is the source term; and $H: \R \rightarrow \R$ is a known function.  	PDE systems of this form appear in many fluid models in different contexts: shallow water models, multiphase flow models, gas dynamic, elastic wave equations, etc.

More precisely, we focus on the extension of high-order Lax-Wendroff methods to systems \eqref{bal_sis}. For linear systems of conservation laws, these methods are based on Taylor expansions in time in which the time derivatives are transformed into spatial derivatives using the governing equations  \cite{LeVeque2007book,Toro2009book,Gideon1971}. The discretization of the spatial derivatives by means of centered high-order differentiation formulas leads then to numerical methods with arbitrary order of accuracy.

The main difficulty to extend Lax-Wendroff methods to nonlinear problems comes from the transformation of time derivatives into spatial derivatives: the use of the equations through the Cauchy-Kovalesky (CK) procedure may be impractical from the computational point of view because it often requires extended symbolic calculus, ended up with inefficient codes. In the context of ADER methods introduced by Toro and collaborators (see \cite{Ader2001, TitarevToro2002, Schwartzkopff2002}), this difficulty has been circumvented   by replacing the CK procedure by local space-time problems that are solved with a Galerkin method: see \cite{DumbserToro2008}, \cite{PNPM}. 

For systems of conservation laws
	\begin{equation}
	    \label{con_sys}
	    U_t + F(U)_x =0,
	\end{equation}
the CK procedure was avoided in \cite{ZBM2017}  by computing time derivatives in a recursive way using high-order centered differentiation formulas combined with  Taylor expansions in time. Nevertheless, if $(2p +1)$-point differentiation formulas are used to compute spatial and temporal derivatives, the resulting method use  $(4p +1)$-point stencils while  Lax-Wendroff methods for linear systems use $(2p +1)$-point ones. 
In \cite{CP2019} a variant of these methods that use $(2p +1)$-point stencils, the so-called Compact Approximated Taylor methods (CAT), was introduced. CAT methods were shown to reduce to the standard high-order Lax-Wendroff methods when applied to linear problems. The technique used to reduce the length of the stencils increases the computational cost of a time step compared to the methods introduced in \cite{ZBM2017}: the Taylor expansions are computed locally, so that the total number of  expansions needed to update the numerical solution is multiplied by $(2p + 1)$. Nevertheless, CAT methods have better stability properties allowing larger time steps, thus compensating the extra cost per time iteration:
see \cite{CP2019}. 

As it happens with high-order Lax-Wendroff methods for linear systems, their extensions to nonlinear problems produce spurious oscillations close to discontinuities and a strategy is needed  to get rid of  them. In \cite{ZBM2017} they were combined with WENO reconstructions to compute the first time derivatives. CAT methods were also combined with WENO in \cite{CP2019} and \cite{CPZ2020} to avoid oscillations near discontinuities. Nevertheless this combination was not optimal: while the best CAT methods are those of even order, WENO methods have odd accuracy order. Moreover, the restriction on the time step imposed by WENO methods may spoil the advantages of the better stability property of CAT methods. To avoid this, in \cite{CPZMR2020} a new version of CAT methods, the so-called Adaptive CAT (ACAT) methods, was introduced  in which the oscillations near discontinuities are cured by adapting the order of accuracy -- and thus the width of the stencils -- to the smoothness of the solution. To do this,  a class of smoothness indicators was introduced. 

The main goal of this article is to extend ACAT methods to balance laws \eqref{bal_sis} and to its 2D counterpart. To do this, we follow the strategy in \cite{Gascon} (see also  \cite{Donat10}) that consists in writing the source term as the derivative of the indefinite integral
$$
       \int_{-\infty}^{x}S(U(\sigma,t))H_x(\sigma)\,d\sigma,
$$
that is then formally treated as a new flux function. Please note that, although this technique is applied here to derive ACAT methods for \eqref{bal_sis}, it can be applied in principle to other extensions of high-order Lax-Wendroff methods to nonlinear  problems using any technique to remove the spurious oscillations.

Systems of balance  laws \eqref{bal_sis} have non-trivial stationary solutions that satisfy the ODE system
  \begin{equation}
        F(U)_x = S(U)H_x.
    \end{equation}   
The objective of {\em well balanced schemes\/} is to preserve exactly or with enhanced accuracy some of these steady state solutions.   In the context of shallow water
equations, Berm\'udez and V\'azquez-Cend\'on introduced in
\cite{BV94} the condition called {\it C-property}: a
scheme is said to satisfy this condition if it preserves the water at rest solutions. 
Since then, many different numerical methods that satisfy this property have been introduced in the literature: see \cite{Bouchut04}, \cite{Xing17} and their references. 
In the framework of finite difference methods,  high-order schemes that satisfy the C-property were introduced in 
 \cite{Caselles09} and  \cite{XingShu1}: while the former was based on the formal writing of the system in conservative form based on the above mentioned technique, the latter relied on the expression of the source term as a function of variables that are constants for the stationary solutions to be preserved: see \cite{XingShu2}.
 A similar technique based on the reconstruction of equilibrium variables has been proposed in 
 \cite{RussoKhe10}.
 In \cite{PP2021} a general technique to derive high-order well-balanced finite-difference methods for systems of balance-laws was introduced. The strategy, inspired on the general technique for finite volume methods discussed in \cite{CP2020}, was as follows:  let $U_i$ be the numerical approximation of the solution $U(x_i, t)$ at the node
 $x_i$ at time $t$ and let  $U^*_i$ be the stationary solution satisfying the Cauchy problem:
 \begin{equation}\label{mot1}
 \begin{cases}
      & \displaystyle F(U^*_i)_x = S(U^*_i)H_x, \\[0.5em]
      & \displaystyle U^*_i(x_i) = U_i.
 \end{cases}
 \end{equation}
 Then, if $U^*_i$ can be found, one has trivially
 \begin{equation}\label{mot2}
 S(U_i) H_x (x_i)  = S(U^*_i(x_i)) H_x (x_i) = F(U^*_i(x_i))_x.
 \end{equation}
 Therefore, locally the system of balance laws can be written in conservation form as follows
 $$
 U_t + (F(U)- F(U^*_i(x)))_x = 0.
 $$
 The numerical method is obtained then by discretizating this conservative form by means of high-order WENO reconstruction of this extended 'flux' function. We will follow this strategy here to derive well-balanced ACAT methods.

This paper is organized as follows: in Section 2, CAT methods for systems of conservation laws \eqref{con_sys} are recalled. In Section 3, these methods are  extended to systems of balance laws \eqref{bal_sis}: after obtaining a first high-order version of the methods for systems of balance laws, the well-balanced property is discussed and a second version is introduced that preserves in principle any stationary solution. In Section  4 the technique introduced in \cite{CPZMR2020}
based on the order-adaption of the  methods is recalled and it is applied to the  2 versions of CAT methods derived in Section 3, which leads to  ACAT and well-balanced ACAT methods for systems of balance laws. Section 5 is devoted to the extension of the methods  to 2D problems. In Section 6 the methods are applied to a number of problems: a 1D scalar balance law, Burgers equation with source term, the 1D shallow water model, and 2D Euler equations with gravitational potential.  Finally, in the last section, we draw some conclusions.

	\section{Compact Approximate Taylor Methods for systems of conservation laws}\label{S:cl}
	First of all, let us recap the expression of CAT methods for a 1D system of conservation laws \eqref{con_sys} with initial condition $U(x,0) = U_0(x)$.	In the CAT$2P$ method, a Taylor expansion in time is used to update the numerical solution: 
	\begin{equation}
        \label{L-W-scheme}
        U_i^{n+1} = U_i^n + \sum_{k=1}^{2P}\frac{(\Delta t)^k}{k!}U^{(k)}_i,
    \end{equation}
    where $\{x_i\}$ are the nodes of a uniform mesh of step $\Delta x;$ $U_i^n$ is an approximation of the value  of the exact solution $U(x_i,t_n)$ at time $t_n=n\Delta t$ at $x_i$; and $U_i^{(k)}$ is an approximation of $\partial^k_tU(x_i,t_n)$.  Following the strategy introduced in \cite{ZBM2017} to derive Approximate Taylor methods, the  Cauchy-Kovalevskaya procedure is  avoided using  the equality 
    \begin{equation}
        \label{CK}
        \partial_t^k U = -\partial_x\partial_t^{k-1}F(U),
    \end{equation}
    satisfied by smooth enough solutions.
    
In the Approximate Taylor methods introduced in \cite{ZBM2017}, global approximations that have to be computed only once at every node were used. Although the use of local approximations may increase the number of calculations, it allows us to prevent the increase of the stencil observed in the methods introduced in this reference: indeed, CAT$2P$ methods are written in conservative form 
    \begin{equation}
        \label{CAT2P}
        U_i^{n+1} = U_i^n + \frac{\Delta t}{\Delta x}\Bigl(F_{i-\ha}^{P} - F_{i+\ha}^{P}\Bigr), 
    \end{equation}
    and only the values
    $$  U_{i-P+1}^n, \dots, U_{i+P}^n  $$
 at the stencil $\mathcal{S}^P_{i+1/2} $ $$
\mathcal{S}^P_{i+1/2} = \{ x_{i-P+1}, \dots, x_{i+P} \},
$$  are used to compute the numerical flux $ F_{i+\ha}^{P}$, so that the stencil used to update the solution at the $i$th node is $\{ x_{i-P}, \dots, x_{i+P}\}$, while in \cite{ZBM2017} the stencil for the same order of accuracy required $4P$ points. More precisely, the time derivatives of the solution will be approximated by applying a formula of numerical differentiation for first order spatial derivatives to some approximations
    \begin{equation}\label{f(k-1)_ij}
{F}^{(k-1)}_{i,j}  \approx \partial_t^{k-1}F(U)(x_{i+j}, t_n), \quad j=-P+1,\dots, P
\end{equation}
that will be computed using recursively Taylor expansions in time. In  this notation, $P$ is an arbitrary positive integer and  $j$ is a local coordinate in the stencil $\mathcal{S}^P_{i+1/2} \equiv  \{-P+1, \dots, P\} $
i.e. $F^{(k)}_{i,j} $ is the approximation of $\partial_t^{k}F(U)$ at time $t_n$ at the  node of local coordinate $j$ of the stencil $\mathcal{S}^P_{i+1/2}$ which is
$x_{i+j}$. These approximations are local in the following sense: let us suppose that $i_1 + j_1 = i_2 + j_2 = l$, 
i.e. $x_l$ belongs to  $\mathcal{S}^P_{i_1+1/2}$ and $\mathcal{S}^P_{i_2+1/2}$ with local coordinates $j_1$ and $j_2$ respectively. Then $F^{(k)}_{i_1,j_1} $ and $F^{(k)}_{i_2,j_2} $ are, in general, two \textit{different} approximations of  $\partial_t^{k}F(U)(x_l, t_n)$.
Moreover, in \cite{CP2019} it has been shown that CAT$2P$  methods reduce to the $2P$-order Lax-Wendroff method for linear systems. This implies the linear stability for these methods under the usual CFL requirement (which is not the case for the methods introduced in \cite{ZBM2017}).

    
 Since numerical differentiation plays a fundamental role in the algorithm, before giving the expression of the numerical flux, let us introduce the notation to describe the formulas that will be used. Given two positive integers $p$, $k$, an index $i$, and a real number $q$, we will represent by
\begin{equation}\label{upwF}
f^{(k)}(x_i + q \Delta x) \approx A^{k,q}_{p}(f_i, \Delta x) = \frac{1}{\Delta x^k} \sum_{j = -p + 1}^p \gamma^{k,q}_{p,j} f(x_{i+j}), \quad \quad k = 0,\ldots,2p-1,
\end{equation}
the interpolatory formula that approximates the $k$-th derivative of a function $f$  at the point $x_i + q \Delta x$ using its values at the $2p$  points $x_{i-p+1}, \dots, x_{i+p}$. For $k = 0$,
 $$
f(x_i+q\Delta x,t_n) = A^{0,q}_{p}(f_i, \Delta x) = \sum_{j = -p + 1}^p \gamma^{0,q}_{p,j} f(x_{i+j}),
$$
represents the value at $x_i + q \Delta x $ of the Lagrange polynomial that interpolates the values of $f$ at the points $x_{i-p+1}, \dots, x_{i+p}$. 

\begin{remark}
The coefficients $ \gamma^{k,q}_{p,j}$  of the differentiation formulas can be recursively computed using the algorithm introduced in \cite{Fornberg1}. See also \cite{CP2019}. 
\end{remark}

The following  notation 
\begin{align}
\label{2.6}
\partial_x^kf(x_i+q\Delta x,t_n) \approx A^{k,q}_{p}(f_{i,\bcdot}^n, \Delta x)   &=   \frac{1}{\Delta x^k} \sum_{j=-p+1}^{p} \gamma^{k,q}_{p,j} f_{i+j}, \\
\label{2.7}
\partial_t^kf(x_i,t_n) \approx A^{k,0}_{p}(f^{n,\bcdot}_i, \Delta t)   &=   \frac{1}{\Delta t^k} \sum_{j=-p+1}^{p} \gamma^{k,q}_{p,j} f_{i+j},
\end{align}
will be used to indicate that the formula is applied to some approximations $f_i$ of $f$ and not to its exact point values $f(x_i)$. 
In cases where there are two or more indices, the symbol $\bcdot$ will be used to indicate to which index (space or time) the differentiation is applied. From now on, since all the formulas are computed at time $t=t_n$ we avoid the extra index $n$ on equations \eqref{2.6} and \eqref{2.7}. For instance, the following approximations will be used in the algorithm adopted to compute the numerical fluxes:

\begin{eqnarray}
& &  \partial^k_t U(x_{i+j}, t_n)  \approx 
-A_{P}^{1,j}\Bigl(F_{i, \bcdot}^{(k-1)},\Delta x\Bigr)
=  \frac{1}{\Delta x} \sum_{l=-P+1}^{P} \gamma^{1,j}_{P,l} F_{i, l}^{(k)}, \label{ex1}\\
& &  \partial^k_t F(U)(x_{i+j}, t_n)  \approx  A_{P}^{k,0}\Bigl(F_{i,j}^{k, \bcdot},\Delta t\Bigr)
 = \frac{1}{\Delta t^k} \sum_{r=-P+1}^{P} \gamma^{k,0}_{P,r} F_{i,j}^{k, n + r},\label{ex2}\\
 \label{ex3}
& & \partial^k_t F(U)\left(x_i + \frac{\Delta x}{2}, t_n \right)  \approx    A^{0, 1/2}_{P}\left({F}_{i,\bcdot}^{(k)}, \Delta x \right)
= \sum_{j=-P+1}^{P} \gamma^{0,1/2}_{P,j} {F}_{i,j}^{(k)}.
 \end{eqnarray}
In \eqref{ex1}, that is the discrete version of \eqref{CK}, numerical differentiation in space is used to approximate the time derivative of the solution at $x_{i+j}$
from the local approximations $F_{i, l}^{(k)}$, $l = -P +1, \dots, P$. 
In \eqref{ex2}, numerical differentiation in time is used to approximate the $k$-th time derivative of $F(U)$ at $x_{i+j}$
from some approximations $F_{i,j}^{k, n+r}$ of $F(U)(x_{i+j}, t_{n+r})$, $r = -P+1, \dots, P$. Finally, in \eqref{ex3},
Lagrange interpolation is used to approximate the value of the $k$-th time derivative of $F(U)$ at $x_i + \Delta x/2$ at 
time $t_n$ from $F_{i, l}^{(k)}$, $l = -P +1, \dots, P$. 
 

Using this notation, the expression of the numerical flux is as follows:
\begin{equation}
    \label{Fp}
    F_{i+\ha}^{P}= \sum_{k=1}^{2P} \frac{(\Delta t)^{k-1}}{k!}A^{0, 1/2}_{P}({F}_{i,\bcdot}^{(k-1)}, \Delta x)
    = \sum_{k=1}^{2P}\frac{(\Delta t)^{k-1}}{k!}\sum_{j=-P+1}^{P}\gamma_{P,j}^{0,\ha}F_{i,j}^{(k-1)},
\end{equation}
where the time derivatives of the flux are computed by the following iterative algorithm (see \cite{CP2019} for details):

  
      \begin{itemize}
        \item Define $F_{i,j}^{(0)}:= F(U_{i+j}^n)$ for all $j = -P+1,\ldots,P.$
        \item For $k=1,\ldots,2P-1:$
            \begin{itemize}
                \item Compute for all $j=-P+1,\ldots,P$ 
                $$ U_{i,j}^{(k)} = -A_{P}^{1,j}\Bigl(F_{i, \bcdot}^{(k-1)},\Delta x\Bigr). $$
                \item Define for all $j,r=-P+1,\ldots,P$
                $$ F_{i,j}^{k,n+r} := F\Bigl(U_{i,j}^{k,n+r}\Bigr), $$ where $U_{i,j}^{k,n+r}$ is the 
                approximation of $U(x_{i+j}, t_{n+r})$ given by the approximate Taylor expansion in time:
                $$ U_{i,j}^{k,n+r} = U_{i+j}^n + \sum_{m = 1}^{k}\frac{(\Delta t)^m}{m!}U_{i,j}^{(m)}. $$
                \item Compute for all $j=-P+1,\ldots,P$
                $$ F_{i,j}^{(k)} = A_{P}^{k,0}\Bigl(F_{i,j}^{k, \bcdot},\Delta t\Bigr). $$
            \end{itemize}
    \end{itemize}

\begin{figure}[!ht]
    	\centering
    	\includegraphics[scale=0.9]{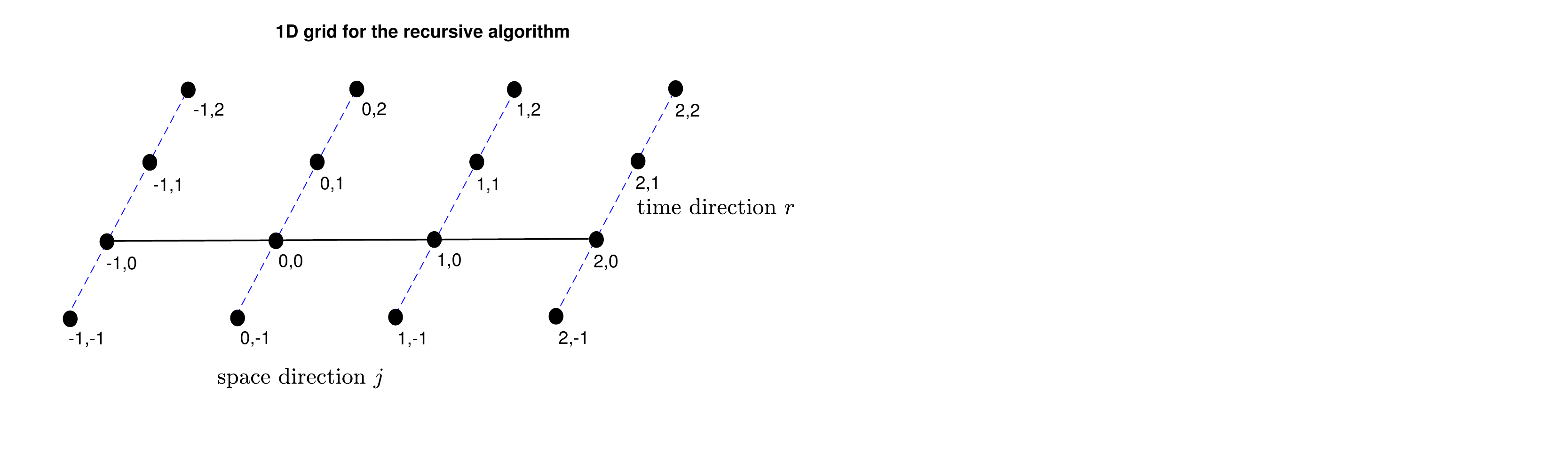}
    	\vspace{-1.6 cm}
    	\caption{Local space-time  grid where approximations of $U$ are computed to calculate $F^P_{i+1/2}$ with $P = 2$. For simplicity a pair $j,r$ represents the point $(x_{i+j}, t_{n+r})$. Taylor expansions in time are used to obtain these approximations following the blue lines. These Taylor expansions are centered in the points lying on the black line.}
    	\label{1D_grid}
    \end{figure}

    \begin{remark}
    Observe that the computation of the numerical flux  $F_{i+\ha}^{P}$ requires the approximation of $U$ at the nodes of a space-time grid of $2P \times 2P$ points: $U_{i,j}^{k,n+r}$,  $j,r = -P+1,\ldots,P$ (see Figure \ref{1D_grid}). The approximations of the solution $U$ at times $(n-P+1)\Delta t$, \dots, $(n-1)\Delta t$ are different from the ones  already computed in the previous steps: $U_{i+j}^{n-p}$, \dots, $U_{i+j}^{n-1}$.  In other words, the discretization in time is not based on a multistep method but in a one-step one: in fact it can be interpreted as a RK method whose stages are $\tilde U^{n+r}_{i,j}$, $r = -p, \dots, p$: see \cite{CPZMR2020}.

    \end{remark}

In \cite{CP2019} it has been shown that CAT$2P$ has order of accuracy $2P$ in time and it is linearly stable under the usual CFL condition. 

        For $P=1$, the expression of CAT2 numerical flux reduces to:
        \begin{equation}
            \label{F_CAT2}
            F_{i+\ha}^1 = \frac{1}{4}\Bigl(F_{i,0}^{1,n+1} + F_{i,1}^{1,n+1} + F(U_{i}^{n}) + F(U_{i+1}^{n})\Bigr).
        \end{equation}
        where, for $j=0,1,$
        \begin{equation*}
            F_{i,j}^{1,n+1} = F\Bigl(U_{i,j}^{1,n+1}\Bigr) = F\Bigl(U_{i+j}^n -\frac{\Delta t}{\Delta x}(F(U_{i+1}^n)-F(U_i^n))\Bigr) 
        \end{equation*}
        Consequently, the corresponding numerical scheme (\ref{CAT2P}) writes as follows:
        \begin{equation}
            \label{CAT2}
            U_{i}^{n+1} = U_i^n + \frac{\Delta t}{4\Delta x}\Bigl(F_{i-1,0}^{1,n+1} + F_{i-1,1}^{1,n+1} - F_{i,0}^{1,n+1} - F_{i,1}^{1,n+1} + F(U_{i-1}^{n}) -  F(U_{i+1}^{n})\Bigr).
        \end{equation}
        \begin{remark}
            The second order CAT2 numerical flux \eqref{F_CAT2} could be seen as a new one-step \textit{Jacobian-free} extension of Lax-Wendroff scheme that is even different from the two-step MacCormack and Richtmyer method.
        \end{remark}

    \section{Compact Approximate methods for systems of balance laws}
    \subsection{CAT2P methods}\label{ss_CAT2Pbl}
    The goal of this section is to extend CAT$2P$ methods to systems of balance laws \eqref{bal_sis}
	with initial condition $U(x,0) = U_0(x)$. Many authors have treated the source term as divergence of flux introducing an elliptic equation, see \cite{russo1990deterministic,russo1990particle}. To avoid the introduction of a new equation we follow the strategy in \cite{Gascon} (see also  \cite{Donat10}) in which the sistem \eqref{bal_sis} is first written in  conservative form through the definition of a ‘combined flux’ formed by the sum of  flux function  $F$  and the indefinite integral of the source term: more precisely,  let us introduce the function $\F$ given by 
    \begin{equation}
        \label{flux}
        \F(U)(x,t)= F(U(x,t))
        -\int_{-\infty}^{x}S(U(\sigma,t))H_x(\sigma)\,d\sigma,
    \end{equation}
    assuming that the integral is finite. Then, the equality
    $$\F(U)_x = F(U)_x - S(U)H_x,$$
allows one to write    the system of balance laws (\ref{bal_sis}) in the form
    \begin{equation}
        \label{con_equ}
        U_t + \F(U)_x = 0.
    \end{equation}
Then, the formal expression of CAT$2P$ is given by:
    \begin{equation}
        \label{2P_scheme}
        U_i^{n+1} = U_i^n + \frac{\Delta t}{\Delta x}\Big( \mathfrak{F}_{i-\ha}^P - \mathfrak{F}_{i+\ha}^P\Big),
    \end{equation}
    where
    \begin{equation}
        \label{con_flux} 
        \mathfrak{F}_{i+\ha}^P = \sum_{k=1}^{2P}\frac{\Delta t^{k-1}}{k!}A_{P}^{0,1/2}\Big(\F_{i,\bcdot}^{(k-1)}\Big).
    \end{equation}
    Here, $\F^{(k)}_{i,j}$ are local  approximation of $\partial_t^{(k)}\F(U)(x_{i+j},t_n)$ that are computed by adapting the algorithm described in Section \ref{S:cl}. Formally, the algorithm is as follows:
         \begin{itemize}
        \item Define 
        \begin{eqnarray*} 
        F_{i,j}^{(0)} &:=&  F(U_{i+j}^n), \quad j = -P+1,\ldots,P; \\
 I_{i,j}^{(0)} &:=&  \int_{-\infty}^{x_{i+j}}S(U(x,t_n))H_x(x)\, dx, \quad j=-P+1, \dots, P.
        \end{eqnarray*}
        
        \item For $k=1,\ldots,2P-1:$
            \begin{itemize}
                \item Compute for all $j=-P+1,\ldots,P$ 
                $$ U_{i,j}^{(k)} = -A_{P}^{1,j}\Bigl(F_{i, \bcdot}^{(k-1)},\Delta x\Bigr)
                +  A_{P}^{1,j}\Bigl(I_{i,\bcdot}^{(k-1)},\Delta x\Bigr). $$
                \item Define for all $j,r=-P+1,\ldots,P$
                \begin{eqnarray*}
                I_{i,j}^{n+r} &:= & \int_{-\infty}^{x_{i+j}}S(U(x,t_{n+r}))H_x(x)\, dx, \\
                F_{i,j}^{k,n+r} &:= & F\Bigl(U_{i,j}^{k,n+r}\Bigr), 
                \end{eqnarray*}
                where $U_{i,j}^{k,n+r}$ is the 
                approximation of $U(x_{i+j}, t_{n+r})$ given by the Taylor expansion in time:
                $$ U_{i,j}^{k,n+r} = U_{i+j}^n + \sum_{m = 1}^{k}\frac{(\Delta t)^m}{m!}U_{i,j}^{(m)}. $$
                \item Compute for all $j=-P+1,\ldots,P$
                $$
                F_{i,j}^{(k)}  =  A_{P}^{k,0}\Bigl(F_{i,j}^{k, \bcdot},\Delta t\Bigr), \quad
                I_{i,j}^{(k)}  =  A_{P}^{k,0}\Bigl(I_{i,j}^{\bcdot}, \Delta t \Bigr).
                $$
            \end{itemize}
    \end{itemize} 
The 'numerical fluxes' are then defined by:
    \begin{equation}
    \label{flux_2P}
         \mathfrak{F}_{i+\ha}^P = F^P_{i+1/2}  - I^P_{i+1/2}
    \end{equation}
    where
    \begin{eqnarray}\label{CAT2P_flux}
        F^P_{i+1/2}  = \sum_{k=1}^{2P} \frac{\Delta t^{k-1}}{k!}A^{0, 1/2}_{P}(F_{i,\bcdot}^{(k-1)}, \Delta x) , \\
        \label{CAT2P_source}
        I^P_{i+1/2}  = \sum_{k=1}^{2P} \frac{\Delta t^{k-1}}{k!}A^{0, 1/2}_{P}(I_{i, \bcdot}^{(k-1)}, \Delta x).
    \end{eqnarray}

  This algorithm is formal, since it requires the computation of integrals that depend on the exact solution in intervals of the form $(-\infty, x_{i+j}]$. In order to be computationally implementable, let us first rewrite it using only integrals in bounded intervals.  To do that, the key point is the following chain of equalities:
    \begin{eqnarray*}
     A_{P}^{1,j}\Bigl(I_{i, \bcdot}^{(k-1)}, \Delta x\Bigr) & = &
     \frac{1}{\Delta x} \sum_{s= -P + 1}^P \gamma_{P,s}^{1,j} I_{i,s}^{(k-1)} \\
 & = &  \frac{1}{\Delta x} \sum_{s= -P + 1}^P \gamma_{P,s}^{1,j} A_{P}^{k-1,0}
 \Bigl(I_{i,s}^{\bcdot}, \Delta t \Bigr)\\
 & =  &  \frac{1}{\Delta x \Delta t^{k-1}} 
 \sum_{s= -P + 1}^P \gamma_{P,s}^{1,j} \sum_{r= -P + 1}^P \gamma_{P, r}^{k-1,0} I_{i,s}^{k, n+r} \\
 & = &  \frac{1}{\Delta x \Delta t^{k-1}} 
  \sum_{r= -P + 1}^P \gamma_{P, r}^{k-1,0}  \sum_{s= -P + 1}^P \gamma_{P,s}^{1,j}  \int_{-\infty}^{x_{i+s}} S(U(x,t_{n+r}))H_x(x) \, dx \\
   & = &  \frac{1}{\Delta x \Delta t^{k-1}} 
  \sum_{r= -P + 1}^P \gamma_{P, r}^{k-1,0}  \sum_{s= -P + 1}^P \gamma_{P,s}^{1,j} \Bigl( \int_{-\infty}^{x_{i+s}} S(U(x,t_{n+r}))H_x(x) \, dx \\
  & & \qquad  - \int_{-\infty}^{x_{i-P +1}} S(U(x,t_{n+r}))H_x(x) \, dx \Bigr)\\
   & = &  \frac{1}{\Delta x \Delta t^{k-1}} 
  \sum_{r= -P + 1}^P \gamma_{P, r}^{k-1,0}  \sum_{s= -P + 1}^P \gamma_{P,s}^{1,j}  \int_{x_{i-P + 1}}^{x_{i+s}} S(U(x,t_{n+r}))H_x(x) \, dx
    \end{eqnarray*}
    where the identity
           $$\sum_{s=-P +1}^P \gamma^{1,j}_{P,s}  = 0, \quad j=-P+1,\ldots,P,$$
            has been used: remember that an interpolatory formula of numerical differentiation 
            that uses $2P$ points is exact at least for polynomials of degree $2P-1$ and thus it is exact for constant polynomials. Therefore, if the formula is applied to the constant polynomial $p \equiv 1$, we have
 $$
 0 = p'(x_{i+j}) = A^{1,j}_{P}(p_{i,\bcdot}, \Delta x) = \sum_{s=-P +1}^P \gamma^{1,j}_{P,s}p(x_{i+s}) 
 = \sum_{s=-P +1}^P \gamma^{1,j}_{P,s}.
 $$
 By introducing the notation
\begin{eqnarray*}
I_{i,j,l}^m & := & \int_{x_{i+j}}^{x_{i+l}} S(U(x,t_{m})) H_x(x) \,dx,\\
I_{i,j,l}^{(k-1)}  & : = &  A_{P}^{k-1,0}\Bigl(I_{i,j,l}^{\bcdot}, \Delta t \Bigr),
\end{eqnarray*}
we obtain 
    \begin{eqnarray*}
     A_{P}^{1,j}\Bigl(I_{i, \bcdot}^{(k-1)}, \Delta x\Bigr) & = &
  \frac{1}{\Delta x \Delta t^{k-1}} 
  \sum_{r= -P + 1}^P \gamma_{P, r}^{k-1,0}  \sum_{s= -P + 1}^P \gamma_{P,s}^{1,j} I_{i,-P+1,s}^{n+r}\\
  & = &
  \frac{1}{\Delta x} 
  \sum_{r= -P + 1}^P \gamma_{P, r}^{k-1,0} I_{i,-P+1,s}^{(k-1)}\\
  & = &  A_{P}^{1,j}\Bigl(I_{i,-P+1, \bcdot}^{(k-1)}, \Delta x\Bigr), 
    \end{eqnarray*}
  where only integrals in intervals of the form $[x_{i-P+1}, x_{i+s}]$  appear.
 Observe that $I_{i-P + 1, i-P + 1}^{n+r} = 0$ for all $k$.
 
 Concerning the expression of the numerical method, observe that:
    \begin{eqnarray*}
      I^P_{i+1/2}  -  I^P_{i-1/2} & = &  
     \sum_{k=1}^{2P} \frac{\Delta t^{k-1}}{k!}\left( A^{0, 1/2}_{P}(I_{i,\bcdot}^{(k-1)}, \Delta x) 
     -  A^{0, 1/2}_{P}(I_{i-1,\bcdot}^{(k-1)}, \Delta x) \right) \\
& = &  \sum_{k=1}^{2P} \frac{\Delta t^{k-1}}{k!} \sum_{j= -P +1}^P \gamma_{P,j}^{0,1/2} \left( I_{i,j}^{(k-1)} - I_{i-1,j}^{(k-1)} \right)\\
& =  & \sum_{k=1}^{2P} \frac{\Delta t^{k-1}}{k!} \sum_{j= -P +1}^P \gamma_{P,j}^{0,1/2} 
A_{P}^{k-1,0} \left(I_{i,j}^{\bcdot} - I_{i-1,j}^{\bcdot}, \Delta t\right).\\
\end{eqnarray*}
Since
\begin{equation}\label{calI}
I_{i,j}^{n + r} - I_{i-1,j}^{n + r}  = I_{i,j-1,j}^{n+r} = \int_{x_{i+j-1}}^{x_{i+j}} S(U(x,t_{n+r})) H_x (x) \, dx,
\end{equation}
 if we define 
 \begin{equation}\label{calI(k-1)}
  \mathcal{I}_{i,j}^{(k-1)}  =  A_{P}^{k-1,0}\Bigl({I}_{i,j-1,j}^{\bcdot}, \Delta t \Bigr),
 \end{equation}
we have
$$
 I^P_{i+1/2}  -  I^P_{i-1/2}  =  
     \sum_{k=1}^{2P} \frac{\Delta t^{k-1}}{k!}A^{0, 1/2}_{P}(\mathcal{I}_{i,\bcdot}^{(k-1)}, \Delta x),
$$
so that \eqref{2P_scheme} can be written in equivalent form
    \begin{equation}
        \label{2P_scheme_fs}
        U_i^{n+1} = U_i^n + \frac{\Delta t}{\Delta x}\Big(F_{i-\ha}^P - {F}_{i+\ha}^P + S_i^P\Big),
    \end{equation}
    where
    \begin{equation} \label{source}
    S_i^P =  \sum_{k=1}^{2P} \frac{\Delta t^{k-1}}{k!}A^{0, 1/2}_{P}(\mathcal{I}_{i,\bcdot}^{(k-1)}, \Delta x).
    \end{equation}
   Observe that only integrals \eqref{calI} in intervals of length $\Delta x$ appear in the expression of the numerical source term.

Finally, in order to have an implementable algorithm, all the integrals appearing in it are approximated  using quadrature formulas combined with the approximations $U_{i,j}^{k, n+r}$ of the exact solution that are available at every stage. To do this, given $i$ and $j=-P+2, \dots, P$, we consider at $[x_{i+j-1}, x_{i+j}]$ the interpolatory quadrature formula
$$
\int_{x_{i+j-1}}^{x_{i+j}} f(x) \, dx \approx \Delta x \sum_{s= -P + 1}^P a^{i,j}_{P,s} f(x_{i+s}) 
$$
whose nodes are $x_{i+s}$, $s = -P+1, \dots, P$.
This formula will be used to approximate the integrals appearing at the $k$-th stage of the algorithm  as follows: given two indices $j_1 < j_2$
\begin{eqnarray*}
I_{i, j_1,j_2}^m \approx \widetilde{I}^{k,m}_{i,j_1,j_2} := \Delta x \sum_{s=j_1+1}^{j_2} \sum_{l = -P + 1}^P a^{i,s}_{P,l} S(U^{k,m}_{i, l})H_x(x_{i+l}).
\end{eqnarray*}

Taking into account these approximations of the integral terms, the algorithm is finally as follows, see Figure \ref{1D_grid} for the case $P=2:$
         \begin{itemize}
        \item Compute 
        \begin{eqnarray*} 
        & & F_{i,j}^{(0)} =  F(U_{i+j}^n), \quad j = -P+1,\ldots,P; \\
        & &  \widetilde{I}_{i,j-1,j}^{(0)}  =  \Delta x \sum_{l= -P+1}^P a^{i,j}_{P,l} S(U^{n}_{i+l})H_x(x_{i+l}), \quad j=-P+2, \dots, P; \\
        & &   \widetilde{I}_{i,-P+1,-P+1}^{(0)}  =   0;\\
        & & \widetilde{I}_{i,-P+1,j}^{(0)}  =   \sum_{s= -P+2}^j   \widetilde{I}_{i,s-1,s}^{(0)}, \quad j=- P+2, \dots, P.
        \end{eqnarray*}
        
        \item For $k=1,\ldots,2P-1:$
            \begin{itemize}
                \item Compute for all $j=-P+1,\ldots,P$ 
                $$ U_{i,j}^{(k)} = -\A_{P}^{1,j}\Bigl(F_{i, \bcdot}^{(k-1)},\Delta x\Bigr)
                +  \A_{P}^{1,j}\Bigl(I_{i,-P+1,\bcdot}^{(k-1)},\Delta x\Bigr). $$
                \item Compute for all $j,r=-P+1, \dots, P$
      $$ U_{i,j}^{k,n+r} = U_{i+j}^n + \sum_{m = 1}^{k}\frac{(r\Delta t)^m}{m!}U_{i,j}^{(m)}. $$          
                \item Compute for all $j, r = -P +1, \ldots, P$

        $$
         F_{i,j}^{k,n+r} =  F\Bigl(U_{i,j}^{k,n+r}\Bigr),
         $$
         \item Compute for all $r = -P+1, \dots, P$, $j = -P+2, \dots, P$
         $$
         \widetilde{I}^{k,n+r}_{i,j-1,j} = \Delta x \sum_{l = -P + 1}^P a^{i,j}_{P,l} S(U^{k, n+r}_{i, l})H_x(x_{i+l}).
         $$

         \item Compute for all $j = -P+2, \dots, P$
         $$
         \widetilde{I}^{(k)}_{i,j-1,j} =  A_{P}^{k,0}\Bigl(\widetilde{I}_{i,j-1,j}^{k,\bcdot}, \Delta t \Bigr).
         $$
         
         \item Compute
         \begin{eqnarray*}
               & &  F_{i,j}^{(k)}  =  A_{P}^{k,0}\Bigl(F_{i,j}^{k, \bcdot},\Delta t\Bigr), \quad j=-P+1,\ldots,P; \\
            & & \widetilde{I}_{i,-P+1,-P+1}^{(k)}  =  0;\\
             & &    \widetilde{I}_{i,-P+1,j}^{(k)}  = \sum_{s=-P+2}^{j}  \widetilde{I}^{(k)}_{i,s-1,s}\quad j=-P+2,\ldots,P.
        \end{eqnarray*}

            \end{itemize}
    \end{itemize} 
    
    Once the algorithm has been executed, the integrals already computed can be used to approximate the source term  as follows:
    
    \begin{itemize}
 
        \item For $k=1,\dots, 2P$ define
        $$
          \widetilde{\mathcal{I}}_{i,j}^{(k-1)}  =  \begin{cases}
        \widetilde{I}^{(k-1)}_{i-1,j,j+1} & \text{if $j = -P+1, \ldots, 0$;} \\
        & \\
        \widetilde{I}^{(k-1)}_{i, j-1, j} & \text{if $j = 1,\dots, P$. }
            \end{cases}
            $$
          
          \item Compute
           $$
    \widetilde{S}_i^P =  \sum_{k=1}^{2P} \frac{\Delta t^{k-1}}{k!}A^{0, 1/2}_{P}(\widetilde{\mathcal{I}}_{i, \bcdot}^{(k-1)}, \Delta x).
  $$
          \end{itemize}
          
Observe that the first $P$ integral terms $\mathcal{I}_{i,j}^{(k-1)}$ appearing in the expression of the numerical source term \eqref{source} are approximated with the values $I_{i-1, j, j+1}^{(k-1)}$, used to compute the flux at the intercell $i-1/2$,  and the $P$ last ones by $I_{i, j-1, j}^{(k-1)}$, used to compute the flux at the intercell 
$i +1/2$. 

The final expression of the numerical method is then
   \begin{equation}
        \label{2P_scheme_fs_qf}
        U_i^{n+1} = U_i^n + \frac{\Delta t}{\Delta x}\Big(F_{i-\ha}^P - {F}_{i+\ha}^P + \widetilde{S}_i^P\Big),
    \end{equation}
    where ${F}_{i+\ha}^P$ is given by \eqref{CAT2P_flux}.

  \subsection{CAT2 for system of balance laws}
  Let us illustrate the above numerical method in the easiest case $P = 1$. In this case, the quadrature formula used to compute integrals in intervals of length $\Delta x$ is  the trapezoidal rule:
  $$
\int_{x_{i}}^{x_{i+1}} f(x) \, dx \approx \frac{\Delta x}{2} \Bigl( f(x_{i}) + f(x_{i+1})\Bigr) .
$$
The numerical method is then as follows: for every $i$
\begin{itemize}
    \item Compute
        $$
    U_{i,j}^{(1)}  =  -\frac{1}{\Delta x}\left(F(U_{i+1}^n) - F(U_i^n) \right)  + \frac{1}{2} \left( S(U_i^n) H_x (x_i) + S(U_{i+1}^n ) H_x(x_{i+1}) \right), \quad j= 0,1.
        $$
        Observe that, there is no dependence on $j$ in the right term because the derivative are computed in the same way.
        
    \item Compute
        $$
    U_{i,j}^{1,n+1} = U_{i+j}^{n} + \Delta t\;U_{i,j}^{(1)}, \quad j=0,1.
    $$
    
\end{itemize}

Then, define
   \begin{equation}\label{CAT2_flux}
    {F}_{i + \ha}^1  := \frac{1}{4} \left(  F(U_{i}^n) + F(U_{i+1}^n) + F(U_{i,0}^{1,n+1})  +F(U_{i,1}^{1,n+1}) \right) 
    \end{equation}
\begin{eqnarray}
 \widetilde{S}_{i}^1&:= & \frac{\Delta x}{8} \Bigl( (S(U_{i-1}^n) + S(U_{i-1,0}^{1, n+1})) H_x (x_{i-1}) +   (S(U_{i}^n) + S(U_{i-1,1}^{1, n+1})) H_x (x_{i})  \label{CAT2_source}\\\nonumber
        &  & +  (S(U_{i}^n) + S(U_{i,0}^{1, n+1})) H_x (x_{i}) + (S(U_{i+1}^n) + S(U_{i,1}^{1, n+1})) H_x (x_{i+1}) \Bigr).   
        \end{eqnarray}
The numerical method is then \eqref{2P_scheme_fs_qf} with $P = 1.$
\begin{remark}
The CAT procedure applied to systems of balance law introduces spurious oscillations near discontinuity points. In Section \ref{adaptive} we introduce the adaptive version of those methods, that prevent formation of spurious oscillations.
\end{remark}

\section{Well-balanced CAT2P for systems of balance laws}
\subsection{WBCAT2P methods}\label{ss:WBCAT2Pbl}
The goal of this section is to derive a well-balanced version of the CAT$2P$ methods introduced in the previous section. The idea is as follows: let us suppose that the initial  condition is given by
$$
U(x,0) = U^*(x),
$$
were $U^*$ is a stationary solution of  \eqref{bal_sis}. Let us introduce then
the function $\Tilde{\F}$ given by
    \begin{align}
        \Tilde{\F}(U)(x,t) = & \F(U)(x,t)) - \F(U^{*})(x) =\nonumber \\ 
        = & F(U(x,t)) - F(U^*(x))   -\int_{-\infty}^{x}\Bigl(S(U(\sigma,t)) - S(U^{*}(\sigma))\Bigr)H_{x}(\sigma)d\sigma.
    \end{align}
    Hence, observing that 
    $$\Tilde{\F}(U)_x = F(U)_x - F(U^*)_x - (S(U)-S(U^*))H_x = F(U)_x - S(U)H_x,
    $$
    the system of balance laws (\ref{bal_sis}) can be formally written in the form
    \begin{equation}
        U_t + \Tilde{\F}(U)_x = 0.
    \end{equation}
    Obviously $\Tilde{\F}(U^*) = 0$, therefore a numerical method based on the discretization of this conservative form
 is expected to preserve $U^*$ exactly.
    
    In  practice, this strategy is applied as follows: once the approximation $U_i^n$ has been obtained, we consider the stationary solution $U^{*}_i$ that satisfies
    $$
U^{*}_i (x_i) = U_i^n,
$$
i.e. $U^{*}_i$ solves the Cauchy problem
    \begin{equation}\label{Cauchyin}
    \left \{
    \begin{array}{l}
    F(U)_x = S(U) H_x \\
    U(x_i) = U_i^n.
    \end{array}
    \right.
    \end{equation}
    Let us assume for simplicity that this Cauchy problem has a unique solution that is explicitly known. Then, the system of balance laws is rewritten in the form 
       \begin{equation}
        \label{con_equ*}
        U_t + \Tilde{\F}_i(U)_x = 0.
    \end{equation}
 where   
    \begin{align}
        \Tilde{\F}_i(U)(x,t) = & \F(U)(x,t)) - \F(U^{*}_i)(x) =\nonumber \\ \label{flux*}
        = & F(U(x,t)) - F(U^*_i(x))   -\int_{-\infty}^{x}\Bigl(S(U(\sigma,t)) - S(U^{*}_i(\sigma))\Bigr)H_{x}(\sigma)\, d\sigma
    \end{align}
 and the  CAT$2P$ method is then applied:
    \begin{equation}
        \label{2P_scheme*}
        U_i^{n+1} = U_i^n + \frac{\Delta t}{\Delta x}\Big( \mathfrak{\Tilde{F}}_{i,i-\ha}^P - \mathfrak{\Tilde{F}}_{i,i+\ha}^P\Big),
    \end{equation}
    where,
    \begin{align}
        \label{con_flux*} 
        \mathfrak{\Tilde{F}}_{i,i+\ha}^P &= \sum_{k=1}^{2P}\frac{\Delta t^{k-1}}{k!}A_{P}^{0,1/2}\Big(\Tilde\F_{i;i,\bcdot}^{(k-1)}\Big), \\
        \mathfrak{\Tilde{F}}_{i,i-\ha}^P &= \sum_{k=1}^{2P}\frac{\Delta t^{k-1}}{k!}A_{P}^{0,1/2}\Big(\Tilde\F_{i;i-1,\bcdot}^{(k-1)}\Big).
    \end{align}
    Here  $\Tilde\F^{(k)}_{i;l,j}$ is an approximation of $\partial_t^{(k)}\Tilde\F_i(U)(x_{l+j},t_n).$ Observe that, in this case, two numerical fluxes  have to be computed at every inter-cell $x_{i+1/2}$,  $\mathfrak{\Tilde{F}}_{i,i+\ha}^P$ and
     $\mathfrak{\Tilde{F}}_{i+1,i+\ha}^P$, whose computation are based respectively on the stationary solutions $U^*_i$ 
     (that satisfies $U^*_i(x_i) = U_i^n$) and  $U^*_{i+1}$ (that satisfies $U^*_{i+1}(x_{i+1}) = U_{i+1}^n$).
    
    Observe that, if the initial  condition $U_0$ is a stationary solution, then at time $t = 0$,
    $U^*_i = U_0$ for all $i$,  so that 
    $$
    \Tilde{\F}_i(U_0) =  \F(U_0) - \F(U^{*}_i) = \F(U_0) - \F(U_0) = 0, \quad \forall i,
    $$
    
    and the numerical method is expected to preserve the initial condition. 
    
    
    The algorithm is then as follows: for every $i$
         \begin{itemize}
         \item Compute the solution $U^*_i(x)$ of the Cauchy problem \eqref{Cauchyin}.
        \item Compute 
        \begin{eqnarray*} 
        & & F_{i;i,j}^{(0)} =  F(U_{i+j}^n) - F(U^*_i(x_{i+j})), \quad j = -P+1,\ldots,P; \\
        & & F_{i;i-1,j}^{(0)} =  F(U_{i-1+j}^n) - F(U^*_i(x_{i-1+j})), \quad j = -P+1,\ldots,P; \\
        & &  \widetilde{I}_{i;i,j-1,j}^{(0)}  = 
        \Delta x \sum_{l= -P+1}^P a^{i,j}_{P,l}\left( S(U^{n}_{i+l}) - S(U^*_i(x_{i+l})) \right) H_x(x_{i+l}), \quad j=-P+2, \dots, P; \\
        & &   \widetilde{I}_{i;i,-P+1,-P+1}^{(0)}  =   0;\\
        & & \widetilde{I}_{i;i,-P+1,j}^{(0)}  =   
        \sum_{s= -P+2}^j   \widetilde{I}_{i;i,s-1,s}^{(0)}, \quad j=- P+3, \dots, P;\\
        & &  \widetilde{I}_{i;i-1,j-1,j}^{(0)}  = 
        \Delta x \sum_{l= -P+1}^P a^{i-1,j}_{P,l}\left( S(U^{n}_{i-1+l}) - S(U^*_i(x_{i-1+l})) \right) H_x(x_{i-1+l}), \quad j=-P+3, \dots, P;\\
         & &   \widetilde{I}_{i;i-1,-P+1,-P+1}^{(0)}  =   0;\\
        & & \widetilde{I}_{i;i-1,-P+1,j}^{(0)}  =   \sum_{s= -P+2}^j   \widetilde{I}_{i;i-1,s-1,s}^{(0)}, \quad j=- P+2, \dots, P.
        \end{eqnarray*}
        
        \item For $k=1,\ldots,2P-1:$
            \begin{itemize}
                \item Compute for all $j=-P+1,\ldots,P$ 
                \begin{eqnarray*}
              & &   U_{i;i,j}^{(k)} = -A_{P}^{1,j}\Bigl(F_{i;i, \bcdot}^{(k-1)},\Delta x\Bigr)
                +  A_{P}^{1,j}\Bigl(I_{i;i,-P+1,\bcdot}^{(k-1)},\Delta x\Bigr);\\ 
              & &  U_{i;i-1,j}^{(k)} = -A_{P}^{1,j}\Bigl(F_{i;i-1, \bcdot}^{(k-1)},\Delta x\Bigr)
                +  A_{P}^{1,j}\Bigl(I_{i;i-1,-P+1,\bcdot}^{(k-1)},\Delta x\Bigr). 
                \end{eqnarray*}
                \item Compute for all $j,r=-P+1, \dots, P$
        \begin{eqnarray*}        
     & &   U_{i;i,j}^{k,n+r} = U_{i+j}^n + \sum_{m = 1}^{k}\frac{(\Delta t)^m}{m!}U_{i;i,j}^{(m)}, \\        
     & &  U_{i;i-1,j}^{k,n+r} = U_{i+j-1}^n + \sum_{m = 1}^{k}\frac{(\Delta t)^m}{m!}U_{i;i-1,j}^{(m)}. 
      \end{eqnarray*}
                \item Compute for all $j, r = -P +1, \ldots, P$
        $$
         F_{i;i,j}^{k,n+r} =  F\Bigl(U_{i;i,j}^{k,n+r}\Bigr),
         \quad 
         F_{i;i-1,j}^{k,n+r} =  F\Bigl(U_{i;i-1,j}^{k,n+r}\Bigr).
         $$
         \item Compute for all $r = -P+1, \dots, P$, $j = -P+2, \dots, P$
         \begin{eqnarray*}
        & &  \widetilde{I}^{k,n+r}_{i;i,j-1,j} = 
         \Delta x \sum_{l = -P + 1}^P a^{i,j}_{P,l} S(U^{k, n+r}_{i;i, l}) H_x(x_{i+l}),\\
       & &   \widetilde{I}^{k,n+r}_{i;i-1,j-1,j} = 
         \Delta x \sum_{l = -P + 1}^P a^{i-1,j}_{P,l} S(U^{k, n+r}_{i;i-1, l}) H_x(x_{i-1+l}).
        \end{eqnarray*}
         \item Compute for all $j = -P+2, \dots, P$
         $$
         \widetilde{I}^{(k)}_{i;i,j-1,j} =  A_{P}^{k,0}\Bigl(\widetilde{I}_{i;i,j-1,j}^{k,\bcdot}, \Delta t \Bigr), \quad
         \widetilde{I}^{(k)}_{i;i-1,j-1,j} =  A_{P}^{k,0}\Bigl(\widetilde{I}_{i;i-1,j-1,j}^{k,\bcdot}, \Delta t \Bigr).
         $$
        
         \item Compute
         \begin{eqnarray*}
               & &  F_{i;i,j}^{(k)}  =  A_{P}^{k,0}\Bigl(F_{i;i,j}^{k, \bcdot},\Delta t\Bigr), \quad j=-P+1,\ldots,P; \\
            & & \widetilde{I}_{i;i,-P+1,-P+1}^{(k)}  =  0;\\
             & &    \widetilde{I}_{i;i,-P+1,j}^{(k)}  = \sum_{s=-P+2}^{j}  \widetilde{I}^{(k)}_{i;i,s-1,s}\quad j=-P+2,\ldots,P; \\
            & &  F_{i;i-1,j}^{(k)}  =  A_{P}^{k,0}\Bigl(F_{i;i-1,j}^{k, \bcdot},\Delta t\Bigr), \quad j=-P+1,\ldots,P; \\
            & & \widetilde{I}_{i;i-1,-P+1,-P+1}^{(k)}  =  0;\\
             & &    \widetilde{I}_{i;i-1,-P+1,j}^{(k)}  = \sum_{s=-P+2}^{j}  \widetilde{I}^{(k)}_{i;i-1,s-1,s}\quad j=-P+2,\ldots,P.
        \end{eqnarray*}

            \end{itemize}
    \end{itemize} 
    
    Once the algorithm has been executed, the integrals already computed can be used to approximate the source term  as follows:
    
    \begin{itemize}
 
        \item For $k=1,\dots, 2P$ define
        $$
          \widetilde{\mathcal{I}}_{i,j}^{(k-1)}  =  \begin{cases}
        \widetilde{I}^{(k-1)}_{i;i-1,j,j+1} & \text{if $j = -P+1, \ldots, 0$;} \\
        & \\
        \widetilde{I}^{(k-1)}_{i;i, j-1, j} & \text{if $j = 1,\dots, P$. }
            \end{cases}
            $$
          
          \item Compute
           \begin{equation}\label{WBCAT2P_source}
    \widetilde{S}_i^P =  \sum_{k=1}^{2P} \frac{\Delta t^{k-1}}{k!}A^{0, 1/2}_{P}(\widetilde{\mathcal{I}}_{i, \bcdot}^{(k-1)}, \Delta x).
  \end{equation}
          \end{itemize}
The final expression of the numerical method is then
   \begin{equation}
        \label{2P_scheme_wb_fs_qf}
        U_i^{n+1} = U_i^n + \frac{\Delta t}{\Delta x}\Big(F_{i;i-\ha}^P - {F}_{i;i+\ha}^P + \widetilde{S}_i^P\Big),
    \end{equation}
    where ${F}_{i,i\pm\ha}^P$ are given by
    \begin{eqnarray} \label{WBCAT2P_flux+}
          F^P_{i;i+1/2} & = & \sum_{k=1}^{2P} \frac{\Delta t^{k-1}}{k!}A^{0, 1/2}_{P}(F_{i;i,\bcdot}^{(k-1)}, \Delta x) , \\
          \label{WBCAT2P_flux-}
          F^P_{i;i-1/2} & = & \sum_{k=1}^{2P} \frac{\Delta t^{k-1}}{k!}A^{0, 1/2}_{P}(F_{i;i-1,\bcdot}^{(k-1)}, \Delta x) .
   \end{eqnarray}
   
   \begin{remark}
   Apparently scheme \ref{2P_scheme_wb_fs_qf} is not exactly conservative because the fluxes at the two sides of a cell edge may be different, since they depend on the local reconstruction. However, the schemes are exactly well balanced (see Theorem \ref{th:WB} below), and as a consequence, in case the source is identically zero, the numerical fluxes at the two sides of a cell edges are equal, and the scheme is conservative.
   \end{remark}
   
   \begin{remark}
   Observe that this algorithm can be used to update the solution at the point $x_i$ at time $t_n$ only if the Cauchy problem \eqref{Cauchyin} has a solution that is defined in the cells of the stencils $S^P_{i \pm 1/2}$ whose analytic expression is known.
   Therefore:
   \begin{itemize} 
   \item If \eqref{Cauchyin} has no solution, the CAT$2P$ method will be used instead. Please note that this choice does not spoil the well-balanced character of the numerical method: in this case, the cell values in the stencil cannot be the point values of a stationary solution (otherwise there would be at least one solution of \eqref{Cauchyin}) and thus there is no local equilibrium to preserve. 
  \item If \eqref{Cauchyin} has more than one solution, a criterion is needed to select one of them: this is the case for the shallow water system that will be discussed in Section \ref{ss_family}.
  \item If \eqref{Cauchyin} has a solution defined in the stencils but it is not possible to find its expression by analytic procedures, it is possible to apply an ODE solver to approximate it, like it has been done in \cite{CastroGomezPares} for finite-volume methods. In all the problems considered in Section \ref{nUm} the analytic expression of the stationary solutions is available either in explicit or implicit form.
  \end{itemize}
   \end{remark}

   \subsection{Well-balanced property}
   Numerical method \eqref{2P_scheme_wb_fs_qf} is fully well-balanced in the following sense:
   \begin{theorem}  \label{th:WB}
   Let $U^*$ be a continuous stationary solution of \eqref{bal_sis}. Then, if numerical method \eqref{2P_scheme_wb_fs_qf} is applied to the initial condition
   $$
   U^0_i = U^*(x_i), \quad \forall i,
   $$
   then
   $$
   U_i^n = U_i^0, \quad \forall i, n.
   $$
   \end{theorem}
   \noindent \textbf{Proof:} The proof is based on induction on the order of the proximate time derivative. Observe first that $U^*$ solves any Cauchy problem \eqref{Cauchyin} for $n = 0$. Therefore, at the first step the solution of \eqref{Cauchyin} is given by
   $$
   U^*_i = U^*, \quad \forall i.
   $$
   Therefore, for every $i$:
        \begin{eqnarray*} 
        & & F_{i;i,j}^{(0)} =  F_{i;i-1,j}^{(0)} = 0, \quad j = -P+1,\ldots,P; \\
        & &  \widetilde{I}_{i;i,j-1,j}^{(0)}  =   \widetilde{I}_{i;i-1,j-1,j}^{(0)}  = 0, \quad j= -P+2, \ldots, P;\\
        & & \widetilde{I}_{i;i,-P+1,j}^{(0)}  =   \widetilde{I}_{i;i-1,-P+1,j}^{(0)}  =  0, \quad j=- P+1, \dots, P;
        \end{eqnarray*}
        and thus
$$
U_{i;i,j}^{(1)} =   U_{i;i-1,j}^{(1)} = 0, \quad j = -P+1,\ldots, P.
$$
As a consequence:
\begin{eqnarray*}     
 & & U_{i;i,j}^{1,r} = U_{i+j}^0,  \quad F_{i;i,j}^{1,r} =  F\Bigl(U_{i+j}^{0}\Bigr), \quad j,r=-P+1, \dots, P;\\
 & &  U_{i;i-1,j}^{1,r} = U_{i-1+j}^0, \quad F_{i;i-1,j}^{1,r} =  F\Bigl(U_{i-1+j}^{0}\Bigr),  \quad j,r=-P+1, \dots, P\\
  & &  \widetilde{I}^{1,r}_{i;i,j-1,j} = 
         \Delta x \sum_{l = -P + 1}^P a^{i,j}_{P,l} S(U^0_{i+l}) H_x(x_{i+l}),\\
       & &   \widetilde{I}^{1,r}_{i;i-1,j-1,j} = 
         \Delta x \sum_{l = -P + 1}^P a^{i-1,j}_{P,l} S(U^{0}_{i-1+l}) H_x(x_{i-1+l}).
 \end{eqnarray*}
 Notice that the values of all these quantities do not depend on $r$. Therefore, when numerical differentiation in time is applied we obtain:
    \begin{eqnarray*}
               & &  F_{i;i,j}^{(1)}  =  F_{i;i-1,j}^{(1)}  = 0, \quad j=-P+1,\ldots,P; \\
                & &  \widetilde{I}_{i;i,j-1,j}^{(1)}  =   \widetilde{I}_{i;i-1,j-1,j}^{(1)}  = 0, \quad j= -P+2, \ldots, P;\\
            & & \widetilde{I}_{i;i,-P+1,j}^{(1)}  =  \widetilde{I}_{i;i-1,-P+1,j}^{(1)}  = 0, \quad j=-P+1,\ldots,P.
        \end{eqnarray*}
Therefore
$$
U_{i;i,j}^{(2)} =   U_{i;i-1,j}^{(2)} = 0, \quad j = -P+1,\ldots, P.
$$   
Repeating the reasoning we obtain
   \begin{eqnarray*}
               & &  F_{i;i,j}^{(k)}  =  F_{i;i-1,j}^{(k)}  = 0, \quad j=-P+1,\ldots,P, \quad k = 0, \dots, 2P; \\
                & &  \widetilde{I}_{i;i,j-1,j}^{(k)}  =   \widetilde{I}_{i;i-1,j-1,j}^{(k)}  = 0, \quad j= -P+2, \ldots, P, \quad k = 0, \dots, 2P;\\
            & & \widetilde{I}_{i;i,-P+1,j}^{(1)}  =  \widetilde{I}_{i;i-1,-P+1,j}^{(1)}  = 0, \quad j=-P+1,\ldots,P, \quad k = 0, \dots, 2P.
        \end{eqnarray*}
        Therefore
        $$
        F^P_{i; i+1/2} = F^P_{i, i-1/2} = S^P_{i+1/2} = 0,
        $$
        and we obtain
        $$
        U^1_i = U^0_i, \quad \forall i,
        $$
        as we wanted to prove.

  \subsection{WBCAT2 for system of balance laws}
  Let us illustrate again this numerical method in the  case $P = 1$: 
  
  For every $i$
\begin{itemize}
    \item Compute the solution $U^*_i$ of the Cauchy problem \eqref{Cauchyin}.
    \item Compute
    \begin{eqnarray*}
    U_{i;i,j}^{(1)}  & = &  -\frac{1}{\Delta x}\Bigl(F(U_{i+1}^n) - F(U_i^n)  -  F(U_{i}^*(x_{i+1})) + F(U_i^*(x_i)) \Bigr)\\
                     & & + \frac{1}{2} \Bigl( (S(U_i^n) - S(U_i^*(x_i))) H_x (x_i) + 
                       (S(U_{i+1}^n ) - S(U_i^*(x_{i+1})))H_x(x_{i+1}) \Bigr), \quad j= 0,1; \\
    U_{i;i-1,j}^{(1)}  & = &  -\frac{1}{\Delta x}\Bigl(F(U_{i}^n) - F(U_{i-1}^n)  -  F(U_{i}^*(x_{i})) + F(U_i^*(x_{i-1})) \Bigr)\\
                     & & + \frac{1}{2} \Bigl( (S(U_{i-1}^n) - S(U_i^*(x_{i-1}))) H_x (x_{i-1}) + 
                       (S(U_{i}^n ) - S(U_i^*(x_{i})))H_x(x_{i}) \Bigr), \quad j= 0,1.
    \end{eqnarray*}
        
    \item Compute
        \begin{eqnarray*}
   & &  U_{i;i,j}^{1,n+1} = U_{i+j}^{n} + \Delta t\;U_{i;i,j}^{(1)}, \quad  j = 0,1, \\
& & U_{i;i-1,j}^{1,n+1} = U_{i+j}^{n} + \Delta t\;U_{i-1;i,j}^{(1)},\quad j=0,1.
    \end{eqnarray*}

\item Define
    \begin{eqnarray}
  & &   {F}_{i,i + \ha}^1  := \frac{1}{4} \left(  F(U_{i}^n) + F(U_{i+1}^n) + F(U_{i;i,0}^{1,n+1})  +F(U_{i;i,1}^{1,n+1}) 
  - 2 F(U^*_i(x_i)) + \right. \\ \nonumber& & \qquad \qquad\quad - \left. 2 F(U^*_i(x_{i+1}))\right), \label{WBCAT2_flux1}\\ 
 & &  {F}_{i,i - \ha}^1  := \frac{1}{4} \left(  F(U_{i-1}^n) + F(U_{i}^n) + F(U_{i;i-1,0}^{1,n+1})  +F(U_{i;i-1,1}^{1,n+1}) 
  - 2 F(U^*_i(x_{i-1}))  \right. \\ \nonumber& & \qquad \qquad\quad - \left. - 2 F(U^*_i(x_{i}))\right),\label{WBCAT2_flux2}\\ \label{WBCAT2_source}
 & & \widetilde{S}_{i}^1 :=  \frac{\Delta x}{8} \Bigl( (S(U_{i-1}^n) + S(U_{i-1,0}^{1, n+1}) - 2S(U^*_i(x_{i-1}))) H_x (x_{i-1})\\\nonumber
 & & \qquad \qquad +   (S(U_{i}^n) + S(U_{i;i-1,1}^{1, n+1}) - 2S(U^*_i(x_i))) H_x (x_{i}) \\\nonumber
 & & \qquad \qquad + (S(U_{i}^n) + S(U_{i;i,0}^{1, n+1}) - 2S(U^*_i(x_i))) H_x (x_{i}) \\ \nonumber
 & & \qquad \qquad + (S(U_{i+1}^n) + S(U_{i; i,1}^{1, n+1}) - 2S(U^*_i(x_{i+1}))) H_x (x_{i+1}) \Bigr).   
        \end{eqnarray}
        
        \end{itemize}
        
The numerical method is then \eqref{2P_scheme_wb_fs_qf} with $P = 1$.

    \subsection{CAT2P methods that preserve a family of stationary solutions}\label{ss_family}

The strategy described in Subsection \ref{ss:WBCAT2Pbl}   can be easily adapted to obtain schemes that only preserve a prescribed set of stationary  solutions: this would be the case if, for instance, the set to be preserved is a $k$-parameter family of stationary solutions,
$$
U^*(x; C_1,\dots, C_k),
$$
with $k < d$ where $d$ is the number of variables. If it is the case, instead of looking for a solution of \eqref{Cauchyin}, one looks for a solution of the following nonlinear system:

  Find $C^i_1, \dots, C^i_k$ such that:
\begin{equation}\label{firststagek}
 u_{j_l}^* (x_i; C^i_1, \dots, C^i_k)  =  u_{i, j_l}, \quad l=1, \dots, k,
\end{equation}
where  $u^*_{j}$, $ u_{i,j}$ denote respectively the $j$-th component of $U^*$ and $U_i$ and $\{ j_1, \dots, j_k\}$ is a set of $k$ indices that is
predetermined in order to have the same number of unknowns and equations in \eqref{firststagek}.
These indices $j_1, \dots, j_k$ are chosen so that systems of equations \eqref{firststagek} have a unique solution, if it is possible. 
Once the problem has been solved, the numerical fluxes and source terms  are computed as in Section \ref{ss:WBCAT2Pbl} with the choice
$$
U^*_i(x) = U^*(x, C^i_1, \dots, C^i_k).
$$

    \section{Adaptive CAT2P  for systems of balance law}\label{adaptive}
    Despite the fact that Compact Approximate Taylor (CAT) schemes are linearly stable in the $L^2$ sense under the usual CFL condition, spurious oscillations may appear next to a discontinuity of the solution, as it happens for the Lax-Wendroff method:  see \cite{CP2019}. In order to get rid of them, we consider here the shock-capturing technique introduced in \cite{CPZMR2020} based on a family of high-order smoothness indicators. The idea is as follows: once the approximations at time $t^n$ have been computed, the candidate stencils
    to compute $\mathfrak{F}_{i+\ha}$ are
    $$ \mathcal{S}^p_{i+1/2} = \{x_{i-p+1},\ldots,x_{i+p}\}, \quad p = 1,\ldots,P.$$ 
    The selected stencil is the one with maximal length among those in which the solution at time $t^n$ is smooth, according to some smoothness indicators $\psi^p_{i+1/2}$ for $p = 1, \dots, P$.  
    If a discontinuity is detected in the stencil $\mathcal{S}^1_{i+1/2}$ a robust first-order numerical method is used. The ingredients of this strategy are described below:
    
    \subsection{First-order numerical method} 
     As first-order robust scheme to be combined with CAT2P methods for balance laws, we select the Lax-Friedrichs method applied to \eqref{con_equ} which leads to the formal expression:
     \begin{equation}\label{LFmeth}
     U_i^{n+1} = U_i^n + \frac{\Delta t}{\Delta x}\Bigl(F^{LF}_{i-1/2} - F^{LF}_{i+1/2} + \widetilde{S}^{LF}_i \Bigr)
     \end{equation}
     where
     \begin{eqnarray}
     F^{LF}_{i+1/2} & = & \frac{1}{2}\left(F(U^n_i) + F(U^n_{i+1}) \right) - \frac{\Delta x}{2 \Delta t} \left(U^n_{i+1} - U^n_i\right), \label{LF_flux}\\
     \widetilde S^{LF}_i & = & \Delta xS(U_i^n)H_x(x_i), \label{LF_source}
     \end{eqnarray}
     where the mid-point rule has been used to approximate the integral corresponding to the source term. 

In the case of the WBCAT2$P$ methods, the Lax-Friedrichs method is formally applied to \eqref{con_equ*} and the mid-point formula is again used to compute the integral terms, which leads again to a method of the form \eqref{LFmeth} where $F^{LF}_{i+1/2}$ is again the standard Lax-Friedrichs numerical flux \eqref{LF_flux} but now
\begin{equation}\label{LFWB_source}
\widetilde{S}^{LF}_i = {F}^{LF,*}_{i;i+1/2} -{F}^{LF,*}_{i;i-1/2},
\end{equation}
where
\begin{equation}\label{LF_flux*}
{F}^{LF,*}_{i;i+1/2}  = \frac{1}{2}\left(F(U_i^*(x_i)) + F(U_i^*(x_{i+1})) \right)
- \frac{\Delta x}{2 \Delta t} \left(U^*_i(x_{i+1}) - U^*_i(x_i)\right)
\end{equation}
and $U^*_i$ is the stationary solution that satisfies
$$
U^*_i(x_i ) = U_i^n.
$$
\eqref{LFWB_source} is a consistent approximation of the integral of the source term, since
$$
{F}^{ LF,*}_{i;i+1/2} -  {F}^{ LF,*}_{i;i-1/2} 
\approx \Delta x F(U^*)_x(x_i) = \Delta x S(U^*_i(x_i)) H_x (x_i) = \Delta x S(U^n_i)H_x(x_i).
$$
On the other hand, the well-balanced property of the method can be trivially checked.

\subsection{ACAT2 and WBACAT2 methods}
      The expression of the ACAT2 numerical method is based on a flux limiter (see \cite{LeVeque2002book,LeVeque2007book,Toro2009book}). Its expression is as follows:
     \begin{equation}\label{ACAT2meth}
     U_i^{n+1} = U_i^n + \frac{\Delta t}{\Delta x}\Bigl(F^{*}_{i-1/2} - F^{*}_{i+1/2} + \widetilde{S}^*_i \Bigr)
     \end{equation}
  where    
    \begin{eqnarray}\label{ACAT2_flux} 
    {F}^*_{i\pm1/2} & = & \varphi^1_{i} \, {F}^1_{i\pm1/2} +(1-\varphi^1_{i}) \, {F}^{LF}_{i\pm1/2},  \\
    \label{ACAT2_source}
    \widetilde{S}_i^* & = & \varphi^1_i \widetilde{S}^1_i + (1 - \varphi^1_i) \widetilde{S}^{LF}_i, 
    \end{eqnarray}
    where ${F}^1_{i\pm1/2}$ and $\widetilde{S}^1$ are given by \eqref{CAT2_flux}  and \eqref{CAT2_source} respectively;   $F^{LF}_{i\pm1/2}$ and $\widetilde{S}^{LF}_i$ are given by \eqref{LF_flux} and \eqref{LF_source} respectively; 
    $$
    \varphi^1_i = \min(\varphi^1_{i-1/2}, \varphi^1_{i+1/2}), 
    $$
    where $\varphi^1_{i+1/2}$ is computed by a  flux limiter, and has the property
    \begin{equation}
    \varphi^1_{i+1/2} \approx \begin{cases}
        1 &\text{if $\{ U_{i-1}^n, \dots, U_{i+2}^n \}$ are 'smooth';} \\
        0 & \text{otherwise.}
    \end{cases}
    \end{equation}

 For scalar problems,     
    standard flux limiter functions $ \varphi^1(r)$, such as minmod, superbee, van Leer \cite{Sweby,Kemm2010}, may  be used:
    \begin{equation}\label{indicador1}
    \varphi_{i+1/2}^1  =  \varphi^1(r_{i+1/2}), 
    \end{equation}
    where
    \begin{equation}\label{wavefuction}
        r_{i+1/2}=\frac{\Delta_{\rm upw}}{\Delta_{\rm loc}}
        = \left\{
        \begin{array}{cl}
            \displaystyle r^-_{i+1/2}:= \frac{u^n_i-u^n_{i-1}}{u^n_{i+1}-u^n_{i}} & \mbox {if } a_{i+1/2}  >0, \\
            \displaystyle  r^+_{i+1/2}:= \frac{u^n_{i+2}-u^n_{i+1}}{u^n_{i+1}-u^n_{i}} & \mbox {if } a_{i+1/2}  \leq 0;
        \end{array}\right.
    \end{equation}
    and  $a_{i+1/2}$ is an estimate of the wave speed such as for instance Roe's intermediate speed:
    $$a_{i+1/2} = \begin{cases}  \displaystyle \frac{f(u^n_{i+1}) - f(u^n_i)}{u^n_{i+1} - u^n_i} & \text{ if $|u^n_i - u^n_{i+1}|>{\rm tol};$}\\
    f'(u^n_i) & \text{otherwise.}
    \end{cases}$$
    where tol is a given small tolerance.
    An alternative that avoids the computation of an intermediate speed was introduced in \cite{Toro2009book} and it consists on defining
    \begin{equation}\label{local_indicators}
        \varphi^1_{i+1/2}=\min( \varphi^1(r_{i+1/2}^+), \varphi^1(r^-_{i+1/2})).
    \end{equation}
    This strategy can be easily extended to systems by computing the flux limiter component by component.

 The well-balanced version of the ACAT2 method has the form
      \begin{equation}\label{WBACAT2meth}
     U_i^{n+1} = U_i^n + \frac{\Delta t}{\Delta x}\Bigl(F^{*}_{i;i-1/2} - F^{*}_{i;i+1/2} + \widetilde{S}^*_i \Bigr)
     \end{equation}
  where    
    \begin{eqnarray}\label{WBACAT2_flux} 
    {F}^*_{i;i\pm 1/2} & = & \varphi^1_{i} \, {F}^1_{i; i \pm 1/2} +(1-\varphi^1_i) \, {F}^{LF}_{i;i\pm1/2},  \\
    \label{WBACAT2_source}
    \widetilde{S}_i^* & = & \varphi^1_{i} \widetilde{S}^{1}_i + (1 - \varphi^1_{i})\widetilde{S}^{LF,*}_i = \varphi^1_{i} \widetilde{S}^{1}_i -   (1 -  \varphi^1_{i})\left(F^{LF,*}_{i;i+\ha} - F^{LF,*}_{i;i-\ha}\right),
    \end{eqnarray}
    where ${F}^1_{i;i\pm 1/2}$ are given by \eqref{WBCAT2_flux1}-\eqref{WBCAT2_flux2};   $F^{LF,*}_{i;i\pm1/2}$ is  given by \eqref{LF_flux*}; and $\widetilde{S}_{i}^{1}$ is given by \eqref{WBCAT2_source}.

   \subsection{Smoothness Indicators}\label{ss:smoothness}
   The smoothness indicators introduced in  \cite{CPZMR2020} will be used here: given a set of the point values $f_j$ of a function 
   $f$ in  the nodes of the stencil $S^p_{i+1/2}$, $p \geq 2$  we define  $\psi_{i+1/2}^p$ as follows:  first define the lateral weights $I_{p,L}$ and $I_{p,R}$ as
    \begin{equation}
        I_{p,L}:=\sum_{j=-p+1}^{-1}(f_{i+1+j}-f_{i+j})^2+\varepsilon,\quad
        I_{p,R}:=\sum_{j=1}^{p-1}(f_{i+1+j}-f_{i+j})^2+\varepsilon,
    \end{equation}	
    where $\varepsilon$ is a small quantity that is added to prevent the lateral weights to vanish when the function is constant. Next, compute the half harmonic mean
    \begin{equation}
        I_{p} =\frac{I_{p,L}I_{p,R}}{I_{p,L}+I_{p,R}}.
    \end{equation}
    
     \begin{figure}[!ht]
    	\centering	
        \vspace{-0.75cm}
    	\includegraphics[scale = 1]{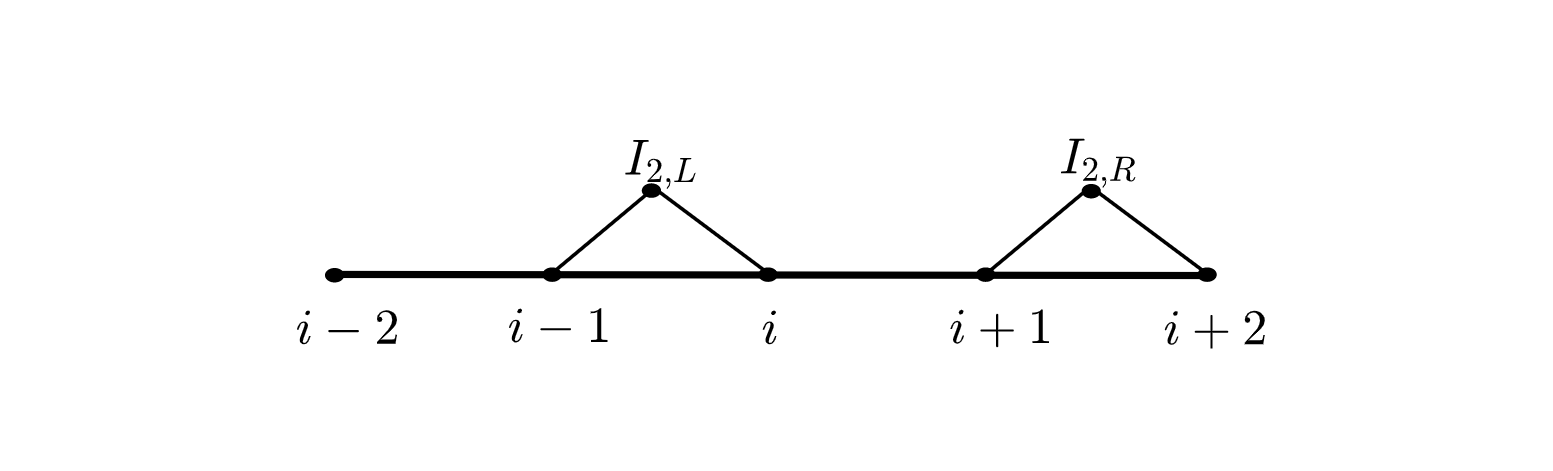}
    	\vspace{-0.8 cm}
    	\caption{Left and right part to compute the high order smoothness indicators where $P=2.$}
    	\label{smooth_fig}
    \end{figure}
    Finally, we define the high order smoothness indicator of the stencil $S_p$ by 
    \begin{equation}\label{indicadores_locales}
        \psi_{i+1/2}^p:=\left( \frac{I_p}{I_p+\tau_p}\right),
    \end{equation}
    where 
    \begin{align}
        \tau_p:=& \left( \Delta^{2p-1}_{i-p +1} f\right)^2.
    \end{align}
    The script $ \Delta^{2p-1}_{i-p +1} f$ represents the undivided difference of 
    $\{ f_{i-p+1}, \dots, f_{i+p}\}$:
    \begin{align}
        \Delta^{2p-1}_{i-p +1} f =&  (2p - 1)! \sum^{p}_{j=-p+1}\,\gamma^{2p-1, 1/2}_{p,j} \, f^n_{i+j}.
    \end{align}
    
       These indicators are such that
    \begin{equation}\label{condiconesS}
    \psi^p_{i+1/2} \approx \left\{
    \begin{array}{cl}
    1 & \mbox { if $\{f_j\}$ are 'smooth' in $\mathcal{S}^p_i$;}\\
    0 & \mbox {otherwise;} 
    \end{array}\right.
    \end{equation}
    see  \cite{CPZMR2020} for a precise statement of this property and its proof. 
       
    \subsection{ACAT2P and WBACAT2P methods}
    Using all these ingredients, the final expression of the ACAT2$P$ method is as follows:
      \begin{equation}
        \label{ACAT2P}
        U_i^{n+1} = U_i^n + \frac{\Delta t}{\Delta x}\Bigl({F}^{\mathcal{A}_i}_{i;i-\ha} - 
        {F}^{\mathcal{A}_i}_{i;i+\ha} + \widetilde{S}_i^{\mathcal{A}_i}\Bigr),
    \end{equation}  
 where
    \begin{equation}\label{numfluxacat}
    {F}^{\mathcal{A}_i}_{i;i\pm1/2} = \begin{cases}
    {F}^{*}_{i\pm1/2} & \text{if $\mathcal{A}_i = \emptyset$;}\\
    {F}^{p_s}_{i\pm1/2} & \text{otherwise, where $p_s = \max(\mathcal{A}_i)$;}
    \end{cases}
    \end{equation}
and
\begin{equation}\label{numsourceacat}
    \widetilde{S}^{\mathcal{A}_i}_{i} = \begin{cases}
    \widetilde{S}^*_i & \text{if $\mathcal{A}_i = \emptyset$;}\\
    \widetilde{S}^{p_s}_i & \text{otherwise, where $p_s = \max(\mathcal{A}_i).$} 
    \end{cases}
    \end{equation}
Here, $\mathcal{A}_i$ is the set of indices given by
    \begin{equation}
    \mathcal{A}_i= \{  p \in \{2, \dots, P \} \ s.t.\ \psi^p_{i-1/2} \approx 1  \text{ and }   \psi^p_{i+ 1/2} \approx 1 \};
    \label{eq:Ai}
    \end{equation}
$F^*_{i+1/2}$ and $\widetilde{S}_i^*$ are the ACAT2 numerical flux and source terms given by \eqref{ACAT2_flux}, (\ref{ACAT2_source}); $F^{p_s}_{i+1/2}$ and $\widetilde{S}_i^{p_s}$ are the ACAT2$p_s$ numerical fluxes and source terms defined in (\ref{flux_2P}), (\ref{source}). Throughout the paper we assume $\psi\approx 1 $ if $\psi\geq 0.9$.

\begin{remark}
Notice that in Equation \ref{eq:Ai} index $p$ starts from 2, since  it is not possible to determine the smoothness of the data in the two-point stencil   $S^1_{i+1/2}$: see \cite{CPZMR2020}.
\end{remark}

Analogously, the final expression of the WBACAT2$P$ method is as follows:
      \begin{equation}
        \label{WBACAT2P}
        U_i^{n+1} = U_i^n + \frac{\Delta t}{\Delta x}\Bigl({F}^{\mathcal{A}_i}_{i;i-\ha} - {F}^{\mathcal{A}_i}_{i;i+\ha} + \widetilde{S}_i^{\mathcal{A}_i}\Bigr),
    \end{equation}  
 where
     \begin{equation}\label{numfluxwbacat}
    {F}^{\mathcal{A}_i}_{i;i+1/2} = \begin{cases}
    {F}^{*}_{i;i+1/2} & \text{if $\mathcal{A}_i = \emptyset$;}\\
    {F}^{p_s}_{i;i+1/2} & \text{otherwise, where $p_s = \max(\mathcal{A}_i);$} 
    \end{cases}
    \end{equation}
    and
      \begin{equation}\label{numsourceacat_wb}
    \widetilde{S}^{\mathcal{A}_i}_{i} = \begin{cases}
    \widetilde{S}^*_i & \text{if $\mathcal{A}_i = \emptyset$;}\\
    \widetilde{S}^{p_s}_i  & \text{otherwise, where $p_s = \max(\mathcal{A}_i);$} 
    \end{cases}
    \end{equation} 
    where  $ {F}^{*}_{i;i\pm 1/2}$ and  $\widetilde{S}^*_i$ are the WBACAT2 numerical fluxes and source terms given by \eqref{WBACAT2_flux},
    \eqref{WBACAT2_source}; $  {F}^{p_s}_{i;i+1/2}$ and   $\widetilde{S}^{p_s}_i $ are WBCAT$2p_s$  numerical fluxes and source term 
    given by \eqref{WBCAT2P_flux-}, \eqref{WBCAT2P_flux+}, \eqref{WBCAT2P_source} with $P = p_s$.

\section{Two-dimensional problems} 
In this section we will focus on the extension of ACAT methods to non-linear two-dimensional systems of hyperbolic balance laws
\begin{equation}
\label{2Dequ}
{U}_t + {F}({U})_x + {G}({U})_y=S_1(U)H_x + S_2(U)H_y.
\end{equation}
As we did for the 1D systems of balance laws \eqref{bal_sis}, let us introduce the functions $\mathcal{F}$ and $\mathcal{G}$ as:
\begin{align}
    \mathcal{F}(U)((x,y),t) &= F(U(x,y,t)) - \int_{-\infty}^xS_1(U)H_\sigma d\sigma; \\
    \mathcal{G}(U)((x,y),t) &= G(U(x,y,t)) - \int_{-\infty}^yS_2(U)H_\tau d\tau, 
\end{align}
assuming that the integrals are finite. Then, the identities
$$ \mathcal{F}(U)_x = F(U(x,y,t))_x - S_1(U)H_x $$ and 
$$ \mathcal{G}(U)_y = G(U(x,y,t))_y - S_2(U)H_y $$
allows one to write the 2D system of balance laws \eqref{2Dequ} in the equivalent conservative form 
\begin{equation}
    \label{2Dcon_equ}
    U_t + \mathcal{F}(U)_x + \mathcal{G}(U)_y = 0. 
\end{equation}
The idea is now extend the ACAT$2P$ schemes to 2D systems written in the form \eqref{2Dcon_equ}.

\begin{figure}[h!]
	\centering
	\begin{picture}(80,175)
	\put(10,4){\makebox(60,165)[c]{		
			\includegraphics[scale= 0.5]{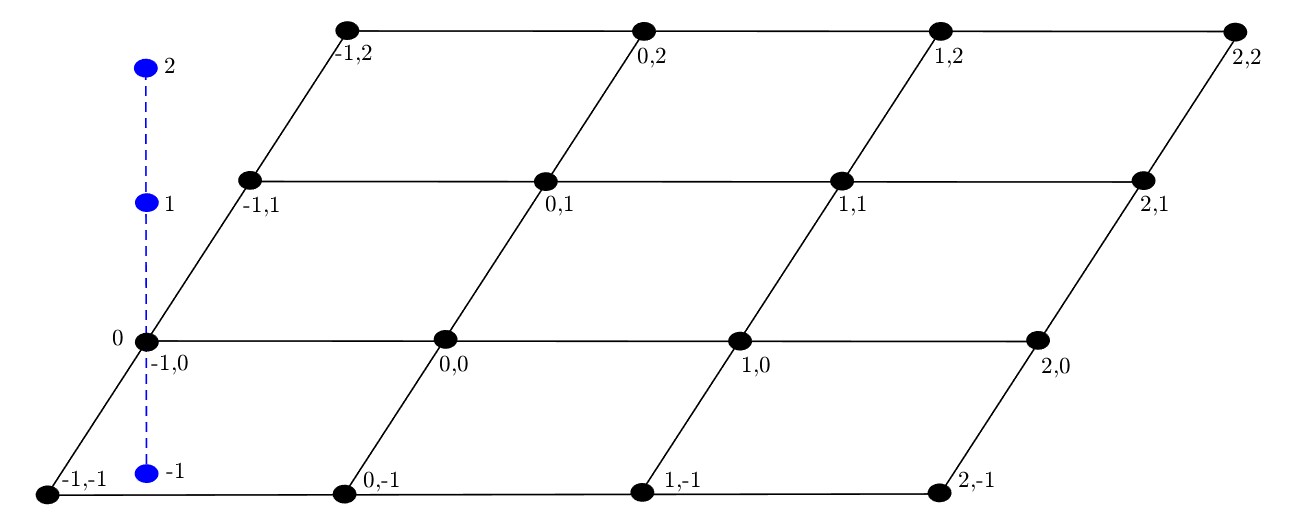}}}
	\end{picture}
	\vspace{-0.5cm}
	\caption{Stencil $S_2$ centered in $\mathbf{x}_\bhalf = (0.5\Delta x, 0.5\Delta y$). The rectangular black grid represents the space dimension while the blue line means the time direction of Taylor approximation. }
	\label{2D_CAT_1}
\end{figure}

The following multi-index notation will be used:
$$
\bi = (i_1, i_2) \in \mathbb{Z} \times \mathbb{Z},
$$
and
$$
\textbf{0} =(0,0) \quad \bone = (1,1)\quad \bhalf= (1/2, 1/2), \quad \bef = (1,0), \quad \bes = (0,1).
$$
We consider Cartesian meshes with nodes
$$\mathbf{x}_{\bi} = (i_1 \Delta x, i_2 \Delta y).$$
Using this notation, CAT$2P$ methods  can be extended as follows:
\begin{equation}
\label{2Dscheme}
U_\bi^{n+1}=U_\bi^n + \frac{\Delta t}{\Delta x}\left[ {F}_{\bi -  \frac{1}{2}\bef}^P- {F}_{\bi +  \frac{1}{2}\bef}^P+ 
\widetilde{S}^P_{1,\bi} \right] +
\frac{\Delta t}{\Delta y}\left[ {G}_{\bi -  \frac{1}{2}\bes}^P- {G}_{\bi +  \frac{1}{2}\bes}^P
+ \widetilde{S}^P_{2,\bi}  \right] 
\end{equation}
where the numerical fluxes ${F}_{\bi +  \frac{1}{2}\bef}^P$,  ${G}_{\bi +  \frac{1}{2}\bes}^P$ and the source terms 
$\widetilde{S}^P_{j,\bi} $, $j =1,2$ are computed using the values of the numerical solution $U_\bi^n$ in the $P^2-$point stencil centered at $\mathbf{x}_{\bi+\bhalf} = ((i_1+\ha)\Delta x, (i_2+\ha)\Delta y)$ 
$$
\mathcal{S}^P_{\bi + \bhalf} = \{ \mathbf{x}_{\bi + \bj}, \quad \bj \in \mathfrak{I}_P \}, 
$$
where 
$$
\mathfrak{I}_P =\{ \bj = (j_1, j_2) \in \mathbb{Z} \times \mathbb{Z}, \quad -P +1 \le j_k \le P, \quad k = 1,2 \}.
$$
See Figure \ref{2D_CAT_1} for an example. 

For the sake of simplicity let us introduce index sets which will be used in the approximation of integrals of the source term
$$
\mathfrak{M}_P =\{ \bj \in\mathfrak{I}_P \quad {\rm such \; that} \quad -P +2 \le j_1 \le P \}
$$
and
$$
\mathfrak{N}_P =\{ \bj \in\mathfrak{I}_P\quad {\rm such \;that} \quad -P +2 \le j_2 \le P \}.
$$
The following algorithm will be used to compute the numerical fluxes ${F}_{\bi +  \frac{1}{2}\bef}^P$,  ${G}_{\bi +  \frac{1}{2}\bes}^P$ and the source terms 
$\widetilde{S}^P_{j,\bi} $, $j =1,2$ of 2D CAT$2P$ method:
\begin{itemize}
    \item Define
          \begin{eqnarray*}
          & &F_{\bi,\bj}^{(0)} := F(U_{\bi+\bj}^n), \quad \bj \in \mathfrak{I}_P; \\
          & &G_{\bi,\bj}^{(0)} :=  G(U_{\bi+\bj}^n), \quad \bj \in \mathfrak{I}_P; \\
          & &  I_{\bi,\bj-\bef,\bj}^{(0)} := \Delta x\sum_{q=-P+1}^{P}a_{P,q}^{i_1,j_1}S_1(U^n_{\bi+q\bef})H_x(x_{\bi+q\bef}), \quad \bj\in\mathfrak{M}_P;\\
           & &I_{\bi,(-P+1,j_2),(-P+1,j_2)}^{(0)} = 0, \quad j_2 = -P+1, \dots, P; \\ 
           & &I_{\bi,(-P+1,j_2),\bj}^{(0)} = \sum_{s = -P+2}^{j_1}I_{\bi,(s-1,j_2),(s,j_2)}^{(0)} \quad \bj\in \mathfrak{M}_P; \\
          & &J_{\bi,\bj-\bes,\bj}^{(0)} := \Delta y\sum_{q=-P+1}^{P}a_{P,q}^{i_2,j_2}S_2(U^n_{\bi+q\bes})H_y(x_{\bi+q\bes}), \quad \bj\in\mathfrak{N}_P; \\
           & &J_{\bi,(j_1,-P+1),(j_1,-P+1)}^{(0)} = 0, \quad j_1= -P+1, \dots, P; \\ 
          & & J_{\bi,(j_1,-P+1),\bj}^{(0)} = \sum_{s = -P+2}^{j_2}J_{\bi,(j_1,s-1),(j_1,s)}^{(0)} \quad \bj\in \mathfrak{N}_P;
          \end{eqnarray*}
    \item For $k = 1,\ldots,2P-1:$
          \begin{itemize}
              \item Compute for all $\bj\in\mathfrak{I}_P$
                    \begin{eqnarray*}
                    U_{\bi,\bj}^{(k)}&=&  -A_{P}^{1,j_1}\Bigl(F_{\bi,(\bcdot,j_2)}^{(k-1)},\Delta x\Bigr) +  A_{P}^{1,j_1}\Bigl(I_{\bi,(-P+1,j_2),(\bcdot,j_2)}^{(k-1)},\Delta x\Bigr) \\
                    & & - A_{P}^{1,j_2}\Bigl(G_{\bi,(j_1,\bcdot)}^{(k-1)},\Delta y\Bigr)  +  A_{P}^{1,j_2}\Bigl(J_{\bi,(j_1,-P+1),(j_1,\bcdot)}^{(k-1)},\Delta y\Bigr).
                    \end{eqnarray*}
                    
               \item Compute for all $\bj\in\mathfrak{I}_P$ and for all $r=-P+1,\ldots,P$
                     $$ U_{\bi,\bj}^{k,n+r} = U_{\bi+\bj}^n + \sum_{m = 1}^{k}\frac{(\Delta t)^m}{m!}U_{\bi,\bj}^{(m)}.$$
               \item Compute for all $\bj\in\mathfrak{I}_P$ and for all $r=-P+1,\ldots,P$
                     $$ F_{\bi,\bj}^{k,n+r} = F(U_{\bi,\bj}^{k,n+r}) \quad {\rm{and}}\quad G_{\bi,\bj}^{k,n+r} = G(U_{\bi,\bj}^{k,n+r}).$$
               \item Compute for all $\bj\in\mathfrak{M}_P$ and for all $r=-P+1,\ldots,P$
                     \begin{eqnarray*}
                         I_{\bi,\bj-\bef,\bj}^{k,n+r} &=& \Delta x\sum_{q=-P+1}^{P}a_{P,q}^{i_1, j_1}S_1(U^n_{\bi+q\bef})H_x(x_{\bi+q\bef}); \\
                         I_{\bi,\bj-\bef,\bj}^{(k)} &=& A^{k,0}_{P}\Bigl(I_{\bi,\bj-\bef,\bj}^{k,\bcdot},\Delta t\Bigr).
                     \end{eqnarray*}            
                \item Compute for all $\bj\in\mathfrak{N}_P$ and for all $r=-P+1,\ldots,P$
                      \begin{eqnarray*}
                         J_{\bi,\bj-\bes,\bj}^{k,n+r} &=& \Delta y \sum_{q=-P+1}^{P}a_{P,q}^{i_2,j_2}S_2(U^n_{\bi+q\bes})H_y(x_{\bi+q\bes}); \\
                         J_{\bi,\bj-\bes,\bj}^{(k)} &=& A^{k,0}_{P}\Bigl(J_{\bi,\bj-\bes,\bj}^{k,\bcdot},\Delta t\Bigr).
                     \end{eqnarray*} 
                \item Compute
                      \begin{eqnarray*}
                          & &F_{\bi,\bj}^{(k)} = A^{k,0}_{P}\Bigl(F_{\bi,\bj}^{k,\bcdot},\Delta t\Bigr), \quad \bj\in\mathfrak{I}_P;\\
                          & &G_{\bi,\bj}^{(k)} = A^{k,0}_{P}\Bigl(G_{\bi,\bj}^{k,\bcdot},\Delta t\Bigr), \quad \bj\in\mathfrak{I}_P;\\
                      & & I_{\bi,(-P+1,j_2),(-P+1,j_2)}^{(k)} = 0, \quad j_2 = -P+1, \dots, P; \\ 
                          & &I_{\bi,(-P+1,j_2),\bj}^{(k)} = \sum_{s = -P+2}^{j_1}I_{\bi,(s-1,j_2),(s,j_2)}^{(k)} \quad \bj\in \mathfrak{M}_P; \\
                          & &J_{\bi,(j_1,-P+1),(j_1,-P+1)}^{(k)} = 0, \quad j_1 = -P + 1, \dots,P; \\ 
                          & &J_{\bi,(j_1,-P+1),\bj}^{(k)} = \sum_{s = -P+2}^{j_2}I_{\bi,(j_1,s-1),(j_1,s)}^{(k)} \quad \bj\in \mathfrak{N}_P; 
                          \end{eqnarray*}
          \end{itemize}
          \item Compute 
                \begin{eqnarray}
                     & & F^P_{\bi+\bhalf\bef} = \sum_{k=1}^{2P} \frac{\Delta t^{k-1}}{k!}A^{0,\ha}_{P}\Bigl(F_{\bi,(\bcdot,0)}^{(k-1)},\Delta x\Bigr); \\ 
                     & & G^P_{\bi+\bhalf\bes} = \sum_{k=1}^{2P} \frac{\Delta t^{k-1}}{k!}A^{0,\ha}_{P}\Bigl(G_{\bi,(0,\bcdot)}^{(k-1)},\Delta y\Bigr) 
                \end{eqnarray}
          Once the algorithm has finished, the integrals will have already been computed and can be used to approximate the source term as follows:
          \item For $k=1,\ldots,2P$ define
                 \begin{eqnarray*}
                     \I_{\bi,j_1}^{(k-1)} &=&
                     \begin{cases}
                         I_{\bi-\bef,j_1\bef,(j_1 + 1)\bef}^{(k-1)} \quad {\rm{if}}\quad j_1=-P+1,\ldots,0;\\
                         I_{\bi,(j_1 - 1)\bef,j_1\bef}^{(k-1)} \quad {\rm{for}}\quad j_1=1,\ldots,P.
                         \end{cases} \\
                         \mathcal{J}_{\bi,j_2}^{(k-1)} &=&
                     \begin{cases}
                         J_{\bi-\bes,j_2\bes,(j_2 + 1)\bes}^{(k-1)} \quad {\rm{for}}\quad j_2=-P+1,\ldots,0;\\
                         J_{\bi,(j_2 - 1)\bes, j_2\bes}^{(k-1)} \quad {\rm{if}}\quad j_2=1,\ldots,P.
                         \end{cases}
                \end{eqnarray*}
          \item Compute 
                \begin{eqnarray}
                     & &\widetilde{S}^P_{1,\bi} = \sum_{k=1}^{2P} \frac{\Delta t^{k-1}}{k!}A^{0,\ha}_{P}\Bigl(\I_{\bi,\bcdot}^{(k-1)},\Delta x\Bigr); \\ 
                     & &\widetilde{S}^P_{2,\bi} = \sum_{k=1}^{2P} \frac{\Delta t^{k-1}}{k!}A^{0,\ha}_{P}\Bigl(\mathcal{J}_{\bi,\bcdot}^{(k-1)},\Delta y\Bigr). 
                \end{eqnarray}
    \end{itemize}

In the case of the well-balanced methods WBCAT$2P$ there is an important difference: if the algorithm described in Subsection \ref{ss:WBCAT2Pbl} (adopting the 2D  above notation) had to be used, the first step to update the numerical solution at the point $\mathbf{x}_\bi$ at a time $t_n$ would be to find a solution to the problem
   \begin{equation}\label{Cauchyin2d}
    \left \{
    \begin{array}{l}
    F(U)_x + G(U)_y = S_1(U) H_x + S_2(U) H_y\\
    U(\mathbf{x}_{\bi}) = U_{\bi}^n.
    \end{array}
    \right.
\end{equation}
This problem is obviously  much more difficult to solve, either exactly or numerically,  than \eqref{Cauchyin} since it is about a nonlinear PDE system instead of an ODE system. Moreover in this case there may exist infinitely many stationary solutions satisfying the condition at only one point $\mathbf{x}_{\bi}$: some extra conditions have to be imposed to determine one of them.

Nevertheless, if the stationary solutions to be preserved constitute a $k$-parameter family
$$
U^*(x, y; C_1,\dots, C_k),
$$
with $k < d$, then the numerical strategy described in Subsection \ref{ss_family} can be followed: this strategy will be used in Subsection \ref{ss:2DEuler} to preserve a family of stationary solutions of the 2D Euler system with gravity.

The extension of the adaptive CAT2$P$ and WBCAT2$P$ to \eqref{2Dequ} is similar. Following what has been done for ACAT2$P$ and WBACAT2$P$ in 1D we define $\A,$ the index set used to select the square stencils according with the smoothness of numerical data, is defined as:
\begin{eqnarray} \label{admissibleset_2D}
\A_\bi & = & \{  p \in \{2, \dots, P \} \text{ s.t. } \psi^p_{\bi\pm\ha\bef}\approx 1 \quad  \rm{and}\quad\psi^p_{\bi\pm\ha\bes}\approx 1 \}
\end{eqnarray} where $\psi^p_{\bi+\frac{1}{2}\bef}$, $\psi^p_{\bi+\frac{1}{2}\bes}$ are the smoothness indicators introduced in Subsection \ref{ss:smoothness} computed direction by direction. Then the 2D adaptive CAT2P are so defined:
     \begin{equation}\label{numfluxacat_f_2D}
    {F}_{\bi;\bi \pm  \ha\bef}^{\mathcal{A}_\bi} = \begin{cases}
    {F}^{*}_{\bi;\bi \pm\ha\bef} & \text{if $\mathcal{A}_\bi = \emptyset$;}\\
    {F}^{p_s}_{\bi \pm  \ha\bef} & \text{otherwise, where $p_s = \max(\mathcal{A}_\bi);$ }
    \end{cases}
    \end{equation}
    \begin{equation}\label{numfluxacat_g_2D}
    {G}^{\A_\bi}_{\bi;\bi \pm  \ha\bes} = \begin{cases}
    {G}^{*}_{\bi;\bi \pm  \ha\bes} & \text{if $\mathcal{A}_\bi = \emptyset$;}\\
    {G}^{p_s}_{\bi\pm  \ha\bes} & \text{otherwise, where $p_s = \max(\mathcal{A}_\bi);$}
    \end{cases}
    \end{equation}
    and
      \begin{equation}\label{numsourceacat_2D}
    \widetilde{S}^{\A_\bi}_{j\bi} = \begin{cases}
    \widetilde{S}^*_{j,\bi} & \text{if $\mathcal{A}_\bi = \emptyset$;}\\
    \widetilde{S}^{p_s}_{j,\bi}  & \text{otherwise, where $p_s = \max(\mathcal{A}_\bi).$} 
    \end{cases}
    \end{equation}
    $F^*_{\bi +  \ha\bef},$ $G^*_{\bi +  \ha\bes},$ and $\widetilde{S}_{j,\bi}^*$ are the ACAT2 numerical flux and source terms given by 1D \eqref{ACAT2_flux}, (\ref{ACAT2_source}); $F^{p_s}_{\bi +  \ha\bef},$ $G^{p_s}_{\bi +  \ha\bes}$ and $\widetilde{S}_{j,\bi}^{p_s}$ are the ACAT2$p_s$ numerical fluxes and source terms defined in 1D (\ref{flux_2P}), (\ref{source}).
    
    Using this notation, ACAT$2P$ methods  can be extended as follows:
\begin{equation*}
U_\bi^{n+1}=U_\bi^n + \frac{\Delta t}{\Delta x}\left[ {F}_{\bi -  \frac{1}{2}\bef}^{\A_\bi}- {F}_{\bi +  \frac{1}{2}\bef}^{\A_\bi}+ 
\widetilde{S}^{\A_\bi}_{1,\bi} \right] +
\frac{\Delta t}{\Delta y}\left[ {G}_{\bi -  \frac{1}{2}\bes}^{\A_\bi}- {G}_{\bi +  \frac{1}{2}\bes}^{\A_\bi} + \widetilde{S}^{\A_\bi}_{2,\bi}  \right] 
\end{equation*}
    
    In a similar way, the 2D adaptive well-balanced WBACAT2P writes as follows:
    \begin{equation*}
U_\bi^{n+1}=U_\bi^n + \frac{\Delta t}{\Delta x}\left[ {F}_{\bi;\bi -  \frac{1}{2}\bef}^{\A_\bi}- {F}_{\bi;\bi +  \frac{1}{2}\bef}^{\A_\bi} + 
\widetilde{S}^{\A_\bi}_{1,\bi} \right] +
\frac{\Delta t}{\Delta y}\left[ {G}_{\bi;\bi -  \frac{1}{2}\bes}^{\A_\bi}- {G}_{\bi;\bi +  \frac{1}{2}\bes}^{\A_\bi} + \widetilde{S}^{\A_\bi}_{2,\bi}  \right],
\end{equation*}
where    
     \begin{equation}\label{numfluxacat_f_2D_wb}
    {F}_{\bi;\bi \pm  \ha\bef}^{\A_\bi} = \begin{cases}
    {F}^{*}_{\bi;\bi \pm  \ha\bef} & \text{if $\mathcal{A}_\bi = \emptyset$;}\\
    {F}^{p_s}_{\bi;\bi \pm  \ha\bef} & \text{otherwise, where $p_s = \max(\mathcal{A}_\bi);$}
    \end{cases}
    \end{equation}
    \begin{equation}\label{numfluxacat_g_2D_wb}
    {G}^A_{\bi;\bi \pm  \ha\bes} = \begin{cases}
    {G}^{*}_{\bi;\bi \pm  \ha\bes} & \text{if $\mathcal{A}_\bi = \emptyset$;}\\
    {G}^{p_s}_{\bi;\bi \pm  \ha\bes} & \text{otherwise, where $p_s = \max(\mathcal{A}_\bi);$}
    \end{cases}
    \end{equation}
    and
      \begin{equation}\label{numsourceacat_2D_wb}
    \widetilde{S}^{\A_\bi}_{j\bi} = \begin{cases}
    \widetilde{S}^*_{j,\bi} & \text{if $\mathcal{A}_\bi = \emptyset$;}\\
    \widetilde{S}^{p_s}_{j,\bi}  & \text{otherwise, where $p_s = \max(\mathcal{A}_\bi).$} 
    \end{cases}
    \end{equation}
    Here, $F^*_{\bi;\bi +  \ha\bef},$ $G^*_{\bi;\bi +  \ha\bes},$ and $\widetilde{S}_{j,\bi}^*$ are the ACAT2 numerical flux and source terms given by 1D \eqref{WBACAT2_flux},\eqref{WBACAT2_source}; $F^{p_s}_{\bi;\bi +  \ha\bef},$ $G^{p_s}_{\bi;\bi +  \ha\bes}$ and $\widetilde{S}_{j,\bi}^{p_s}$ are the ACAT2$p_s$ numerical fluxes and source terms defined in 1D \eqref{WBCAT2P_flux-}, \eqref{WBCAT2P_flux+}, \eqref{WBCAT2P_source}.

\section{Numerical results}\label{nUm}
In this section we apply ACAT$2P$ and WBACAT$2P$, $P = 1,2$ methods to several 1D and 2D problems: the 1D linear transport equation and Burgers equation with source term, the 1D shallow water system, and 2D Euler equations with a gravitational potential.
The Minmod flux limiter \cite{Roesuperbee} is used in ACAT2  and the smoothness indicators \eqref{indicadores_locales} are used for ACAT4: no loss of precision for first order critical points has been observed in any of the test problems considered here due to the use of $\psi^2_{i+\ha}$. Fornberg's algorithm \cite{Fornberg1} is used to compute the coefficients of the numerical differentiation formulas.

    \subsection{Linear Equation}
    We consider the linear scalar balance law
    \begin{equation}
        \label{bal_scalar}
        u_t+ u_x=u,
    \end{equation}
that has the form \eqref{bal_sis} with  $H(x) = x$. The analytic solution of the initial value problem  with condition
$$
u(x,0) = u_0(x)
$$
is given by:
    \begin{equation}
        \label{Exact_lin}
        u(x,t) = u_0(x-t)e^{t}.
    \end{equation} 
    The stationary solutions solve the ODE
    \begin{equation*}
        u_x = u.
    \end{equation*}
    Hence, the set of stationary solutions is 
    $$ u^*(x) = Ce^{x}, \quad C\in\mathbb{R}.$$

    \subsubsection{Order test}\label{Order test}
    Following  \cite{PP2021}, we consider (\ref{bal_scalar}) with initial condition
    \begin{equation}
        \label{condsuave_1}
        u_0(x)=
        \begin{cases}
           0 \quad\quad \;\;\;\mathrm{if}\quad  x<0;\\
           p(x) \quad\;\;\mathrm{if}\quad 0\le x \le1; \\
           1 \quad \quad\;\;\;\mathrm{if} \quad x>1; 
        \end{cases}
    \end{equation} 
    \begin{figure}[!ht]
    	\centering	
    	\hspace{-1.2cm}
    	\includegraphics[scale = 0.5]{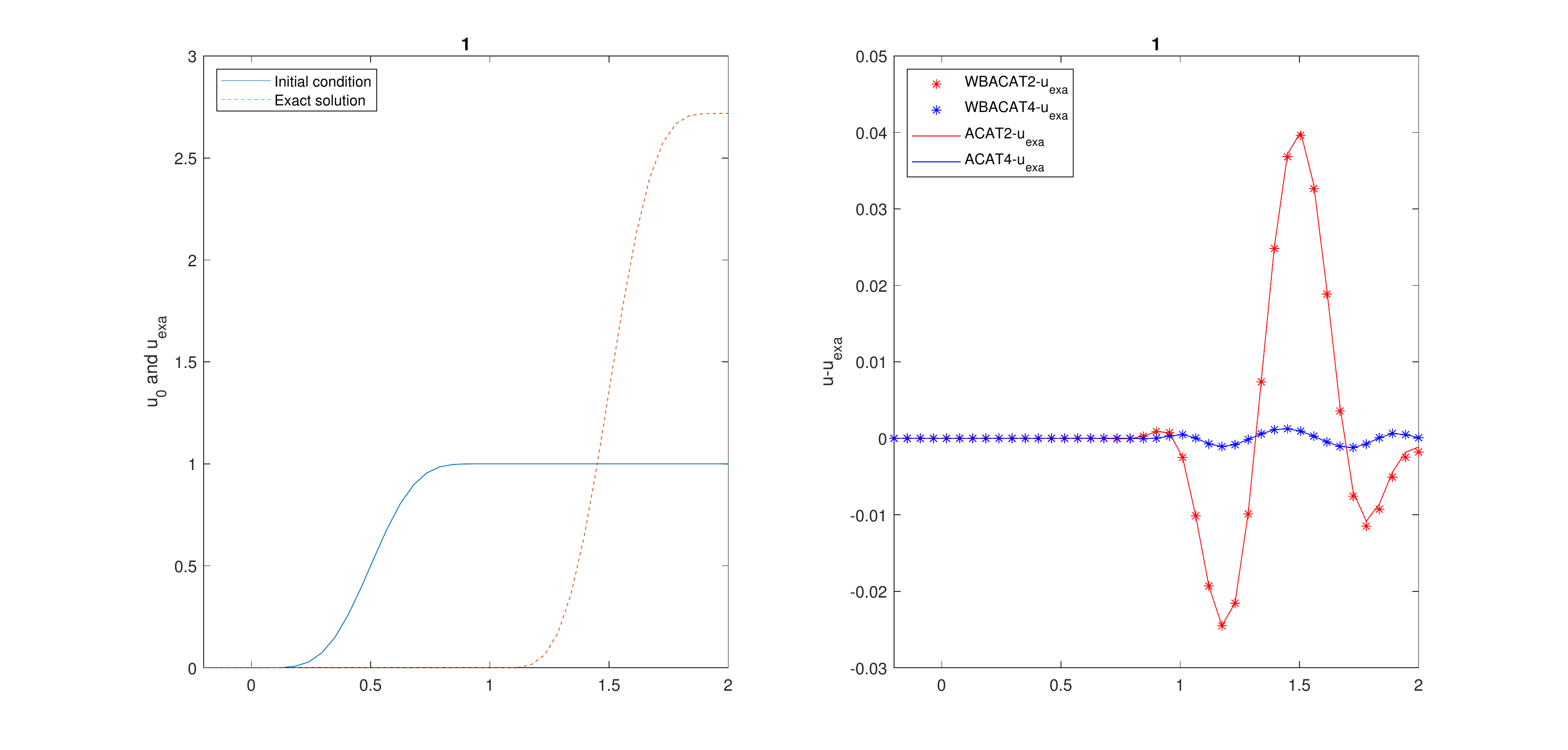}
    	\vspace{-0.8 cm}
    	\caption{Test \ref{Order test}. Initial condition and exact solution (left); difference between the numerical and the exact solutions computed with CAT2, CAT4, WBCAT2 and WBCAT4 at $t = 1$ using a mesh of 41 points (right) on the interval $[-0.2,2]$ and CFL$=0.9.$}
    	\label{Test_1_1}
    \end{figure}
    where $p(x)$ is the polynomial that satisfies  $p(0)=0$, $p(1)=1$, $p^k(0)=p^k(1)=0$, $k=1,\ldots, 5$:
    $$ p(x) = x^6\Bigl(\sum_{k=0}^{5}(-1)^k\binom{5+k}{k}(x-1)^k\Bigr)$$ 
    (see Figure \ref{Test_1_1}). 
    The methods  ACAT2, ACAT4, WBACAT2, WBACAT4 have been applied to  (\ref{bal_scalar}) with initial condition (\ref{condsuave_1})  in the spatial interval $[-0.2,2],$ with CFL$=0.9.$ Dirichlet boundary conditions are considered to the left and free boundary conditions to the right based on the use of ghost cells. Figure \ref{Test_1_1} shows the numerical solutions obtained at time $t=1$ on the interval $[-0.2,2]$ using $40$ mesh points.
    
    Tables \ref{table_1_1}-\ref{table_1_2} show the $L^1$-errors and the empirical order of convergence corresponding to the standard and Adaptive CAT2P and WBCAT2P with $p=1,2.$ As it can be seen, all the schemes keep the expected order and the errors corresponding to methods of the same order are almost identical. In the first case, the smoothness indicators of ACAT2P and WBACAT2P have been fixed to 1, hence, the Adaptive CAT2P coincides exactly with the standard CAT2P scheme; regarding the second case, no restrictions are imposed on the smoothness indicators, requiring an increase in the number of points to capture the theoretical order. No further restrictions are required for the time step.  
    \begin{table}[htbp]
        \begin{center}
            \begin{tabular}{|c|c|c|c|c|c|c|c|c|}
                \hline                        &\multicolumn{2}{c|}{CAT2}&\multicolumn{2}{c|}{WBCAT2}&  \multicolumn{2}{c|}{CAT4} &  \multicolumn{2}{c|}{WBCAT4}  \\
                Points   &  Order & Error      & Order & Error        & Order  & Error     & Order& Error  \\ \hline
                6       & -    &  1.79E-1     & -     & 1.80E-1      & -      & 6.03E-2   & -    &  6.03E-2   \\   
                \hline 
                11      & 1.57 &  6.00E-2     & 1.56  & 6.05E-2      & 2.86   & 8.31E-3   & 2.86 &  8.31E-3\\   
                \hline  
                21      & 1.91 &  1.60E-2     & 1.92  & 1.60E-2      & 3.58   & 6.94E-4   & 3.58 &  6.94E-4 \\   
                \hline  
                41      & 2.02 &  3.93E-3     & 2.02  & 3.94E-3      & 3.88   & 4.69E-5   & 3.89 &  4.69E-5 \\   
                \hline  
                81      & 2.03 &  9.63E-4     & 2.02  & 9.66E-4      & 4.01   & 2.90E-6   & 4.00 &  2.90E-7 \\   
                \hline  
            \end{tabular}
            \vspace{2mm}
            \caption{Test \ref{Order test}: Errors in $L^1$ norm and convergence rates for CAT2, CAT4, WBCAT2 and WBCAT4 at time $t=1.$}
            \label{table_1_1}
            \end{center}
    \end{table}
    \begin{table}[htbp]
        \begin{center}
            \begin{tabular}{|c|c|c|c|c|c|c|c|c|}
                \hline                        &\multicolumn{2}{c|}{ACAT2}&\multicolumn{2}{c|}{WBACAT2}&  \multicolumn{2}{c|}{ACAT4} &  \multicolumn{2}{c|}{WBACAT4}  \\
                Points   &  Order & Error      & Order & Error        & Order  & Error     & Order& Error  \\ \hline
                16       & -    &  2.00E-1     & -     & 1.99E-1      & -      & 1.57E-2   & -    &  1.61E-2   \\   
                \hline 
                31      & 1.51 &  7.01E-2      & 1.50  & 7.04E-2      & 3.01   & 1.95E-3   & 2.99 &  2.02E-3\\   
                \hline  
                61      & 1.85 &  1,94E-2      & 1.86  & 1.94E-2      & 3.62   & 1.58E-4   & 3.60 &  1.67E-4 \\   
                \hline  
                121      & 2.19 &  4.25E-3     & 2.19  & 4.25E-3      & 8.10   & 5.71E-7   & 8.20 &  5.71E-7 \\   
                \hline  
                241      & 2.00 &  1.05E-3     & 2.00  & 1.06E-3      & 3.99   & 3.57E-8   & 4.00 &  3.56E-8 \\   
                \hline
                481      & 2.00 &  2.64E-4     & 2.00  & 2.64E-4      & 4.00   & 2.22E-9   & 4.00 &  2.22E-9 \\   
                \hline
            \end{tabular}
            \vspace{2mm}
            \caption{Test \ref{Order test}: Errors in $L^1$ norm and convergence rates for ACAT2, ACAT4, WBACAT2 and WBACAT4 at time $t=1.$}
            \label{table_1_2}
            \end{center}
    \end{table}

    \subsection{Burgers Equation}
    In this section we consider the scalar Burgers equation with  source term:
    \begin{equation}
        \label{bur_scal}
        u_t + \Bigl(\frac{u^2}{2}\Bigr)_x = u^2H_x(x).
    \end{equation}
    The stationary solutions solve now the ODE
    $$u_x = uH_x(x) $$ 
    whose general solution is 
    $$u^*(x) = Ce^{H(x)}, \quad C\in\mathbb{R}.$$
    
    \subsubsection{Preservation of a stationary solution}\label{Burg_1}
    As first test let us check numerically the well-balanced property. For this reason, we consider an oscillatory  
    \[
        H(x) = x + 0.1\sin(10x)
    \] 
    and consequently the stationary solutions become $u^*(x) = C e^{H(x)}.$
    We solve (\ref{bur_scal}) with initial condition 
    $$u_0(x) = e^{H(x)}$$
    in the interval $[-1,1]$ using $100$ mesh points and CFL$=0.9.$ As boundary conditions, the stationary solution is imposed at ghost cells. 
    
    Table \ref{table_2_1_1} shows the $L^1$-errors and the empirical order of convergence corresponding to ACAT$2P$, and the error of WBACAT$2P$ methods, with $P=1,2$ at time $t=8.$ While the non well-balanced schemes give accurate solutions according to their order, the well-balanced methods capture the stationary solution with machine precision.
    \begin{table}[htbp]
        \begin{center}
            \begin{tabular}{|c|c|c|c|c|c|c|}
                \hline                        &\multicolumn{2}{c|}{ACAT2}&WBACAT2 &\multicolumn{2}{c|}{ACAT4}&WBACAT4  \\
                Points   &  Order & Error      & Error    &  Order & Error      & Error\\ \hline
                100      & -    &  1.82E-3     & 2.66E-17 & -    &  1.93E-5     & 3.99E-17\\   
                \hline 
                200      & 2.05 &  4.32E-4     & 3.11E-17 & 4.09 &  1.13E-6     & 6.21E-17\\   
                \hline  
                400      & 2.03 &  1.07E-4     & 2.44E-17 & 4.05 &  6.84E-8     & 2.22E-17\\   
                \hline  
                800      & 2.01 &  2.63E-5     & 2.78E-17 & 4.03 &  4.19E-9     & 3.55E-17\\   
                \hline  
                1600     & 2.00 &  6.55E-6     & 1,99E-17 & 4.01 &  2.59E-10    & 5.21E-17\\   
                \hline  
            \end{tabular}
            \vspace{2mm}
            \caption{Test \ref{Burg_1}. Test of well-balance property. Errors in $L^1$ norm and convergence rates at time $t=8$ for ACAT2-4. The errors for  WBACAT2-4 are due to round-off.}
            \label{table_2_1_1}
            \end{center}
    \end{table}
    
    Checked numerically the well-balanced property for the Burgers equation with source we will now focus on some experiment to test the accuracy of the numerical solutions where initial conditions are: a perturbation of the stationary solution; and a smooth perturbation of Test \ref{Order test}.

    \subsubsection{Perturbation of a stationary solution}\label{Burg_4}
\begin{figure}[!ht]
    	\centering	
    	\hspace{-1.2cm}
    	\includegraphics[scale = 0.5]{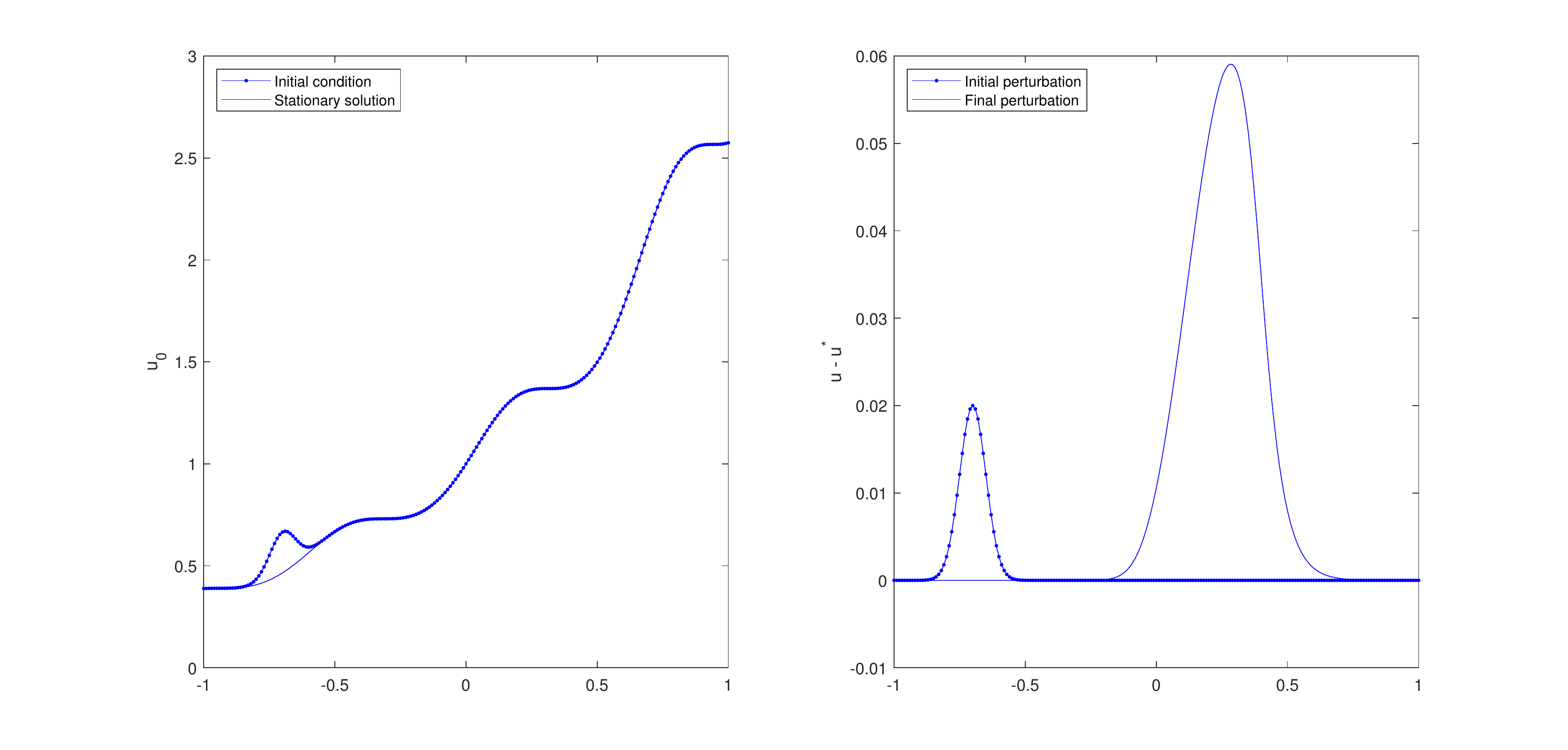}
    	\caption{Test \ref{Burg_4}. Initial condition and stationary solutions (left). Difference between reference and stationary solution at initial and final time (right). The perturbation of the initial condition (left) is amplified by 10 times in order to see clearly the perturbation. The reference solution is computed with WBACAT4 using 1000 mesh points and CFL$=0.9$ at time $t=1.2.$}
    	\label{Test_2_4_1}
    \end{figure}

    Let us consider (\ref{bur_scal}) with oscillatory $H$ given by 
    \[
        H(x) = x + 0.1\sin(10x)
    \]
    and initial condition
       \begin{equation}
       \label{u0_burg_4}
       u_0(x) = e^{H(x)} + 0.02e^{-200(x+0.7)^2},
    \end{equation}
    that is a small smooth perturbation of the stationary solution $u^*(x) = e^{H(x)}$: see Figure \ref{Test_2_4_1}. 
    We solve the problem in the interval $[-1,1]$ using 200 mesh points and CFL$=0.9.$ As boundary conditions the stationary solution is imposed at left ghost point and free boundary at right.
    
    \begin{figure}[!ht]
    	\centering	
    	\hspace{-1.2cm}
    	\includegraphics[scale = 0.5]{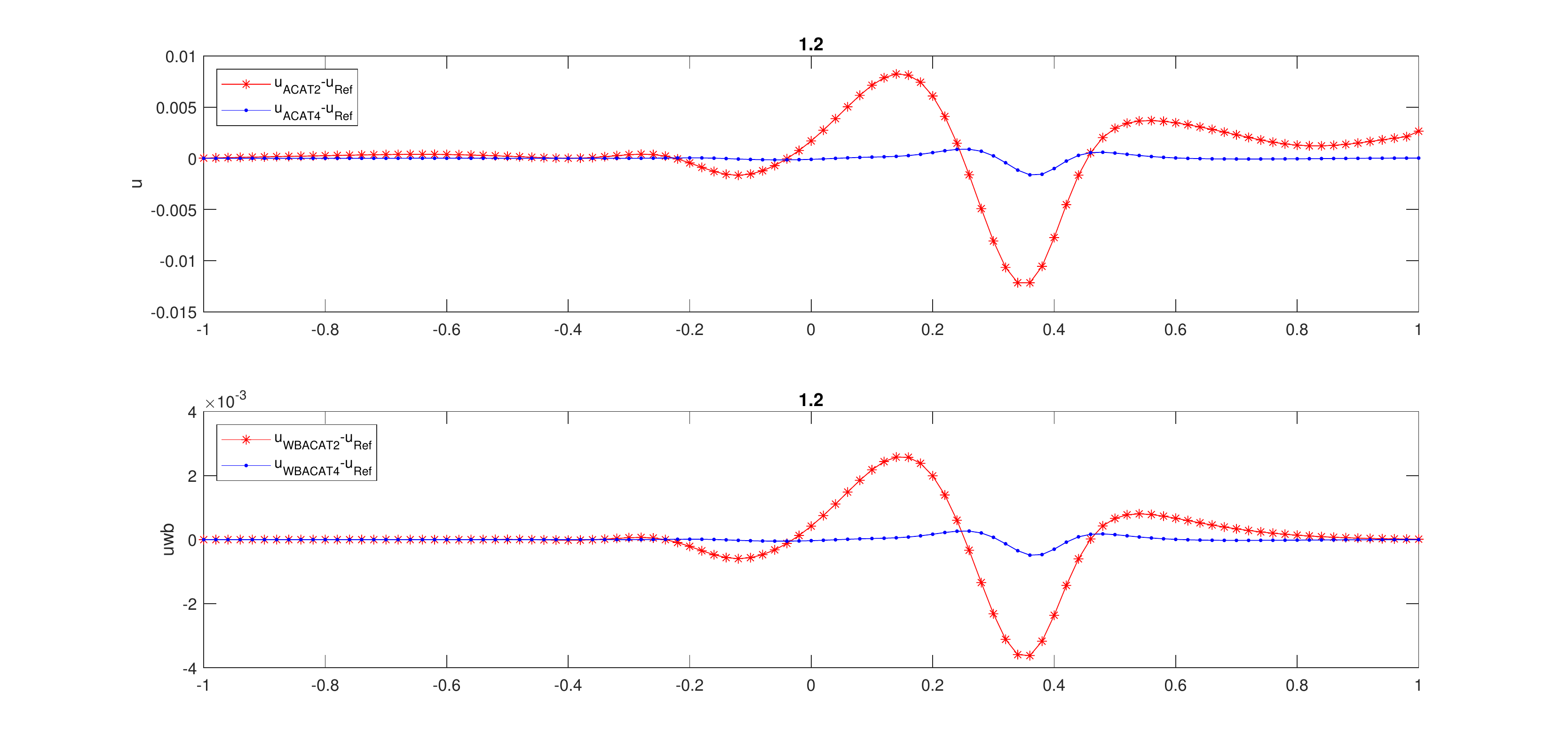}
    	\vspace{-0.8 cm}
    	\caption{Test \ref{Burg_4}. Difference between numerical solutions computed with ACAT$2P$ (top) and WBCAT$2P,$ (bottom) $P=1,2,$ and reference solution at $t = 1.2$ using a mesh of 200 points and CFL$=0.9.$ For the reference solution the WBACAT4 is adopted with a mesh of 1000 points.}
    	\label{Test_2_4_2}
    \end{figure}
    
   Figure \ref{Test_2_4_2} shows that, all the schemes are able to evolve the perturbation in according with the order. Nevertheless, the well-balanced methods, WBACAT2 and WBACAT4, are able to capture more precisely the evolution of the perturbation with a smaller error than the relative non well-balanced schemes. \\ For the errors in $L^1-$norm and convergence rates, we adopt as initial condition
    \begin{equation*}
       u_0(x) = e^{H(x)} + 0.2e^{-200(x+0.7)^2},
    \end{equation*} in other word, a bigger perturbation is considered reducing the final time to $t = 0.2$.\\ Table \ref{table_2_2} shows that the non well-balanced methods introduce a bigger error in comparison with the well-balanced approach but an increasing of points is necessary to achieve the theoretical order. This phenomenon is partly attributable to reconstruction partly to smoothness indicators, because they fail to detect the theoretical regularity. 
    
    \begin{table}[htbp]
        \begin{center}
            \begin{tabular}{|c|c|c|c|c|c|c|c|c|}
                \hline                        &\multicolumn{2}{c|}{ACAT2}&\multicolumn{2}{c|}{WBACAT2}&  \multicolumn{2}{c|}{ACAT4} &  \multicolumn{2}{c|}{WBACAT4}  \\
                Points   &  Order & Error      & Order & Error        & Order  & Error     & Order& Error  \\ \hline
                80       & -    &  6.44E-2     & -     & 1.07E-2      & -      & 1.53E-3   & -    &  1.17E-4\\   
                \hline 
                160      & 1.15 &  2.88E-2     & 1.74  & 3.21E-3      & 2.11   & 3.10E-4   & 2.03 &  2.86E-5\\   
                \hline  
                320      & 1.45 &  1.05E-2     & 1.89  & 8.61E-4      & 4.28   & 1.59E-5   & 4.45 &  1.31E-6 \\   
                \hline  
                640      & 1.62 &  3.42E-3     & 1.96  & 2.21E-4      & 4.14   & 8.99E-7   & 3.85 &  8.99E-8 \\   
                \hline  
                1280     & 1.86 &  9.41E-4     & 1.99  & 5.58E-5      & 4.01   & 5.57E-8   & 3.96 &  5.76E-9 \\   
                \hline
            \end{tabular}
            \vspace{2mm}
            \caption{Test \ref{Burg_4}: Errors in $L^1$ norm and convergence rates for ACAT2, ACAT4, WBACAT2 and WBACAT4 at time $t=0.2$ and CFL$=0.9.$
              }
            \label{table_2_2}
            \end{center}
    \end{table}
    
    This experiment is very interesting because, on the one hand it shows that the  well-balanced methods work well near the stationary solution, on the other hand, it happens that a small perturbation of the initial state, although rather smooth, may result in a loss of accuracy in the adaptive order reconstruction.

    \subsubsection{Preservation of a stationary solution with oscillatory $H$}\label{Burg_3_osc}
\begin{figure}[!ht]
    	\centering	
    	\hspace{-1.2cm}
    	\includegraphics[scale = 0.5]{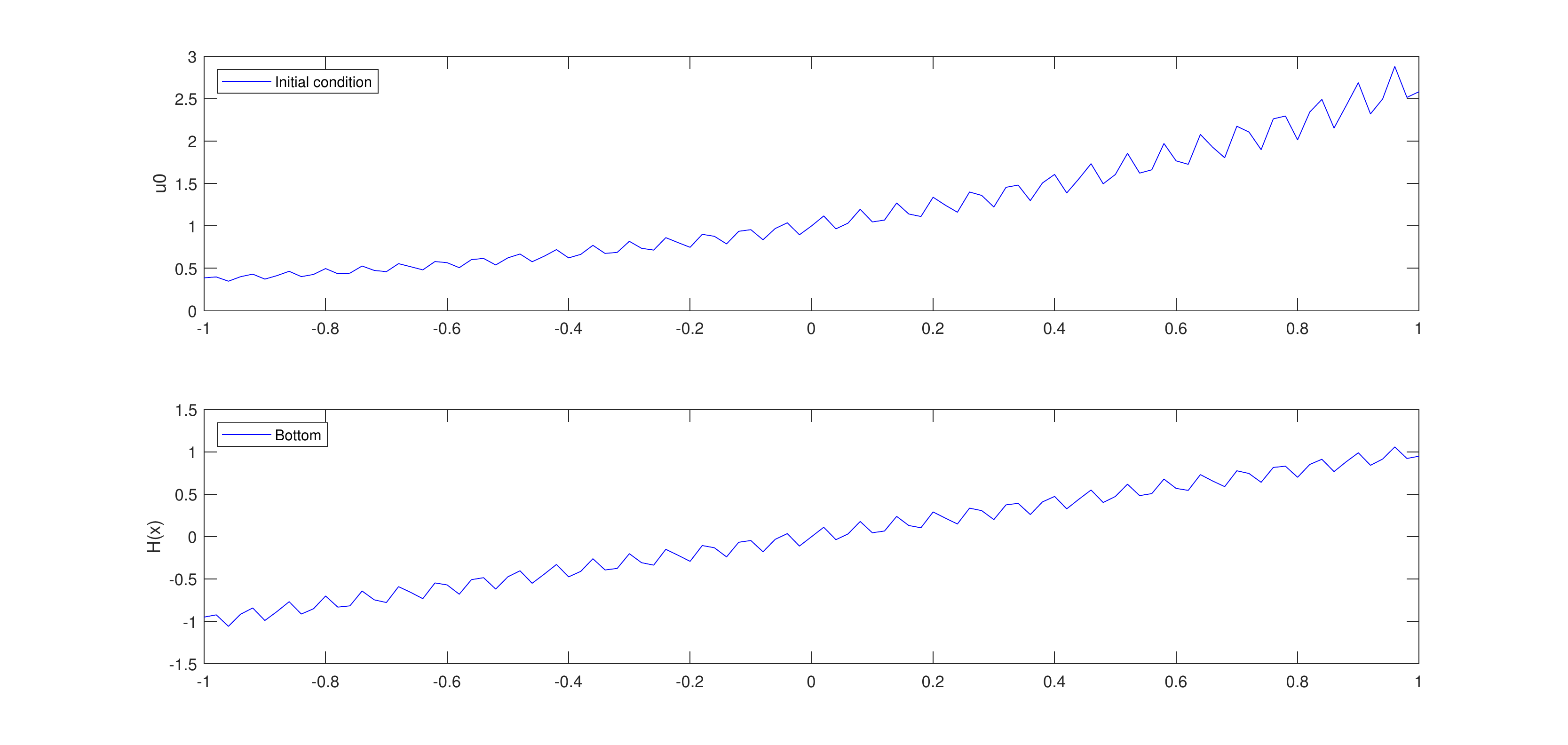}
    	\vspace{-0.8 cm}
    	\caption{Test \ref{Burg_3_osc}. Initial condition (top) and $H$ (down).}
    	\label{Test_2_3_1}
    \end{figure}

     Following  \cite{PP2021}, we consider (\ref{bur_scal}) with 
    \begin{equation}
        \label{H_osc}
        H(x) = x + \frac{1}{10}\sin(100x),
    \end{equation}
    and we take as initial condition the stationary solution
    $$u^*(x) = e^{H(x)}$$ 
    (see Figure \ref{Test_2_3_1}). 
    We solve the problem in the interval $[-1,1]$ using 100 mesh points and $CFL =0.9.$ With this choice of mesh points, the period of the oscillations of $H$ is close to $\Delta x.$ As boundary conditions the stationary solution is imposed again at ghost points.
    
    \begin{figure}[!ht]
    	\centering	
    	\hspace{-1.2cm}
    	\includegraphics[scale = 0.5]{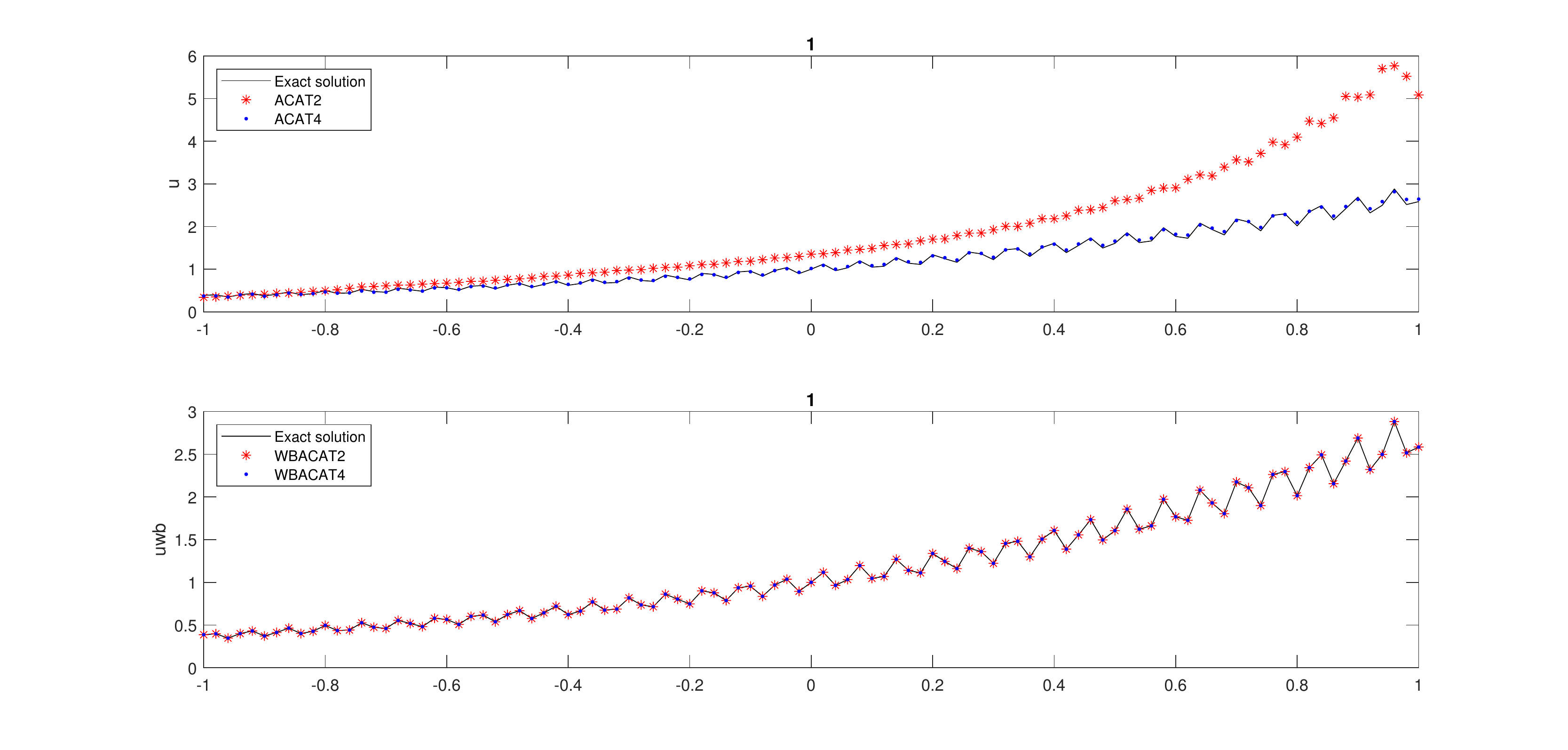}
    	\vspace{-0.8 cm}
    	\caption{Test \ref{Burg_3_osc}. Exact and  numerical stationary solutions computed with ACAT$2P$ and WBCAT$2P$ at $t = 1$ using a mesh of 100 points and CFL$=0.9$: non well-ballanced (top) and well-balanced (bottom).}
    	\label{Test_2_3_2}
    \end{figure}
    \begin{figure}[!ht]
    	\centering	
    	\hspace{-1.2cm}
    	\includegraphics[scale = 0.5]{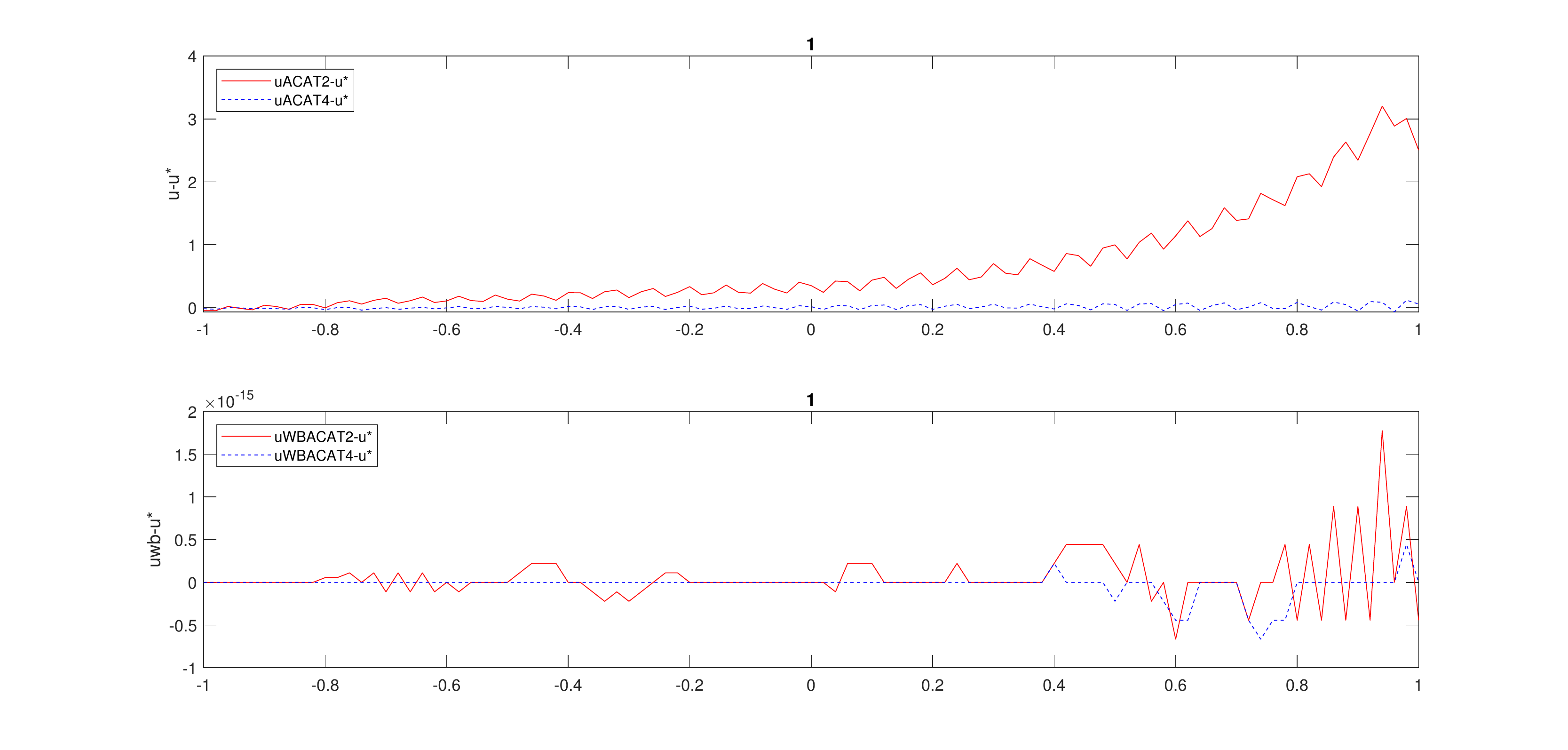}
    	\vspace{-0.8 cm}
    	\caption{Test \ref{Burg_3_osc}. Top: differences between the exact and the numerical stationary solutions computed with ACAT$2P$, $P = 1, 2$, at time $t=1$ using 100 mesh points and $CFL =0.9$. Bottom: differences between the exact and the numerical solutions computed with WBACAT$2P$, $P = 1, 2$, at time $t=1$ using 100 mesh points and $CFL =0.9$.  }
    	\label{Test_2_3_3}
    \end{figure}
    Figure \ref{Test_2_3_2} shows that, while WBACAT2 and WBACAT4  capture the stationary solution with machine precision, this is not the case for ACAT2 and ACAT4. Figure \ref{Test_2_3_3} displays the differences between the numerical solution and the stationary solution obtained at time $t=1$ (top). In this case, the results provided by WBACAT2 and WBACAT4 (bottom) are very similar and they are able to capture the machine precision; while, the results provided by ACAT2 and ACAT4 (top) show that the non well-balanced methods are not able to detect the stationary solution with high precision even if a 4 order method is applied.
    
    As a final check we consider the behaviour of the  methods in the case of  an initial condition of class $\mathcal{C}^5$ which is far from the stationary solution.

    \subsubsection{Order Test} \label{Burg_5}
\begin{figure}[!ht]
    	\centering	
    	\hspace{-1.2cm}
    	\includegraphics[scale = 0.5]{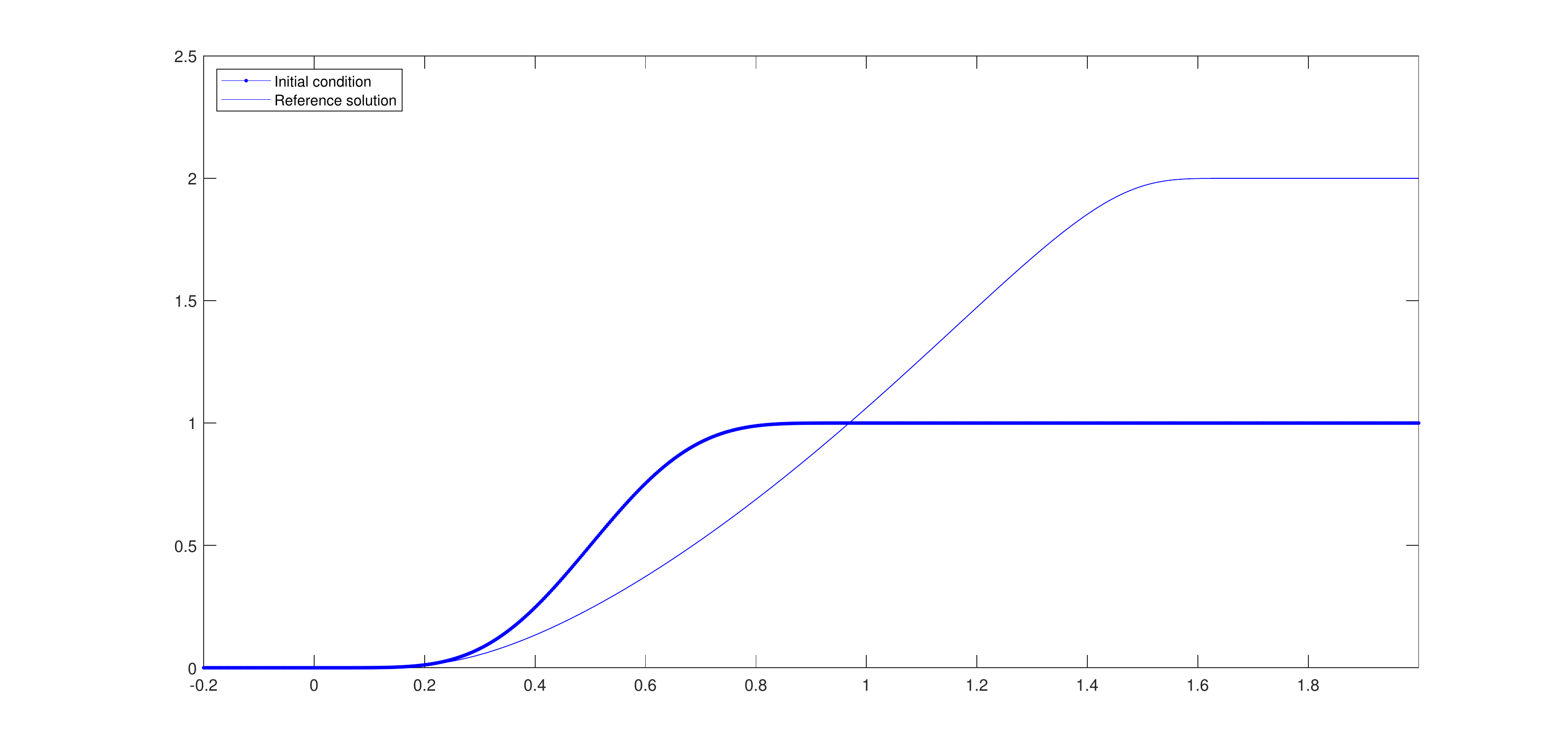}
    	\vspace{-0.8 cm}
    	\caption{Test \ref{Burg_5}. Initial condition and Reference solution obtained with WBACAT4 using a 2560 mesh points and CFL$=0.9$ at time $t=0.5.$}
    	\label{Test_Bu_3_ini}
    \end{figure}

    Let us consider \eqref{bur_scal} with $H(x) = x$ and initial condition \eqref{condsuave_1} (see Figure \ref{Test_Bu_3_ini}) 
    \begin{equation}
        \label{condsuave_2}
        u_0(x)=
        \begin{cases}
           0 \quad\quad\;\;\;\mathrm{if}\quad  x<0;\\
           p(x) \quad\;\;\mathrm{if}\quad 0\le x \le1; \\
           1 \quad\quad\;\;\;\mathrm{if} \quad x>1; 
        \end{cases}
    \end{equation}
    where  $$ p(x) = x^6\Bigl(\sum_{k=0}^{5}(-1)^k\binom{5+k}{k}(x-1)^k\Bigr).$$ 
    We solve the problem in the interval $[-0.2,2]$ using 80 mesh point and CFL$=0.9.$ As boundary condition free boundary is imposed at ghost points.
    
    \begin{figure}[!ht]
    	\centering	
    	\hspace{-1.2cm}
    	\includegraphics[scale = 0.5]{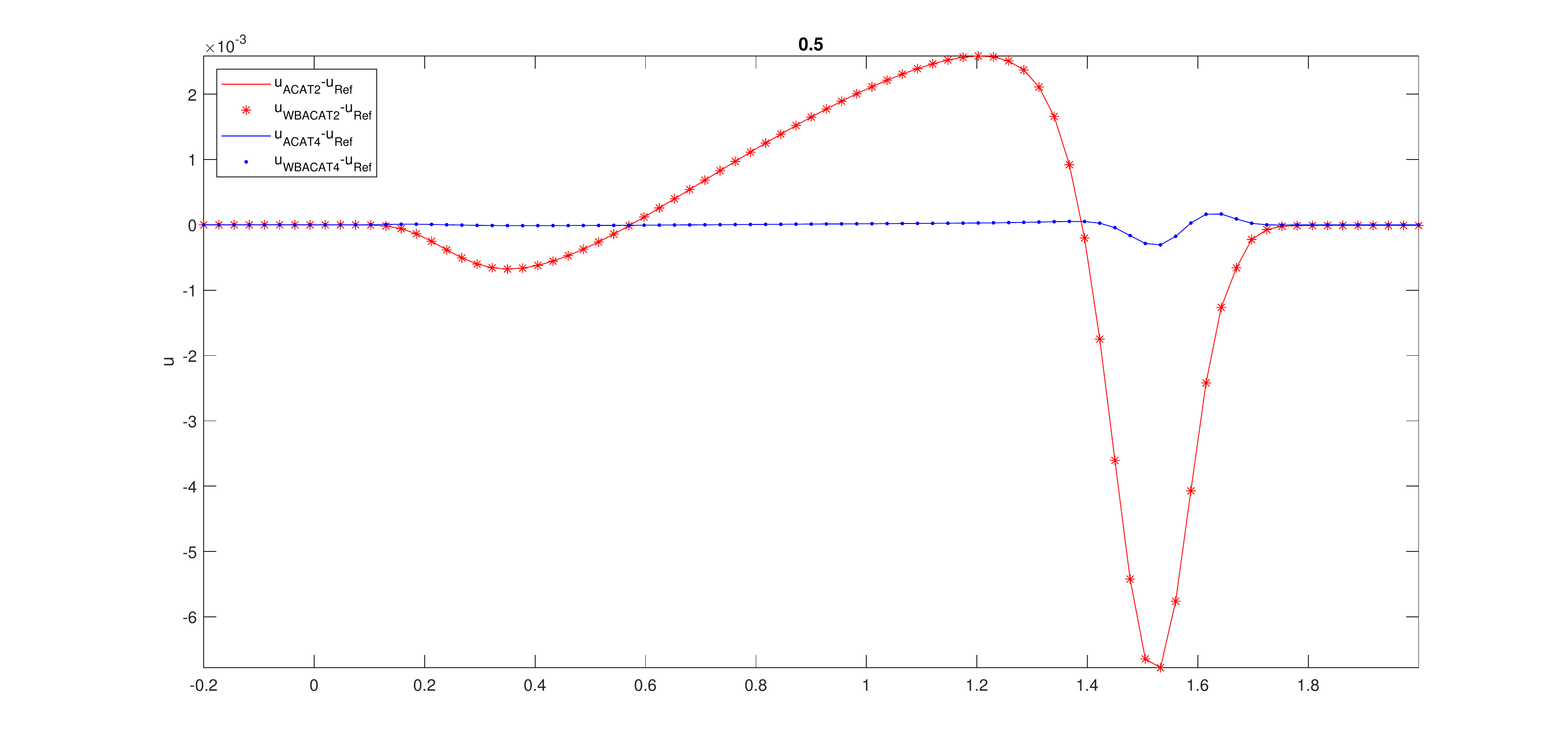}
    	\vspace{-0.8 cm}
    	\caption{Test \ref{Burg_5}. Differences between numerical solutions computed with ACAT$2P$ and WBCAT$2P$, $P=1,2$, and the reference solution at $t = 0.5$ using a mesh of 80 points and CFL$=0.9.$ For the reference solution a mesh of 2560 has been adopted.}
    	\label{Test_Bu_3}
    \end{figure}
    
    \begin{table}[htbp]
        \begin{center}
            \begin{tabular}{|c|c|c|c|c|c|c|c|c|}
                \hline                        &\multicolumn{2}{c|}{ACAT2}&\multicolumn{2}{c|}{WBACAT2}&               \multicolumn{2}{c|}{ACAT4} &  \multicolumn{2}{c|}{WBACAT4}  \\
                Points   &  Order & Error      & Order & Error        & Order  & Error     & Order& Error    \\ 
                \hline
                80       & -    &  1.74E-3     & -     & 1.91E-3      & -      & 2.41E-4   & -    &  2.42E-4 \\   
                \hline 
                160      & 1.93 &  4.11E-4     & 1.92  & 5.05E-4      & 2.40   & 4.56E-5   & 2.41 &  4.57E-5 \\   
                \hline  
                320      & 1.93 &  1.08E-4     & 1.94  & 1.32E-4      & 3.17   & 5.05E-6   & 3.18 &  5.06E-6 \\   
                \hline  
                640      & 1.98 &  2.74E-5     & 1.98  & 3.34E-5      & 3.68   & 3.95E-7   & 3.68 &  3.94E-7 \\   
                \hline  
                1280     & 1.99 &  6.91E-6     & 1.99  & 8.41E-6      & 3.91   & 2.60E-8   & 3.92 &  2.61E-8 \\   
                \hline
                2560     & 2.00 &  1.73E-6     & 2.00  & 2-11E-6      & 3.98   & 1.65E-9   & 3.98 &  1.65E-9 \\   
                \hline
            \end{tabular}
            \vspace{2mm}
            \caption{Test \ref{Burg_5}: Errors in $L^1$ norm and convergence rates for ACAT2, ACAT4, WBACAT2 and WBACAT4 at time $t=0.5$ and CFL$=0.9.$}
            \label{table_Bu_3}
            \end{center}
    \end{table}
    
    Figure \ref{Test_Bu_3} and Table \ref{table_Bu_3} show us how all the methods manage to produce solutions in agreement with each other obtaining the expected order. 
    
    This experiment shows that, when the initial condition is far from the stationary solution, well-balanced and non well-balanced methods produce essentially the same results, with the expected order of accuracy.

    \subsection{Shallow water model}
    In this section we will focus on the one-dimensional hyperbolic shallow water model
    \begin{equation}
    \label{Shallow_water}    
    \begin{cases}
        h_t + q_x =0\\
     \displaystyle   q_t + \left(\frac{q^2}{h} + \frac{g}{2}h^2 \right)_x = ghH_x,
    \end{cases}
    \end{equation}
that can be written in the form (\ref{bal_sis}) with
    \begin{equation*}
        U = 
        \begin{bmatrix}
        h \\ q
        \end{bmatrix},\quad
        F(U) =
        \begin{bmatrix}
        q \\ \displaystyle \frac{q^2}{h} + \frac{g}{2}h^2
        \end{bmatrix}, \quad
        S(U) = 
        \begin{bmatrix}
        0 \\ gh
        \end{bmatrix}.
    \end{equation*}
    The variable $x$ refers to the axis of the channel and $t$ is time; $q(x,t)$ and $h(x,t)$ represent the discharge and the water thickness; $g,$ the acceleration due to gravity; $H(x),$ the depth measured from a fixed level of reference; furthermore, the following relation is verified $q(x,t) =h(x,t)u(x,t),$ with $u$ the depth average horizontal velocity. The eigenvalues of the Jacobian matrix $J(U)$ of the flux $F(U)$ are
    $$ \lambda_1 = u-\sqrt{gh} \quad {\rm{and}} \quad \lambda_2 = u+\sqrt{gh}. $$
    The local Froude number is defined by
    $$
    Fr = \frac{|u|}{\sqrt{gh}}.
    $$
    The flow is said to be supercritical if $Fr > 1$ for all $x\in[a,b],$ critical if $Fr = 1$ for all $x\in[a,b]$ and subcritical if $Fr < 1$ for all $x\in[a,b].$
  
    The stationary solution of the shallow water system (\ref{Shallow_water}) are implicitly given by 
    \begin{equation}
        \label{SW_stat_sol}
        q = Q \quad{\rm{and}} \quad \ha\frac{q^2}{h^2} + gh -gH = C,
    \end{equation}
    where $Q$ and $C$ are arbitrary constants \cite{CLP2013}. In order to implement the well-balanced methods, given $U^n_i = [h^n_i, q^n_i]^T$ one has to find the stationary solution $U^*_i = [q^*_i, h^*_i(x)]^T$ that solves \eqref{Cauchyin}: it  is implicitly given by
$$
q^*_i(x) = q^n_i, \quad \frac{1}{2}\frac{(q^n_i)^2}{{h^*_i(x)}^2} + gh^*_i(x) - gH = C_i, \quad \forall x,
$$ 
with
$$
C_i = \frac{1}{2} \frac{(q^n_i)^2}{(h_i^n)^2} + gh^n_i - gH(x_i).
$$
Therefore, at a point $x_j$ of the stencil, one has $q_i^*(x_j) = q^n_i$ and $h_i^*(x_j)$ has to be a positive root of the polynomial:
$$
P_{i,j}(h) = h^3 - \left(\frac{C_i}{g} + H(x_j)\right) h ^2+ \frac{1}{2g}(q_i^n)^2.
$$
This polynomial can have two, one, or zero positive roots. In the first case, one of the roots corresponds to a supercritical state and the other one to a subcritical state: a criterion is necessary to select one root or the other. We follow here a criterion similar to the one chosen in \cite{CLP2013} in the context of finite volume methods: the solution whose regime (sub or supercritical) is the same as the one of $U^n_i$ is selected. A careful implementation is needed to capture transcritical stationary solutions: see for instance the discussion in \cite{CLP2013} or \cite{PP2021}. 
    
    \subsubsection{Preservation of a subcritical stationary solution} \label{SW_1}
\begin{figure}[!ht]
    	\centering	
    	\hspace{-1.2cm}
    	\includegraphics[scale = 0.5]{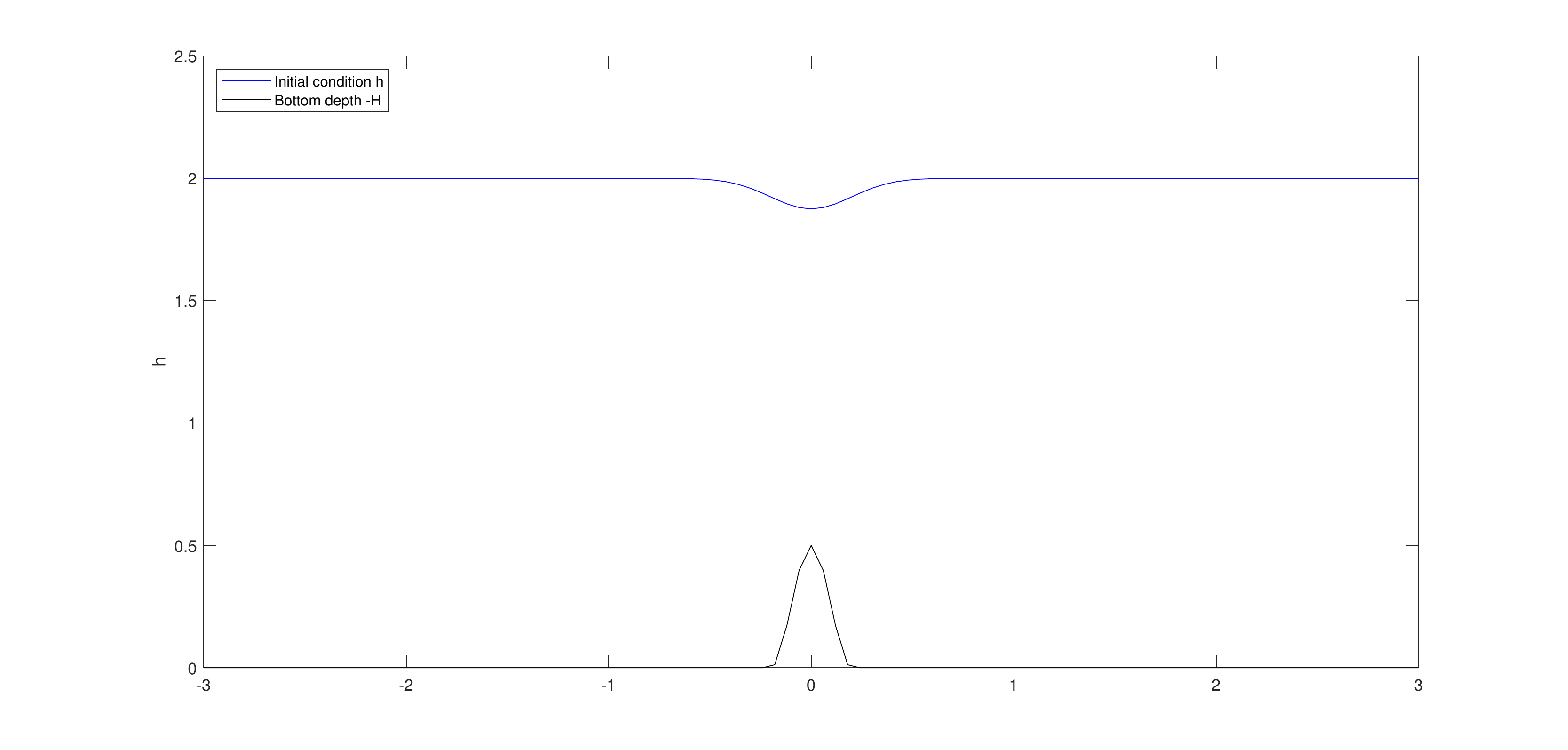}
    	\vspace{-0.8 cm}
    	\caption{Test \ref{SW_1}. Discrete initial condition with 100 mesh points. Free surface and bathymetry.}
    	\label{Test_3_1_1}
    \end{figure}

    Let us consider the shallow water model in the space interval $[-3,3]$ with  bottom depth given by
    \begin{equation}
        \label{H_SW_1}
        H(x) = 
        \begin{cases}
            -0.25(1+\cos(5\pi x)) \quad \rm{if} \; -0.2\le x\le 0.2;\\
            0 \quad \rm{otherwise};
        \end{cases}
    \end{equation}
    and initial condition given by the subcritical stationary solution $U^*$ that satisfies
    $$ q^* = 2.5,  \quad h^*(-3) = 2 $$ 
    (see Figure \ref{Test_3_1_1}). The numerical methods are applied to this problem using 200 mesh points and CFL$=0.8$. At the boundaries, the stationary solution is imposed at ghost points. 
    
    \begin{figure}[!ht]
    	\centering	
    	\hspace{-1.2cm}
    	\includegraphics[scale = 0.5]{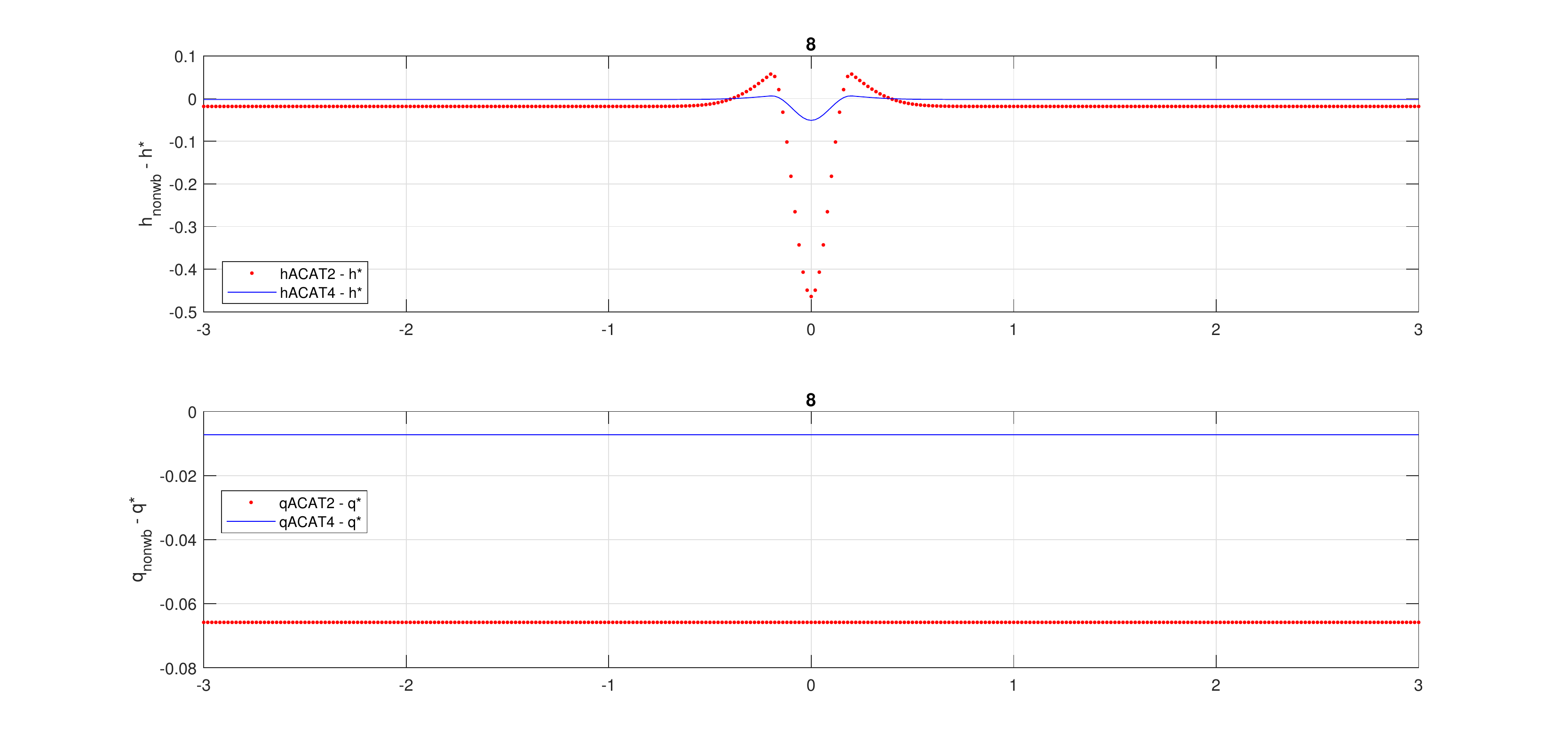}
    	\vspace{-0.8 cm}
    	\caption{Test \ref{SW_1}. Difference between the numerical solution for $h$ (top) and $q$ (bottom) obtained with  ACAT methods and the exact stationary one, at time $t=4$ using 200 mesh points and CFL$=0.8.$	}
    	\label{Test_3_1_2}
    \end{figure}
    Figures \ref{Test_3_1_2} shows the differences between the exact and the non well-balanced numerical solutions obtained at time $t=4.$ As it can be seen, the well-balanced methods capture the stationary solution to machine accuracy. This behaviour is confirmed by   Table \ref{table3_1_1} that shows the $L^1-$errors corresponding to WBACAT$2P$,  $P =1,2$, using  $50,$ $100,$ $200$ and $400$ mesh points at time $t = 4.$
    
    \begin{table}[htbp]
        \begin{center}
            \begin{tabular}{|c|c|c|c|c|c|c|}
                \hline   &\multicolumn{3}{c|}{WBACAT2}&  \multicolumn{3}{c|}{ WBACAT4}  \\
                Points   &  h & q        & u   & h          & q    & u   \\
                \hline
                50      & 2.93E-16       &  1.07E-16     & 2.66E-16     & 2.39E-16    & 5.32E-17   & 1.87E-16  \\ 
                \hline 
                100       &  3.46E-16     & 7.99E-17  & 1.86E-16     & 2.13E-16   & 0 &  1.20E-16\\   
                \hline  
                200      &  3.40E-16       & 0      & 2.46E-16     & 3.99E-17   & 0 &  1.99E-17 \\   
                \hline  
                400       &  1.77E-16     & 0      & 1.20E-16     & 0         & 5.99E-17 &  2.98E-17 \\   
                \hline  
            \end{tabular}
            \vspace{2mm}
            \caption{Test \ref{SW_1}. Errors in $L^1$ norm  for WBACAT$2P$,  $P =1,2$, at time $t=4.$}
            \label{table3_1_1}
            \end{center}
    \end{table}
     The introduction of spurious oscillations with the not well-balanced schemes involves a order reduction since the high order smoothness indicators are not able to detect a priori the real smoothness of the solution. This behaviour is highlighted in the next experiments.
    
    \subsubsection{Perturbation of a subcritical stationary solution} \label{SW_2}
\begin{figure}[!ht]
    	\centering	
    	\hspace{-1.2cm}
    	\includegraphics[scale = 0.5]{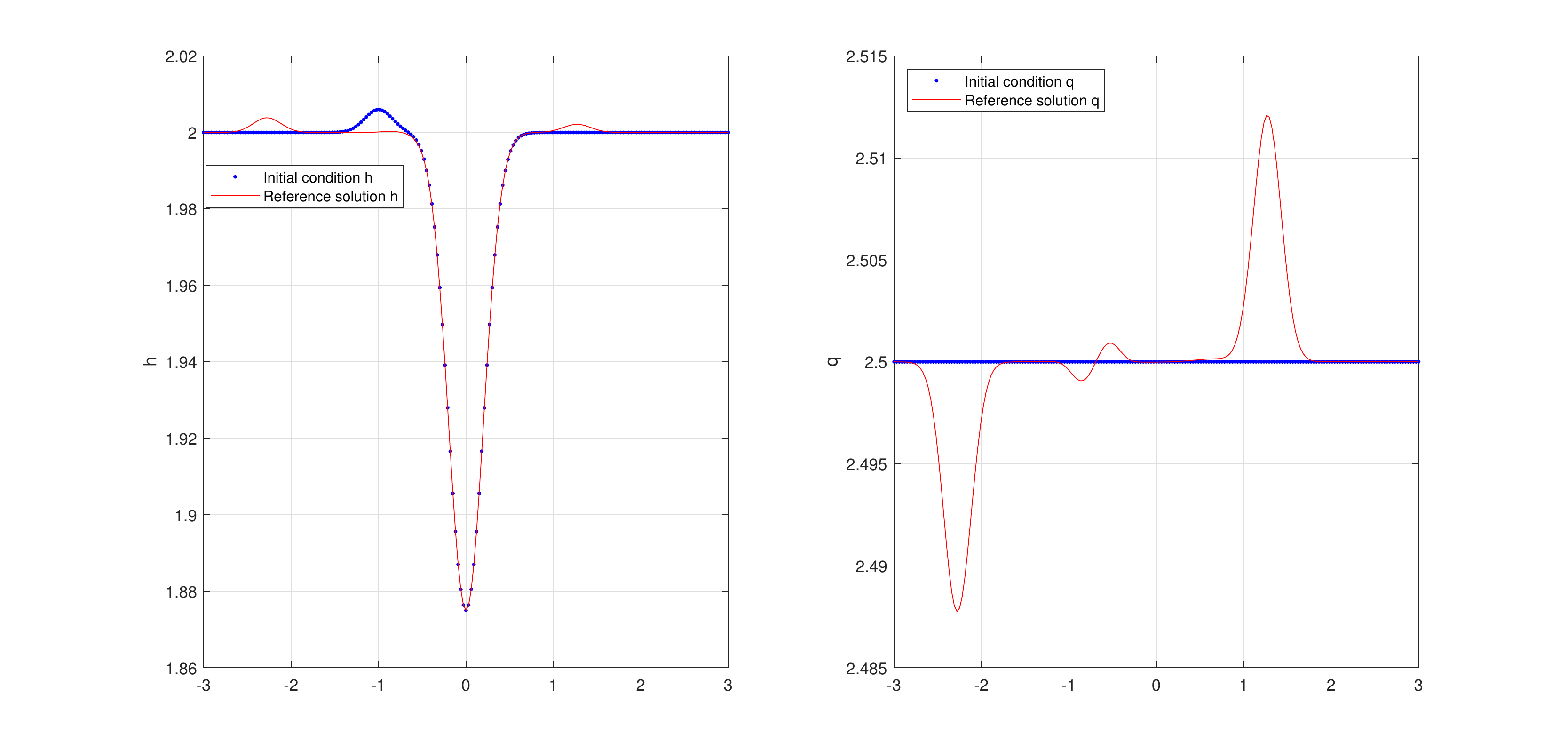}
    	\vspace{-0.8 cm}
    	\caption{Test \ref{SW_2}. Initial condition and reference solution obtained with WBACAT4 computed at time $t=0.4$ using $2000$ mesh points and CFL$=0.8:$ $h$ (left); $q$ (right). In the plot of $q$ there appear the left and right traveling waves, as well as a small left moving reflected wave.}
    	\label{Test_3_2_1}
    \end{figure}

    The setting of this test is similar to the previous one but now the initial condition is a smooth perturbation of the subcritical stationary solution $U^*$ (see Figure \ref{Test_3_2_1}) considered there:
    $$U_0=\begin{bmatrix}
    h^*+ 0.006e^{(-20(x+1)^2)}\\q^*
    \end{bmatrix}.$$ 
    
    \begin{figure}[!ht]
    	\centering	
    	\hspace{-1.2cm}
    	\includegraphics[scale = 0.5]{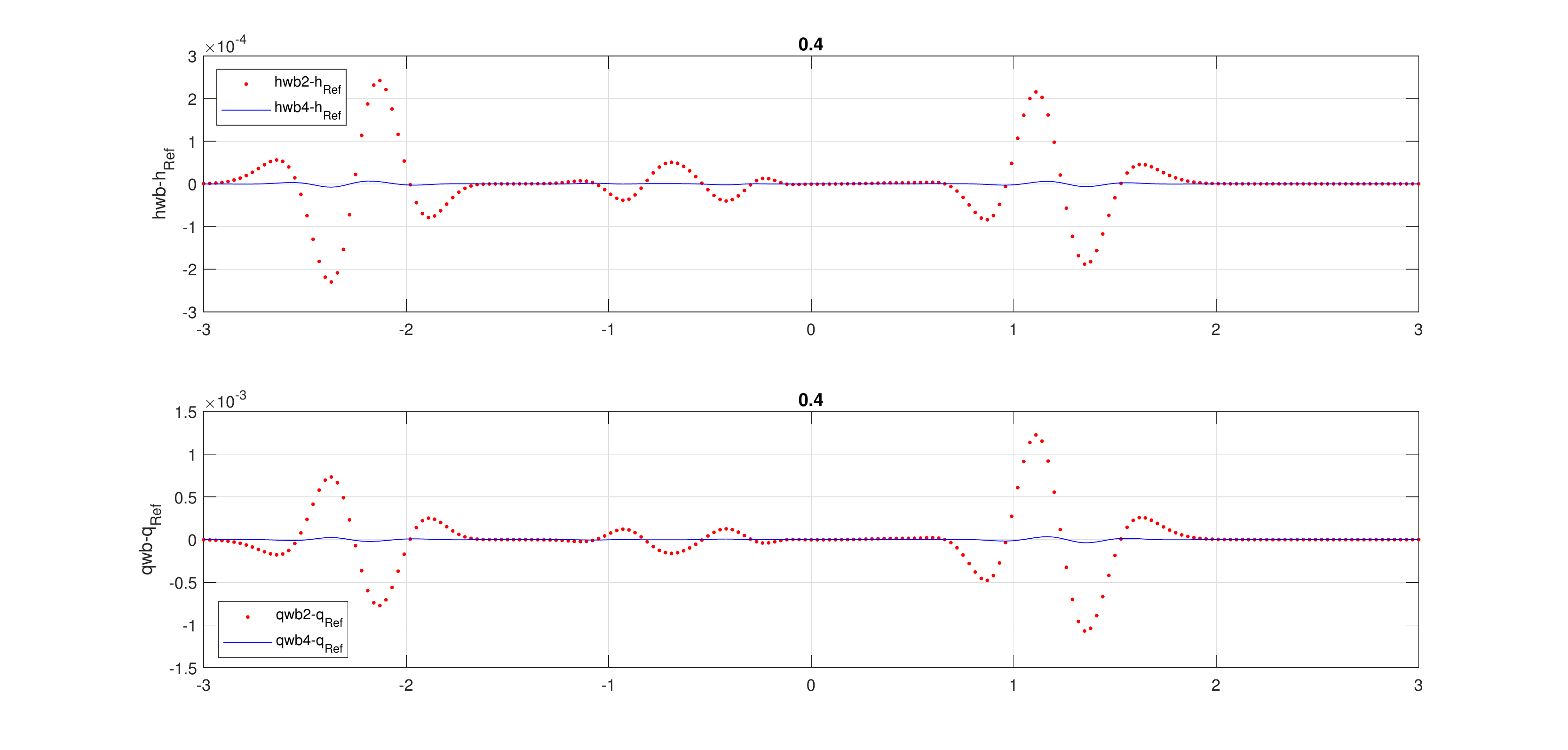}
    	\vspace{-0.8 cm}
    	\caption{Test \ref{SW_2}. Difference between Reference and numerical solutions obtained with WBACAT$2P$, $P = 1,2$, computed at time $t=0.4$ using $200$ mesh points and CFL$=0.8:$ $h$ (top); $q$ (bottom). The Reference solution is computed with WBACAT4 adopting a $2000$ mesh points.}
    	\label{Test_3_2_2}
    \end{figure}
    \begin{figure}[!ht]
    	\centering	
    	\hspace{-1.2cm}
    	\includegraphics[scale = 0.5]{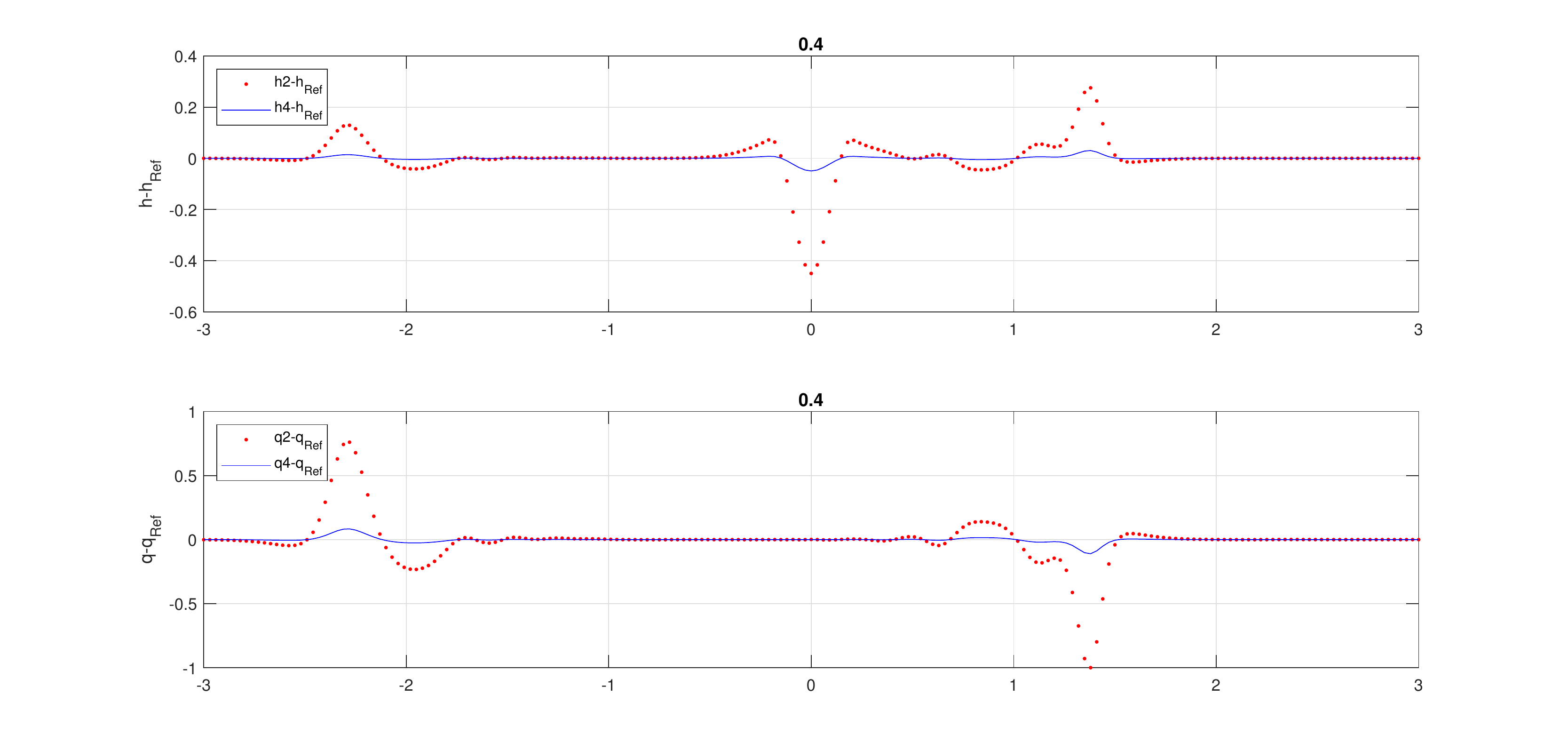}
    	\vspace{-0.8 cm}
    	\caption{Test \ref{SW_2}. Difference between Reference and numerical solutions obtained with ACAT$2P$, $P = 1,2$, computed at time $t=0.4$ using $200$ mesh points and CFL$=0.8:$ $h$ (top); $q$ (bottom). The Reference solution is computed with WBACAT4 adopting a $2000$ mesh points}
    	\label{Test_3_2_3}
    \end{figure}
    The numerical solutions are computed on the interval $[-3,3]$ using $200$ mesh points at time $t = 0.4$ with CFL $=0.8.$ As boundary conditions the subcritical stationary solution is imposed at ghost points.
   
    Figures \ref{Test_3_2_2} and \ref{Test_3_2_3} show the errors obtained by the differences between reference solution and the numerical solutions computed with well-balanced and not well-balanced methods at time $t=0.4.$ The Reference solution considered is WBACAT4 adopting a $2000$ mesh points. As expected, WBACAT$2P$,  $P = 1,2$, capture better the waves generated by the initial perturbation than ACAT$2P$, $P = 1,2$. 
    \begin{table}[htbp]
        \begin{center}
            \begin{tabular}{|c|c|c|c|c|c|c|c|c|}
                \hline                        &\multicolumn{2}{c|}{ACAT2}&\multicolumn{2}{c|}{WBACAT2}&               \multicolumn{2}{c|}{ACAT4} &  \multicolumn{2}{c|}{WBACAT4}  \\
                Points   &  Order & Error      & Order & Error        & Order  & Error     & Order& Error    \\ 
                \hline
                50       & -    &  2.18E-3     & -     & 1.64E-5      & -      & 2.07E-4   & -    &  1.81E-5 \\   
                \hline 
                100      & 0.98 &  1.10E-3     & 1.95  & 4.22E-6      & 1.11   & 9.61E-5   & 2.25 &  3.79E-6 \\   
                \hline  
                200      & 1.23 &  4.96E-4     & 1.94  & 1.09E-6      & 1.59   & 3.19E-5   & 4.72 &  1.44E-7 \\   
                \hline  
                400      & 1.48 &  1.68E-4     & 1.97  & 2.77E-7      & 1.84   & 8.87E-6   & 4.33 &  8.22E-9 \\   
                \hline  
                800      & 1.53 &  5.82E-5     & 1.97  & 7.07E-8      & 1.93   & 2.31E-6   & 4.07 &  6.54E-10\\   
                \hline
            \end{tabular}
            \vspace{2mm}
            \caption{Test \ref{SW_2}: Errors in $L^1$ norm and convergence rates related to $h$ for ACAT2, ACAT4, WBACAT2 and WBACAT4 at time $t=0.15$ and CFL$=0.8.$}
            \label{table3_2_1}
            \end{center}
    \end{table}  
    \begin{table}[htbp]
        \begin{center}
            \begin{tabular}{|c|c|c|c|c|c|c|c|c|}
                \hline                        &\multicolumn{2}{c|}{ACAT2}&\multicolumn{2}{c|}{WBACAT2}&               \multicolumn{2}{c|}{ACAT4} &  \multicolumn{2}{c|}{WBACAT4}  \\
                Points   &  Order & Error      & Order & Error        & Order  & Error     & Order& Error    \\ 
                \hline
                50       & -    &  8.21E-3     & -     & 5.54E-5      & -      & 7.85E-4   & -    &  6.15E-5 \\   
                \hline 
                100      & 1.00 &  4.11E-3     & 1.63  & 1.79E-5      & 1.22   & 3.35E-4   & 2.86 &  1.69E-5 \\   
                \hline  
                200      & 1.25 &  1.73E-3     & 1.94  & 4.68E-6      & 1.67   & 1.05E-4   & 4.78 &  6.17E-7 \\   
                \hline  
                400      & 1.51 &  6.06E-4     & 1.98  & 1.19E-6      & 1.87   & 2.86E-5   & 4.31 &  4.65E-8 \\   
                \hline  
                800      & 1.54 &  2.08E-4     & 1.98  & 3.01E-7      & 1.98   & 7.31E-6   & 4.05 &  1-41E-9\\   
                \hline
            \end{tabular}
            \vspace{2mm}
            \caption{Test \ref{SW_2}: Errors in $L^1$ norm and convergence rates related to $q$ for ACAT2, ACAT4, WBACAT2 and WBACAT4 at time $t=0.15$ and CFL$=0.8.$}
            \label{table3_2_2}
            \end{center}
    \end{table}
    In addition, Tables \ref{table3_2_1}-\ref{table3_2_2} show how well-balanced methods manage to reach the expected order, behavior not respected by non-well-balanced methods. In this case, as seen above, this phenomenon is partly attributable to not well-balanced reconstruction partly to smoothness indicators. In fact, the not well-balanced method at first step introduces a spurious error which implies a loss of numerical smoothness resulting in degradation of the order.

    \subsubsection{Smooth initial condition with flat bottom}\label{SW_3}
\begin{figure}[!ht]
    	\centering	
    	\hspace{-1.2cm}
    	\includegraphics[scale = 0.5]{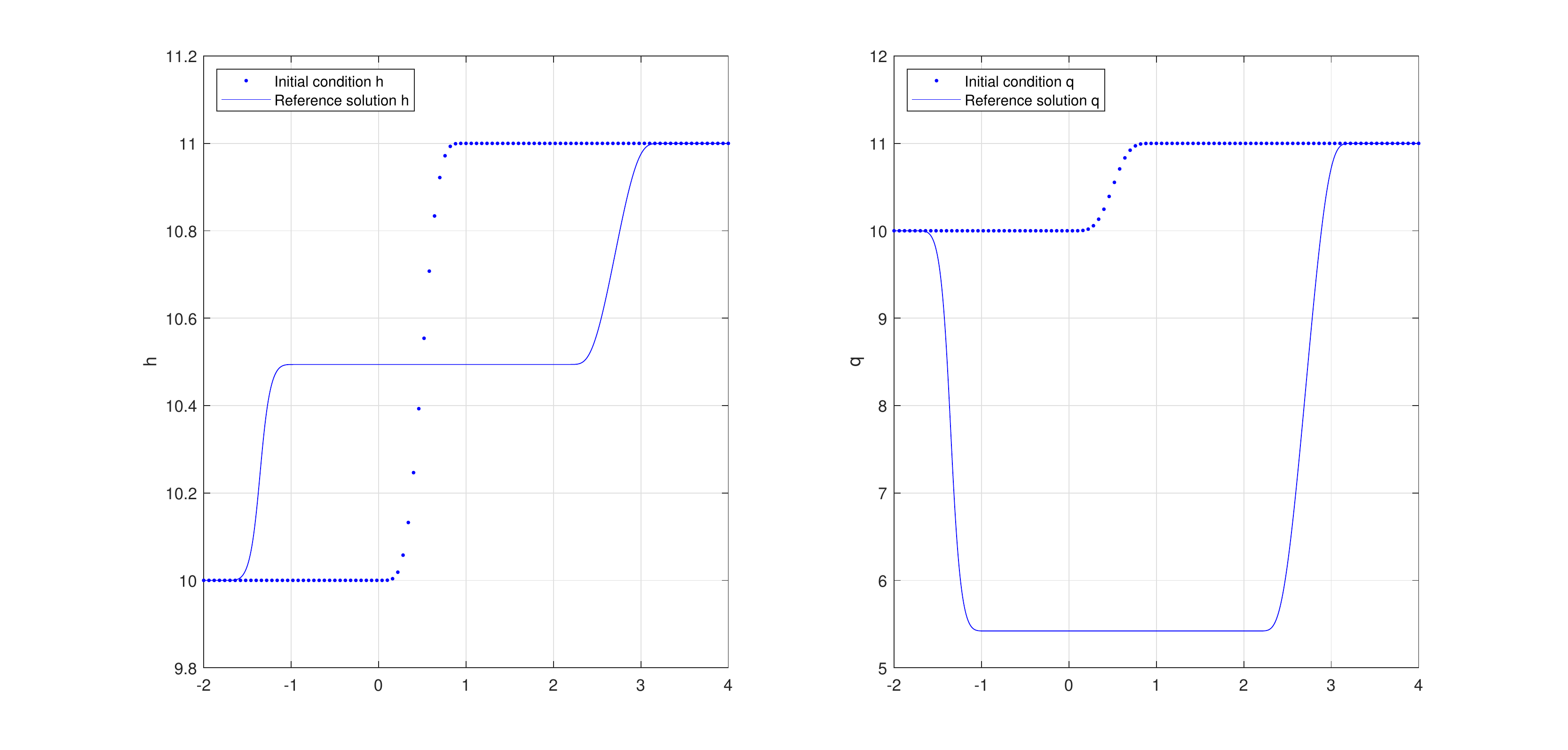}
    	\vspace{-0.8 cm}
    	\caption{Test \ref{SW_3}. Initial condition and reference solution obtained by WBACAT4 at time $t=0.2$ using $3200$ mesh points and CFL$=0.8;$ $h$ (top); $q$ (bottom).}
    	\label{Test_SW_smooth_in}
    \end{figure}

    We now check that in the case of flat bottom and smooth solution, well-balanced and non well-balanced schemes give the same result, all with the expected order of accuracy.
    In order to obtain these results we consider the Shallow water equation \eqref{Shallow_water} with flat bottom and smooth initial condition \eqref{condsuave_1} (see Figure \ref{Test_SW_smooth_in}):
    \begin{equation}
        U_0(x) = \begin{bmatrix}
        h_0(x) \\ q_0(x)
        \end{bmatrix},
    \end{equation}
    where 
    \begin{equation*}
        h_0(x) = q_0(x) =
            \begin{cases}
           0 \quad\quad \;\;\;\mathrm{if}\quad  x<0;\\
           p(x) \quad\;\;\mathrm{if}\quad 0\le x \le1; \\
           1 \quad \quad\;\;\;\mathrm{if} \quad x>1; 
        \end{cases}
    \end{equation*}
    and $$ p(x) = x^6\Bigl(\sum_{k=0}^{5}(-1)^k\binom{5+k}{k}(x-1)^k\Bigr).$$
    
    \begin{figure}[!ht]
    	\centering	
    	\hspace{-1.2cm}
    	\includegraphics[scale = 0.5]{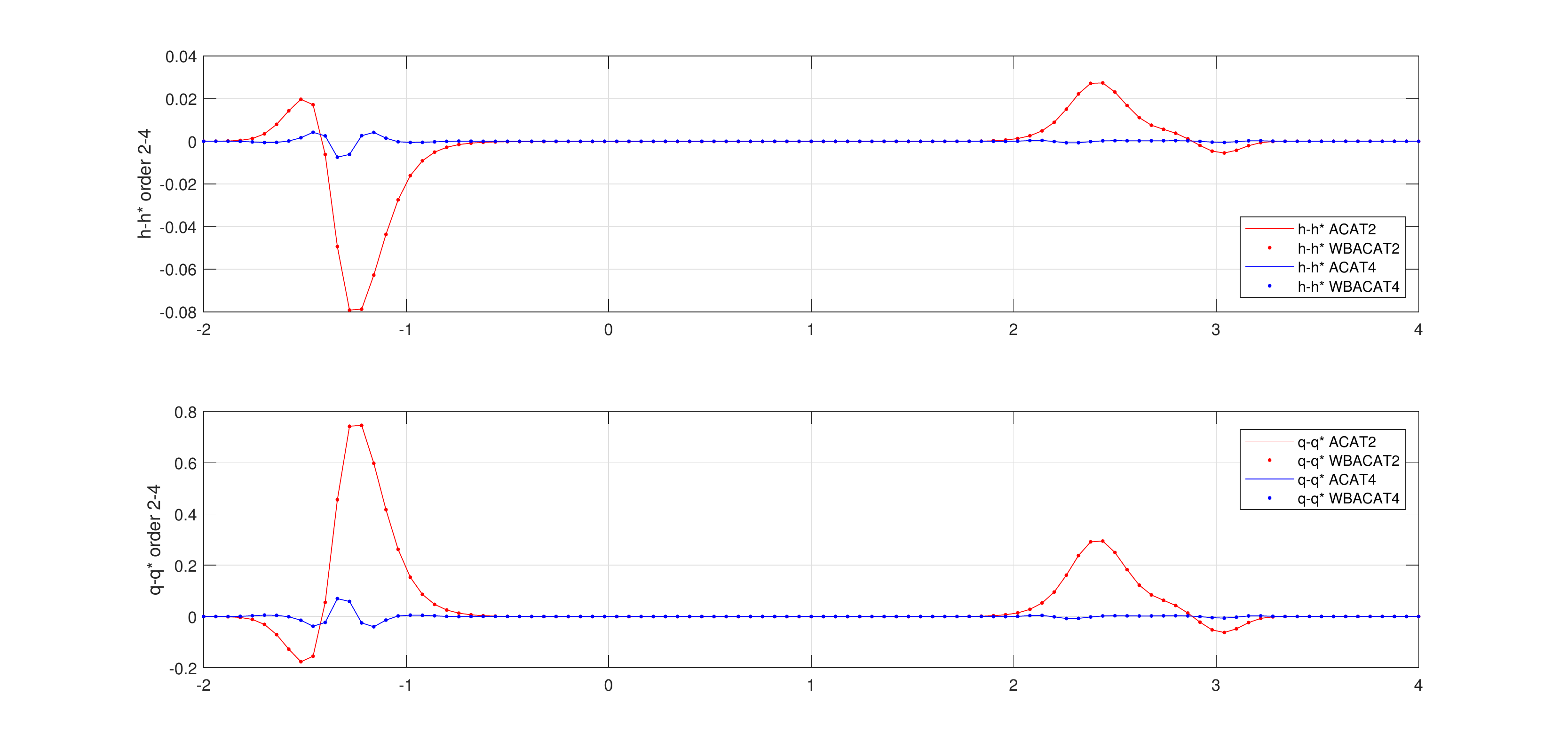}
    	\vspace{-0.8 cm}
    	\caption{Test \ref{SW_3}. Differences between numerical solutions obtained with ACAT$2P$ and WBACAT$2P$, $P = 1,2$, computed at time $t=0.2$ using $100$ mesh points and CFL$=0.8$ and the reference solution.  $h$ (top); $q$ (bottom). For the reference solution a 3200 mesh points has been adopted. }
    	\label{Test_SW_smooth}
    \end{figure}
    
    \begin{table}[htbp]
        \begin{center}
            \begin{tabular}{|c|c|c|c|c|c|c|c|c|}
                \hline                        &\multicolumn{2}{c|}{ACAT2}&\multicolumn{2}{c|}{WBACAT2}&               \multicolumn{2}{c|}{ACAT4} &  \multicolumn{2}{c|}{WBACAT4}  \\
                Points   &  Order & Error      & Order & Error        & Order  & Error     & Order& Error    \\ 
                \hline
                200      & -    &  7.27E-5     & -     & 7.25E-5      & -      & 2.25E-6   & -    &  2.25E-6 \\   
                \hline  
                400      & 2.01 &  1.80E-5     & 2.01  & 1.80E-5      & 3.98   & 1.42E-7   & 3.98 &  1.42E-7 \\   
                \hline  
                800      & 2.00 &  4.50E-6     & 2.00  & 4.50E-6      & 3.99   & 8.91E-9   & 3.99 &  8,91E-9 \\   
                \hline  
                1600     & 2.00 &  1.13E-6     & 2.00  & 1.13E-6      & 4.00   & 5.56E-10  & 4.00 &  5.57E-10\\   
                \hline
            \end{tabular}
            \vspace{2mm}
            \caption{Test \ref{SW_3}: Errors in $L^1$ norm and convergence rates related to $h$ for ACAT2, ACAT4, WBACAT2 and WBACAT4 at time $t=0.2$ and CFL$=0.9.$}
            \label{table3_3_1}
            \end{center}
    \end{table}
    \begin{table}[htbp]
        \begin{center}
            \begin{tabular}{|c|c|c|c|c|c|c|c|c|}
                \hline                        &\multicolumn{2}{c|}{ACAT2}&\multicolumn{2}{c|}{WBACAT2}&               \multicolumn{2}{c|}{ACAT4} &  \multicolumn{2}{c|}{WBACAT4}  \\
                Points   &  Order & Error      & Order & Error        & Order  & Error     & Order& Error    \\ 
                \hline
                200      & -    &  8.75E-4     & -     & 8.74E-4      & -      & 2.67E-5   & -    &  2.67E-5 \\   
                \hline  
                400      & 2.01 &  2.17E-4     & 2.01  & 2.17E-4      & 3.98   & 1.69E-6   & 3.98 &  1.69E-6 \\   
                \hline  
                800      & 2.00 &  5.42E-5     & 2.00  & 5.42E-5      & 3.99   & 1.06E-7   & 3.99 &  1.06E-7 \\   
                \hline  
                1600     & 2.00 &  1.35E-5     & 2.00  & 1.35E-5      & 4.00   & 6.61E-9   & 4.00 &  6.61E-9 \\   
                \hline
            \end{tabular}
            \vspace{2mm}
            \caption{Test \ref{SW_3}: Errors in $L^1$ norm and convergence rates related to $q$ for ACAT2, ACAT4, WBACAT2 and WBACAT4 at time $t=0.2$ and CFL$=0.9.$}
            \label{table3_3_2}
            \end{center}
    \end{table}
    The numerical solutions are computed on the interval $[-2,4]$ using $100$ mesh points and CFL $=0.9$ at time $t=0.1,$ while for the reference solution a $3200$ mesh points is adopted. As boundary conditions free boundary is imposed at ghost points.
    
    Figure \ref{Test_SW_smooth} and Table \ref{table3_3_1}-\ref{table3_3_2} show that all methods have a similar behavior by reproducing similar results. In particular, all schemes are accurate as expected. 

    With this experiment we have proven that all methods, both well-balanced and not, have a similar behavior when they are very far from the stationary condition; while, well-balanced reconstructions reproduce better results, both in accuracy and numerical convergence, when we are close to the stationary solution.

    \subsection{2D Euler system with gravity}\label{ss:2DEuler}
    As last test, let us consider the 2D system of compressible Euler equations with a gravity
    \begin{equation}
    \label{2D Euler with gravity}    
    \begin{cases}
        \rho_t + (\rho u)_x + (\rho v)_y=0,\\
        (\rho u)_t + (\rho u^2 + p)_x +(\rho uv)_y= -\rho H_x,\\
        (\rho v)_t +(\rho uv)_x + (\rho v^2 + p)_y= -\rho H_y,\\
        E_t + (u(E+p))_x + (v(E+p))_y = -\rho uH_x - \rho vH_y.
    \end{cases}
    \end{equation}
    Here, $\rho$ is the density; $u$, the velocity in $x-$direction; $v$, the velocity in $y-$direction; $p$, the pressure; $E$, the energy per unit volume excluding the gravitational energy; and $H(x,y)$, the gravitational potential \cite{Klin-Puppo-Semp}. The pressure is supposed to satisfy the equation of state
    $$ p = (\gamma - 1)\Bigl(E - \ha\rho(u^2+v^2)\Bigr), $$
    where $\gamma$ is the ratio between specific heats at constant pressure and volume, which is taken to be constant. System (\ref{2D Euler with gravity}) can be written in the form (\ref{2Dequ}) with
    \begin{equation*}
        U = 
        \begin{bmatrix}
        \rho \\ \rho u \\ \rho v \\ E
        \end{bmatrix},\quad
        F(U) =
        \begin{bmatrix}
        \rho u\\ \rho u^2 + p \\ \rho uv \\ u(E+p)
        \end{bmatrix}, \quad
        G(U) =
        \begin{bmatrix}
        \rho v\\ \rho uv \\\rho v^2 + p \\ v(E+p)
        \end{bmatrix}, \quad
        S(U) = 
        \begin{bmatrix}
        0 \\ -\rho H_x \\ -\rho H_y \\ -\rho uH_x - \rho vH_y
        \end{bmatrix}.
    \end{equation*}
    
    Hydrostatic stationary solutions satisfy 
    $$ u = 0, \quad v = 0, \quad \nabla p = -\rho\nabla H.$$ 
    A two-parameter family of isothermal stationary solution is given by
    \begin{equation}
        \label{2D_stat_sol_family}
        \rho^*(\mathbf{x}) = C_1e^{-H(\mathbf{x})}\ge 0; \quad p^*(\mathbf{x}) = C_2\rho^*(\mathbf{x})\ge 0; \quad u^* = v^*=0; \quad E^*(x)= \frac{p^*(x)}{\gamma-1}.
    \end{equation}
    Given  
    $$ U_{\bi}^n = [\rho_{\bi}^n, \rho_{\bi}^n u_{\bi}^n, \rho_{\bi}^n v_{\bi}^n, E_{\bi}^n]^T, $$
    the stationary solution $U^*_{\bi}$ selected by applying the technique described in Subsection \ref{ss_family} is then
    \begin{equation}
        \label{2D_stat_sol}
        \rho_{\bi}^*(\mathbf{x}) = \rho_{\bi}^ne^{-(H(\mathbf{x})- H(\mathbf{x}_{\bi}))}; 
        \;p^*_{\bi}(\mathbf{x}) = \rho_{\bi}^ne^{-(H(\mathbf{x}) - H(\mathbf{x}_{\bi}))}; \; u_{\bi}^* = v_{\bi}^* =0; \; E_{\bi}^*(x) = \frac{p^*_{\bi}(x)}{\gamma-1}.
    \end{equation}
    
    \subsection{Preservation of a continuous stationary solution}\label{2D_Eg_1}
    Following \cite{Kappeli-Mishra}-\cite{Klin-Puppo-Semp}-\cite{Grosheintz-Kappeli2019}-\cite{Grosheintz-Kappeli2020} we consider Euler equations in the 2D domain $[0,1]\times[0,1]$ with two different  gravitational potentials
    $$
    H_1(x,y) = x+y, \quad H_2(x,y) = \dfrac{1}{\sqrt{(x-\frac{1}{3})^2 + (y+\ha)^2}}
    $$
    and  initial condition
    \begin{align}
        \label{2D_Eg_1_init_cond}
        \rho(\mathbf{x},0) = e^{-H(\mathbf{x})}; \quad p(\mathbf{x},0) = e^{-H(\mathbf{x})}; \quad u(\mathbf{x},0) = v(\mathbf{x},0) =0.
    \end{align}
    
    \begin{table}[htbp]
        \begin{center}
            \begin{tabular}{|c|c|c|c|c|c|c|}
                \hline   &\multicolumn{6}{c|}{2D density}                    \\
                         & \multicolumn{2}{c|}{2D ACAT2}  & \multicolumn{2}{c|}{2D ACAT4}         & 2D WBACAT2        & 2D WBACAT4 \\
                Points  & Error & Order & Error & Order & Error & Error\\
                \hline
                20$\times$20  &  4.87E-6 & -    &  7.85E-9  & -     & 2.72E-17   & 2.96E-18  \\ 
                \hline 
                40$\times$40  &  1.91E-6 & 1.35 &  1.01E-9  & 2.95  & 2.19E-17   & 2.39E-18  \\  
                \hline  
                80$\times$80  &  5.62E-7 & 1.76 &  8.54E-11 & 3.56  & 1.82E-17   & 2.39E-18  \\  
                \hline  
                160$\times$160&  1.43E-7 & 1.98 &  5.61E-12 & 3.93  & 2.37E-18   & 2.64E-18  \\  
                \hline
            \end{tabular}
            \vspace{2mm}
            \caption{Test \ref{2D_Eg_1}: 2D Euler equation with gravity and gravitational potential $H_1$. Errors in $L^1$ norm for density at time $t=0.3.$ }
            \label{table2D_ACAT2_1}
            \end{center}
    \end{table}
    
    \begin{table}[htbp]
        \begin{center}
            \begin{tabular}{|c|c|c|c|c|c|c|}
                \hline   &\multicolumn{6}{c|}{2D density}                    \\
                         & \multicolumn{2}{c|}{2D ACAT2}  & \multicolumn{2}{c|}{2D ACAT4}         & 2D WBACAT2        & 2D WBACAT4 \\
                Points  & Error & Order & Error & Order & Error & Error\\
                \hline
                20$\times$20  &  3.85E-5 & -    &  3.87E-5 & -     & 2.50E-17   & 3.77E-17  \\ 
                \hline 
                40$\times$40  &  1.58E-5 & 1.28 &  5.16E-6 & 2.91  & 3.23E-17   & 3.59E-17  \\  
                \hline  
                80$\times$80  &  4.78E-6 & 1.72 &  4.45E-7 & 3.53  & 3.33E-17   & 3.33E-17  \\  
                \hline  
                160$\times$160&  1.23E-6 & 1.96 &  2.89E-8 & 3.95  & 3.15E-17   & 3.23E-17  \\  
                \hline
            \end{tabular}
            \vspace{2mm}
            \caption{Test \ref{2D_Eg_1}: 2D Euler equation with gravity and gravitational potential $H_2$. Errors in $L^1$ norm for density at time $t=0.3.$ }
            \label{table2D_ACAT2_2}
            \end{center}
    \end{table}
    
    We solve numerically the equations using a $(21\times21)-$point mesh and CFL$=0.9$. As boundary condition the exact solution is imposed to all sides through the ghost points. 
     Tables \ref{table2D_ACAT2_1} and  \ref{table2D_ACAT2_2} exhibit the errors in $L^1-$norm for  ACAT$2P$, WBACAT$2P$, $P = 1,2$ corresponding to $H = H_1$ and $H = H_2$ respectively. As it can be seen, the differences between the solutions given by well-balanced and no well-balanced methods are bigger for $H = H_2$: please note that, in the case  of the linear potential $H_1,$ the stationary solution is essentially 1D while this is not true in the case $H = H_2.$

    \subsection{Perturbation of the stationary solution}\label{2D_Eg_2}
    We consider now Euler equations in the 2D domain $[0,1]\times [0,1]$ with the gravitational potential $ H_2$ 
    and an initial condition that represents a perturbation of the hydrostatic stationary considered in the previous test case:
    \begin{align}
        \label{2D_Eg_2_init_cond}
        \rho(\mathbf{x},0) &= e^{-H(\mathbf{x})} + 0.008e^{-200(x-0.5)^2 - 200(y-0.5)^2}; \quad p(\mathbf{x},0) = e^{-H(\mathbf{x})}+ 0.008e^{-200(x-0.5)^2 - 200(y-0.5)^2};\nonumber\\ u(\mathbf{x},0) &= v(\mathbf{x},0) =0.
    \end{align}
    
    \begin{table}[htbp]
        \begin{center}
            \begin{tabular}{|c|c|c|c|c|c|c|c|c|}
                \hline   &\multicolumn{8}{c|}{2D density}                    \\
                         & \multicolumn{2}{c|}{2D ACAT2}  & \multicolumn{2}{c|}{2D ACAT4}         &   \multicolumn{2}{c|}{2D WBACAT2}       &  \multicolumn{2}{c|}{2D WBACAT4}\\
                Points  & Error & Order & Error & Order & Error & Order & Error & Order\\
                \hline
                20$\times$20  &  4.49E-5 & -    &  5.93E-6 & -     & 8.27E-6  & -     & 4.28E-7 & -     \\ 
                \hline 
                40$\times$40  &  2.47E-5 & 0.86 &  1.71E-6 & 1.79  & 4.41E-6  & 0.91  & 4.92E-8 & 2.71  \\ 
                \hline  
                80$\times$80  &  1.21E-5 & 1.03 &  4.37E-7 & 1.97  & 2.08E-6  & 1.08  & 7.05E-9 & 2.80  \\ 
                \hline  
                160$\times$160&  5.45E-6 & 1.15 &  9.85E-8 & 2.15  & 8.13E-7  & 1.36  & 9.47E-10 & 2.90  \\ 
                \hline
                320$\times$320&  2.43E-6 & 1.17 &  2.13E-8 & 2.21  & 2.45E-7  & 1.73  & 1.22E-10& 2.96  \\ 
                \hline
            \end{tabular}
            \vspace{2mm}
            \caption{Test \ref{2D_Eg_2}: 2D Euler equation with gravity and gravitational potential $H_2$. Errors in $L^1$ norm for density at time $t=0.2.$ }
            \label{table2D_ACAT2_3}
            \end{center}
    \end{table}
    
   Table \ref{table2D_ACAT2_3} shows error in $L^1-$norm and convergence rates for the numerical solutions obtained with ACAT$2P$ and WBACAT$2P$ and the reference solution at time $t=0.2,$ with $P=1,2.$ As happened for Shallow water, in case that a small perturbation of the stationary solution is considered as initial condition, the well-balanced schemes manage to capture the solution with a better accuracy than standard methods. This phenomena is shown on Tables \ref{table2D_ACAT2_3}. 
    
\subsection{Acoustic propagation}\label{Test_low_per}
    As last experiment we consider the Euler equations in the 2D domain $[0,2]\times [0,2]$ with the gravitational potential $ H_3,$
    $$ H_3(x,y) = \dfrac{1}{\sqrt{(x-0.4)^2 + (y+0.1)^2}}, $$ 
    and an initial condition that represents a very small perturbation of the hydrostatic stationary considered in the previous test case:
    \begin{align}
        \label{2D_Eg_2_init_cond_low}
        \rho(\mathbf{x},0) = e^{-H(\mathbf{x})} + 0.000001e^{-200(x-1)^2 - 200(y-1)^2}; \quad p(\mathbf{x},0) = \rho(\mathbf{x},0); \quad u(\mathbf{x},0) = v(\mathbf{x},0) =0.
    \end{align}
    
    \begin{figure}[!ht]
    	\centering	
    	\hspace{-1.2cm}
    	\includegraphics[scale = 0.5]{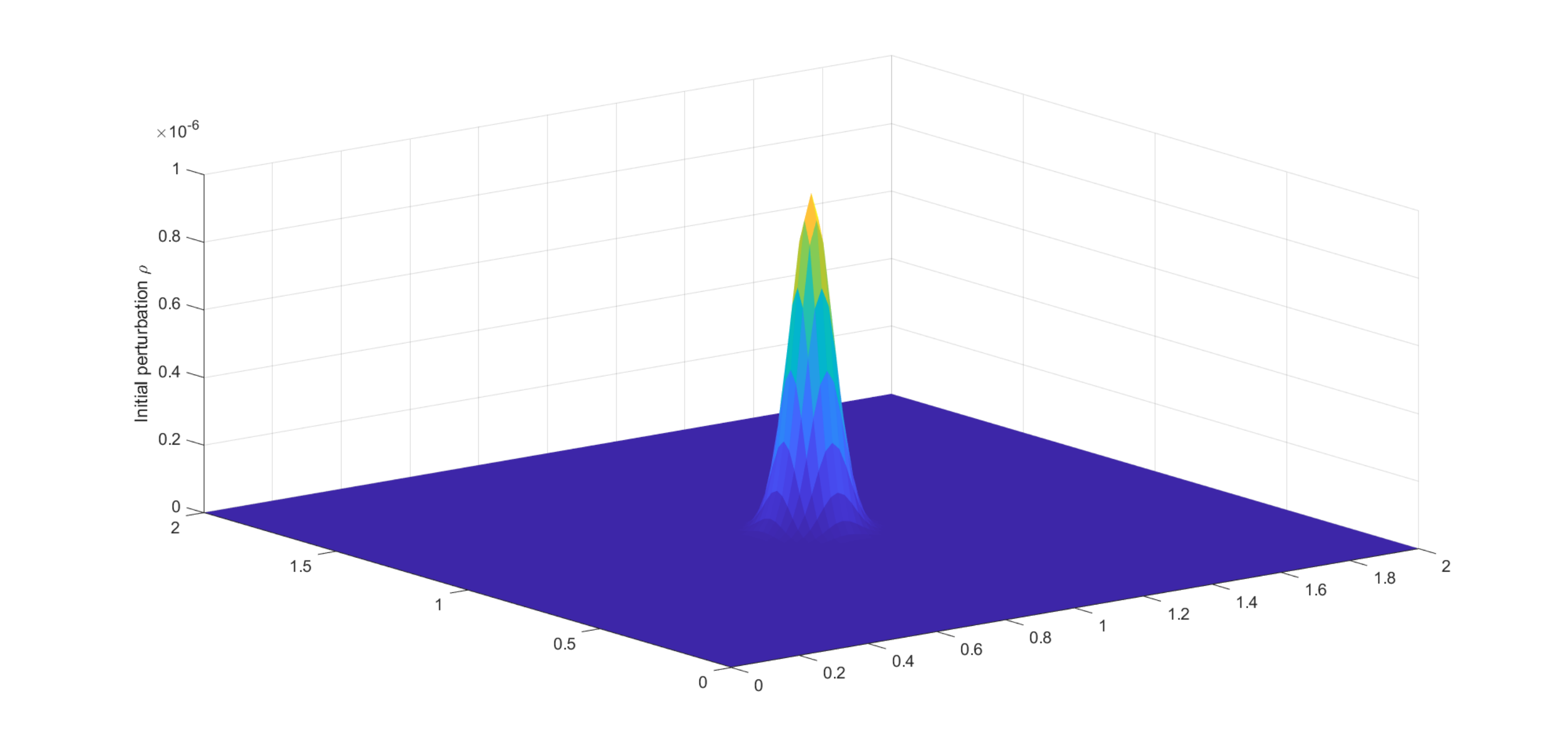}
    	\vspace{-0.8 cm}
    	\caption{Test \ref{Test_low_per}: Euler equations with gravitational potential $H_3$. Initial perturbation using a $101\times101$ mesh points.}
    	\label{Test_5_2_1}
    \end{figure}
    \begin{figure}[!ht]
    	\centering	
    	\hspace{-1.2cm}
    	\includegraphics[scale = 0.5]{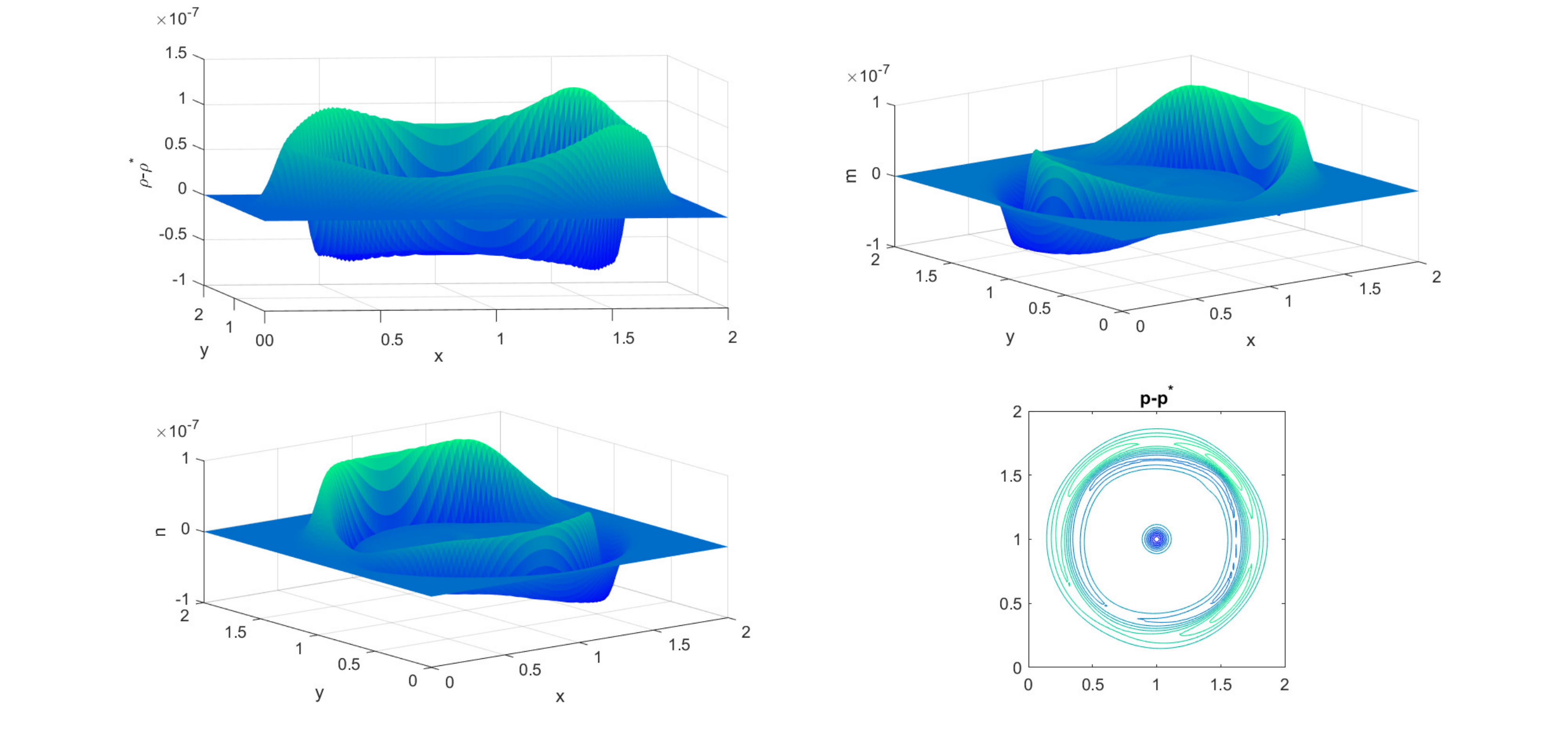}
    	\vspace{-0.8 cm}
    	\caption{Test \ref{Test_low_per}: Euler equations with gravitational potential $H_3$. Difference between the numerical solutions and the stationary solution computed at  time $t = 0.75$ with  WBACAT2  using $101\times101$ mesh points and CFL$=0.8.$  }
    	\label{Test_5_1_1}
    \end{figure}
    \begin{figure}[!ht]
    	\centering	
    	\hspace{-1.2cm}
    	\includegraphics[scale = 0.5]{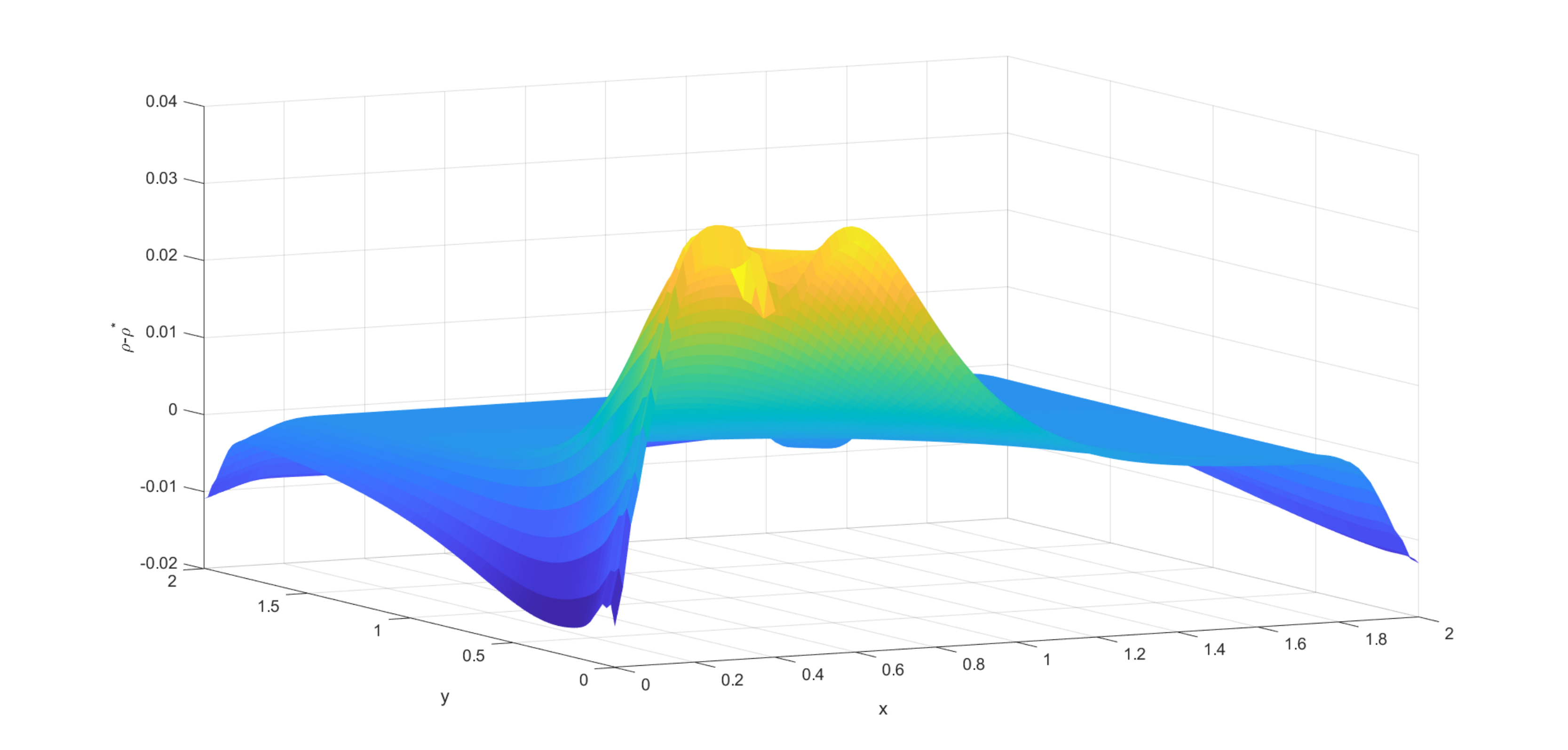}
    	\vspace{-0.8 cm}
    	\caption{Test \ref{Test_low_per}: Euler equations with gravitational potential $H_3$. Difference between the density solutions and the stationary solution computed at  time $t = 0.75$ with  ACAT2  using $101\times101$ mesh points and CFL$=0.8.$  }
    	\label{Test_5_1_2}
    \end{figure}
    
    Figures \ref{Test_5_1_1} shows the difference between the numerical solutions and the stationary solution computed at  time $t = 0.75$ with  WBACAT2 a using $101\times101$ mesh points and CFL$=0.8.$ As expected, the singularity of the gravitational potential modifies the thickness of the corona relative to the signal propagation, thinning it close to the singularity. $H_3$ has a singularity on position $(0.4,-0.1).$ \\ As we can see in Figure \ref{Test_5_1_2}, the non well-balanced method ACAT2 is not able to capture the evolution of the wave generated by the initial perturbation, since the numerical errors are much bigger than the wave amplitude: one would need a space step at least two orders of magnitude lower in order to have a truncation error of the same order of the signal, making computation with non well-balanced method absolutely impractical.

\section{Conclusion}
This paper deals with the construction, analysis, implementation and testing of well-balanced \\ ACAT2$P$ methods to the treatment of hyperbolic systems of balance laws. The starting point is to re-write the systems as conservation laws, by subtracting to the flux a primitive of the the source term. Meanwhile, the well-balanced property has been obtained rewriting the systems as conservation laws, by subtracting to the flux the flux of the stationary solution, and adding to the primitive of the source the source computed at stationary solution.
The methods are developed for systems in one and two space dimensions, and can be extended to 3D.
In principle the procedure allows the construction of well-balanced schemes of arbitrary order, although the computational complexity quickly increases with the order of accuracy.  
We prove that the constructed schemes are exactly well-balanced.

The use of suitable limiters allow an effective treatment of discontinuous solutions. Several test cases have been performed for the scalar equation and for systems in one and two space dimensions. In all cases we observe that stationary solutions are preserved within machine precision, allowing very accurate results when the solution is a small deviation from equilibrium. 

The main advantage of the method consists in its generality: it allows the automatic construction of very high order well-balanced schemes. 

There are still a few things that require improvement and generalization. 
First, a careful complexity analysis of the methods is needed, together with an improvement of the computational efficiency, possibly by parallel computing techniques, which should be possible thanks to the local computations required by these methods. 

Second, we shall explore new limiting strategies, which will allow more accurate results when the solution is smooth.

An open problem is how to couple ACAT methodology with IMEX method for the treatment of problems with stiff source. 

All these issues are  subject of current investigation. 

\section*{Acknowledgements} 
 This research has received funding from the European Union’s Horizon 2020 research and innovation program, under the Marie Sklodowska-Curie grant agreement No 642768. E. Macca was partially supported by GNCS Research Project "Approssimazione numerica di problemi di natura iperbolica ed applicazioni". E. Macca and G.Russo would like to thank the Italian Ministry of Instruction, University and Research (MIUR) to support this research with funds coming from PRIN Project 2017 (No. 2017KKJP4X entitled “Innovative numerical methods for evolutionary partial differential equations and applications”). The  research  of C. Parés  was  partially  supported  by  the  Spanish  Government(SG),  the  European  Regional  Development Fund(ERDF), the Regional Government of Andalusia(RGA), and the University of M\'alaga(UMA) through the projects  of  reference  RTI2018-096064-B-C21  (SG-ERDF),  UMA18-Federja-161  (RGA-ERDF-UMA), and P18-RT-3163 (RGA-ERDF). E. Macca and G. Russo are members of the INdAM Research group GNCS.
 \newpage

\bibliographystyle{acm}
\bibliography{biblio}
	
\end{document}